\documentclass[a4paper,12pt]{article}

\usepackage {amssymb,latexsym}
\usepackage {amsmath}

\def\R{\mathbb{R}}

\usepackage{xypic}
\xyoption{all}
\input{xypic}
\newdir{ >}{{}*!/-8pt/\dir{>}}

\newtheorem{theorem}{Theorem}[section]
\newtheorem{prop}{Proposition}[section]
\newtheorem{lemma}{Lemma}[section]

\newtheorem{defi}{Definition}[section]

\newtheorem{Exem}{{\bf Example}}[section]
\newenvironment{exem}{\begin{Exem} \strut\\ \normalfont}{\hfill$\Diamond$\end{Exem}}
                                                     \medskip

\newtheorem{Rem}{{\bf Remark}}[section]
\newenvironment{rem}{\begin{Rem} \strut\\ \normalfont}{\hfill$\Diamond$\end{Rem}}
                                                     \medskip
\newenvironment{prooof}{
                        \noindent{\bf\small Proof: }\small}
                                       {\hfill {$\mathbf \Box$}\\}

\makeatletter
\@addtoreset{equation}{section}

\makeatother

\begin{document}

\thispagestyle{empty}
\pagestyle{empty}
\title
   {\LARGE \bf Structure theory of Rack-Bialgebras }


\author{Charles Alexandre \\
Universit\'{e} de Strasbourg
  \thanks{Char.Alexandre@gmail.com}\\
\and Martin Bordemann \\
Universit\'{e} de Mulhouse
  \thanks{Martin.Bordemann@uha.fr}\\
\and Salim Rivi\`ere\\
Universit\'{e} de Nantes
\thanks{salim.riviere@hotmail.fr} \\
\and Friedrich Wagemann\\
     Universit\'{e} de Nantes
     \thanks{wagemann@math.univ-nantes.fr}}

\maketitle
\thispagestyle{empty}

\vspace{1cm}
\footnotesize
\vspace{-2cm}
\begin{abstract}
 In this paper we focus on a certain self-distributive multiplication on
 coalgebras, which leads to so-called rack bialgebra. Inspired by semi-group 
 theory (adapting the Suschkewitsch theorem), we do some structure theory for 
 rack bialgebras and cocommutative Hopf dialgebras. We also construct canonical 
 rack bialgebras (some kind of enveloping algebras) for any Leibniz algebra
 and compare to the existing constructions.  
 
 We are motivated
 by a differential geometric procedure which we call the Serre functor:
 To a pointed differentible manifold with multiplication is associated
 its distribution space supported in the chosen point. For Lie groups,
 it is well-known that this leads to the universal enveloping algebra
 of the Lie algebra. For Lie racks, we get rack-bialgebras, for Lie digroups,
 we obtain cocommutative Hopf dialgebras.  
\end{abstract}

\normalsize

\tableofcontents
 \thispagestyle{empty}

\pagestyle{plain}

\section{Introduction}

All manifolds considered in this manuscript are assumed
to be Hausdorff and second countable.

Basic Lie theory relies heavily on the fundamental links between 
associative algebras, Lie algebras and groups. Some of these links are the 
passage from an associative algebra $A$ to its underlying Lie algebra $A^{\rm Lie}$
which is the vector space $A$ with the bracket $[a,b]:=ab-ba$. On the other hand, 
to any Lie algebra ${\mathfrak g}$ one may associate its universal enveloping algebra
$\mathsf{U}({\mathfrak g})$ which is associative. Groups arise as groups of units in 
associative algebras. To any group $G$, one may associate its group algebra $KG$
which is associative. The theme of the present article is to investigate links 
of this kind for more general objects than groups, namely for racks and digroups.   

Recall that a \emph{pointed rack} is a pointed set $(X,e)$ together with a 
binary operation 
$\rhd:X\times X\to X$ such that for all $x\in X$, the map $y\mapsto x\rhd y$ 
is bijective and such that for all $x,y,z\in X$,
the self-distributivity and unit relations
$$
   x\rhd(y\rhd z)\,=\,(x\rhd y)\rhd(x\rhd z),~~~e\triangleright x = x,~~~
   \mathrm{and}~~~x\triangleright e = e
$$
are satisfied. 
Imitating the notion of a Lie group, the smooth version of a pointed rack is called
{\it Lie rack}.  

An important class of examples of racks are the so-called {\it augmented racks}, 
see \cite{FR92}. 
An augmented rack is the data of a group $G$, a $G$-set $X$ and a map $p:X\to G$ 
such that for all $x\in X$ and all $g\in G$, 
$$p(g\cdot x)\,=\,gp(x)g^{-1}.$$
The set $X$ becomes then a rack by setting $x\rhd y\,:=\,p(x)\cdot y$. 
In fact, augmented racks are the Drinfeld center (or the Yetter-Drinfeld modules)
in the monoidal category of $G$-sets over the (set-theoretical) Hopf algebra $G$, see for example 
\cite{KW14}. Any rack may be augmented in many ways, for example by using the canonical 
morphism to its associated group (see \cite{FR92}) or to its group of bijections or 
to its group of automorphisms.   

In order to 
formalize the notion of a rack, one needs the diagonal map 
${\rm diag}_X:X\to X\times X$ given by $x\mapsto(x,x)$. Then the
self-distributivity relation reads in terms of maps  
\begin{eqnarray*}
\lefteqn{\mathbf{m}\circ (\mathrm{id}_M \times \mathbf{m})}  \\
  & =  & \mathbf{m}\circ (\mathbf{m}\times \mathbf{m})\circ 
     (\mathrm{id}_M\times\tau_{M,M}\times \mathrm{id}_M) \circ 
     (\mathrm{diag}_M\times\mathrm{id}_M\times \mathrm{id}_M). 
\end{eqnarray*}
Axiomatizing this kind of structure, one may start with a coalgebra $C$ 
and look for rack operations on this fixed coalgebra, see \cite{CCES08}
and \cite{L12}.

A natural framework where this kind of structure arises (as we show in 
Section 3) is by taking point-distributions (i.e. applying the 
{\it Serre functor}) over (resp. to) the pointed manifold given by a Lie rack.
We dub the arising structure as {\it rack bialgebra}.
  
Lie racks are intimately related to
{\it Leibniz algebras} ${\mathfrak h}$, i.e. a vector space ${\mathfrak h}$ with a 
bilinear bracket $[,]:{\mathfrak h}\otimes{\mathfrak h}\to{\mathfrak h}$ such that
for all $X,Y,Z\in{\mathfrak h}$, $[X,-]$ acts as a derivation:
\begin{equation}    \label{LeibnizIdentity}
[X,[Y,Z]]\,=\,[[X,Y],Z]+[Y,[X,Z]].
\end{equation}  
Indeed, Kinyon showed in \cite{Kin07} that the tangent space at $e\in H$ of a Lie rack $H$ carries 
a natural structure of a Leibniz algebra, generalizing the relation between a Lie group 
and its tangent Lie algebra. Conversely, every (finite dimensional real or complex) 
Leibniz algebra ${\mathfrak h}$
may be integrated into a Lie rack $R_{\mathfrak h}$ (with underlying manifold ${\mathfrak h}$) 
using the rack product
\begin{equation}   \label{rack_structure}
X\blacktriangleright Y\,:=\,e^{{\rm ad}_X}(Y),
\end{equation}
noting that the exponential of the inner derivation ${\rm ad}_X$ for each $X\in{\mathfrak h}$
is an automorphism. 

Another closely related algebraic structure is that of dialgebras. A {\it dialgebra} is 
a vector space $D$ with two (bilinear) associative operations $\vdash:D\times D\to D$ 
and $\dashv:D\times D\to D$
which satisfy three compatibility relations, namely for $a,b,c\in D$:
$$(a\vdash b)\vdash c\,=\,(a\dashv b)\vdash c,\,\,\,\,a\dashv (b\vdash c)\,=\,a\dashv (b\dashv c),
\,\,\,\,(a\vdash b)\dashv c\,=\,a\vdash (b\dashv c).$$
A dialgebra $D$ becomes a Leibniz algebras via the formula
$$[a,b]\,=\,a\vdash b - b\dashv a.$$
In this sense $\vdash$ and $\dashv$ are two halves of a Leibniz bracket.  Loday and Goichot 
have defined an enveloping dialgebra of a Leibniz algebra, see \cite{Goi2001}, \cite{Lod2001}.

One main point of this paper is the link between rack bialgebras and 
cocommutative Hopf dialgebras.  In Theorem 2.5, we adapt 
Suschkewitsch's Theorem in semi-group theory to the present context. The classical result 
(see Appendix B) treats semi-groups
with a left unit $e$ and right inverses (analoguous results in the left-context), called 
{\it right groups}. Suschkewitsch shows that such a right group $\Gamma$ decomposes as a product
$\Gamma=\Gamma e\times E$ where $E$ is the set of all idempotent elements.  

Its incarnation here
shows that a cocommutative right Hopf algebra ${\mathcal H}$ decomposes as a tensor product 
${\mathcal H}{\bf 1}\otimes E_{\mathcal H}$ where $E_{\mathcal H}$ is the subspace of generalized
idempotents.    

Furthermore, we will show in Theorem 2.6 how to associate to any 
augmented rack bialgebra an augmented cocommutative Hopf dialgebra. In Theorem 2.7, we 
investigate what Suschkewitsch's decomposition gives for a cocommutative Hopf dialgebra $A$.
It turns out that $A$ decomposes as a tensor product $E_A\otimes H_A$ of $E_A$ with 
$H_A$ which may be identified to the associative quotient $A_{\rm ass}$ of $A$. 
This result permits to show that the Leibniz algebra of primitives in $A$ is 
a hemi-semi-direct product (see \cite{Kin07}), and thus always split.
In this way we arrive once again at the result that Lie digroups may serve only to integrate
{\it split} Leibniz algebras which has already been observed by Covez in his master thesis
\cite{Cov1}.

Let us comment on the {\bf content of the paper}:\\

All our bialgebra notion are based on the standard theory of coalgebras, some 
features of which as well as our notions are recalled in Appendix A.  
Rack bialgebras and augmented rack bialgebras are studied in Section 2.
Connected, cocommutative Hopf algebras give rise to a special case of rack
bialgebras. In Section 2.2, we associate to any Leibniz algebra ${\mathfrak h}$ an
augmented rack bialgebra $\mathsf{UAR}^{\infty}({\mathfrak h})$ and study the functorial 
properties of this association. This rack bialgebra plays the role of an enveloping 
algebra in our context. It turns out that a truncated, non-augmented version 
$\mathsf{UR}({\mathfrak h})$
is a left adjoint of the functor of primitives ${\rm Prim}$.

We also study the ``group-algebra'' functor associating to a rack $X$ its 
rack bialgebra $K[X]$. Like in the classical framework, $K[-]$ is left adjoint to the 
functor ${\rm Slike}$ associating to a track bialgebra its rack of set-like elements. 
The relation between rack bialgebras and the other algebraic 
notion discussed in this paper is summarized in the diagram (see the end of Section 2.2)
of categories and functors:  

$$
  \xymatrix{ {\rm Lie} \ar@<1ex>[r]^{\mathsf{U}} \ar@{^{(}->}[d]^{i} & 
{\rm Hopf} \ar@{^{(}->}[d]^{j} 
\ar[l]^{\rm Prim} \ar@<1ex>[r]^{\rm Slike} & {\rm Grp} \ar[l]^{K[-]} 
\ar@{^{(}->}[d] \\
{\rm Leib} \ar@<1ex>[r]^{\hspace{-.5cm}\mathsf{UAR}^{\infty}} & {\rm RackBialg}   
\ar[l]^{\hspace{-.5cm}\rm Prim} \ar@<1ex>[r]^{\hspace{.5cm}\rm Slike} & 
{\rm Racks} \ar[l]^{\hspace{.5cm}K[-]}
} 
$$
   
In Section 2.3, we develop the structure theory for rack bialgebras and 
cocommutative Hopf dialgebras, based on Suschkewitsch's Theorem. 
Section 2.3 contains Theorem 2.5, Theorem 2.6 and Theorem 2.7 whose content
we have described above.  

Recollecting basic knowledge about the Serre functor $F$ is the subject of
Section 3. In particular, we show in Section 3.2 that $F$ is a strong monoidal functor
from the category of pointed manifolds $\mathcal{M}f*$ to the category of coalgebras,
based on some standard material on coalgebras (Appendix A). 
In Section 3.3, we apply $F$ to Lie groups, Lie semi-groups, Lie digroups, and to Lie racks and 
augmented Lie racks, and study the additional structure which we obtain on the coalgebra.
The case of Lie racks motivates the notion of rack bialgebra.

Recall that for a Leibniz algebra ${\mathfrak h}$, the vector space ${\mathfrak h}$ becomes 
a Lie rack $R_{\mathfrak h}$ with the rack product
$$X\blacktriangleright Y\,=\,e^{{\rm ad}_X}(Y).$$
In Theorem 3.8, we show that the rack bialgebra $\mathsf{UAR}^{\infty}({\mathfrak h})$
associated to ${\mathfrak h}$ coincides with the rack bialgebra $F(R_{\mathfrak h})$.  \\

{\bf Acknowledgements:}
F.W. thanks Universit\'e de Haute Alsace for an invitation during which the 
shape of this research project was defined. At some point, the subject of 
this paper was joint work of S.R. and F.W. with Simon Covez, and we express our 
gratitude to him for his contributions. Starting from this, S.R. \cite{Riv}
came independently to similar results which we incorporated in this paper.  
M.B. thanks Nacer Makhlouf for his question about the relations of rack-bialgebras
to dialgebras, and Gw\'{e}nael Massuyeau for his question about universals.

\section{Several bialgebras}
  \label{SecSeveralBialgebras}

In the following, let $K$ be an associative commutative unital ring
containing all the rational numbers. The symbol
$\otimes$ will always denote the tensor product of $K$-modules over
$K$. For any coalgebra $(C,\Delta)$ over $K$, we shall use 
Sweedler's notation $\Delta(a)=\sum_{(a)}a^{(1)}\otimes a^{(2)}$
for any $a\in A$. 
See also Appendix \ref{AppCoalgebras} for a survey on definitions and notations
in coalgebra theory.

The following sections will all deal with the following type of 
\emph{nonassociative bialgebra}: Let
$(B,\Delta,\epsilon,\mathbf{1},\mu)$ be a $K$-module such that
$(B,\Delta,\epsilon,\mathbf{1})$ is a \emph{coassociative counital coaugmented
coalgebra} (a $C^3$-coalgebra), and such that the linear map 
$\mu:B\otimes B\to B$
(the multiplication) is a morphism of $C^3$-coalgebras (it satisfies in 
particular 
$\mu(\mathbf{1}\otimes\mathbf{1})=\mathbf{1}$). We shall call this situation
a \emph{nonassociative $C^3I$-bialgebra} (where $I$ stands for $\mathbf{1}$
being an \emph{idempotent} for the multiplication $\mu$). For another
nonassociative $C^3I$-bialgebra $(B',\Delta',\epsilon',\mathbf{1}',\mu')$
a $K$-linear map $\phi:B\to B'$ will be called a \emph{morphism of
nonassociative $C^3I$-bialgebras} iff it is a morphism of $C^3$-coalgebras
and is multiplicative in the usual sense 
$\phi\big(\mu(a\otimes b)\big)=\mu'\big(\phi(a)\otimes\phi(b)\big)$
for all $a,b\in B$. The nonassociative $C^3I$-bialgebra
$(B,\Delta,\epsilon,\mathbf{1})$ is called \emph{left-unital} 
(resp.~\emph{right-unital}) iff for all $a\in B$ $\mu(\mathbf{1}\otimes a)=a$
(resp.~$\mu(a\otimes \mathbf{1})=a$).\\
Moreover, consider the associative algebra $A:=\mathrm{Hom}_K(B,B)$ equipped with the 
composition of $K$-linear maps, and the identity map $\mathrm{id}_B$
as the unit element. There is an associative convolution multiplication $*$ in the $K$-module
$\mathrm{Hom}_K(B,A)$ of all $K$-linear maps $B\to \mathrm{Hom}_K(B,B)$, see Appendix 
\ref{AppCoalgebras}, eqn (\ref{EqDefConvolution}) for a definition with
$\mathrm{id}_B\epsilon$ as the unit element. For a  given 
nonassociative $C^3I$-bialgebra $(B,\Delta,\epsilon,\mathbf{1},\mu)$ we can 
consider the map $\mu$ as a map $B\to \mathrm{Hom}_K(B,B)$ in two ways:
as \emph{left multiplication map} 
$L^\mu:b\mapsto \big(b'\mapsto \mu(b\otimes b')\big)$ or
as \emph{right multiplication map} 
$R^\mu:b\mapsto \big(b'\mapsto \mu(b'\otimes b)\big)$. We call 
$(B,\Delta,\epsilon,\mathbf{1},\mu)$ a \emph{left-regular} 
(resp.~\emph{right-regular})
nonassociative $C^3I$-bialgebra iff the map $L^\mu$ (resp.~the map $R^\mu$)
has a convolution inverse, i.e. iff there is a $K$-linear map 
$\mu': B\otimes B\to B$ (resp.~$\mu'': B\otimes B\to B$) 
such that $L^\mu*L^{\mu'}=\mathrm{id}_B\epsilon=
L^{\mu'}*L^\mu$ (resp.~$R^\mu*R^{\mu''}=\mathrm{id}_B\epsilon=
R^{\mu''}*R^\mu$), or on elements $a,b\in B$ for the left regular case
\begin{equation} \label{EqDefLeftRightRegularC3IBialgebras}
     \sum_{(a)}\mu\big(a^{(1)}\otimes \mu'(a^{(2)}\otimes b)\big)=
       \epsilon(a)b= \mu'\big(a^{(1)}\otimes \mu(a^{(2)}\otimes b)\big).
 \end{equation}
 Note that every associative unital Hopf algebra 
 $(H,\Delta,\epsilon,\mathbf{1},\mu,S)$
 (where $S$ denotes the antipode, i.e. the convolution inverse of the identity
 map in $\mathrm{Hom}_K(H,H)$) is right- and left-regular by setting
 $\mu'=\mu\circ (S\otimes \mathrm{id}_H)$ and 
 $\mu''=\mu\circ (\mathrm{id}_H\otimes S)$.
\begin{lemma}\label{LRegularityOfLeftOrRightMultiplications}
  Let $(B,\Delta,\epsilon,\mathbf{1},\mu)$ be a nonassociative $C^3I$-bialgebra.
  \begin{enumerate}
   \item If $B$ is left-regular (resp.~right-regular), then the corresponding
     $K$-linear map $\mu':B\otimes B\to B$ is unique, and in case $\Delta$ is
     cocommutative, $\mu'$ is map of 
     $C^3$-coalgebras.
   \item If $(B,\Delta,\epsilon,\mathbf{1},\mu)$ is left-unital (rep.~right-unital)
     and its underlying $C^3$-coalgebra is connected, then 
     $(B,\Delta,\epsilon,\mathbf{1},\mu)$ is always
     left-regular (resp.~right-regular).
  \end{enumerate}
\end{lemma}
\begin{prooof}
 {\it 1.} In any monoid (in particular in the convolution monoid)
  two-sided inverses are always unique.
  Moreover, as can easily be checked, a $K$-linear map $\phi:B\otimes B\to B$ 
  is a morphism of coalgebras iff for each $b\in B$
  \begin{equation}\label{EqCompPhiMorphC3CoalgByLeftMultiplication}
     \sum_{(b)}\big(L^\phi(b^{(1)})\otimes L^\phi(b^{(2)})\big)\circ \Delta
        =\Delta\circ L^\phi(b).
  \end{equation}
  (and analogously for right multiplications).
  Both sides of the preceding equation, seen as maps of $b$, are in
  $\mathrm{Hom}_K\big(B,\mathrm{Hom}_K(B,B\otimes B)\big)$. Since
  $\mathrm{Hom}_K(B,B\otimes B)$ is an obvious 
  right $\mathrm{Hom}_K(B,B)$-module, the $K$-module 
  $\mathrm{Hom}_K\big(B,\mathrm{Hom}_K(B,B\otimes B)\big)$ is a right 
  $\mathrm{Hom}_K\big(B,\mathrm{Hom}_K(B,B)\big)$-module with respect to
  the convolution. Define the $K$-linear map $F^{\mu'}:B\to
  \mathrm{Hom}_K(B,B\otimes B)$ by 
  \[
      F^{\mu'}(b):=\sum_{(b)}\big(L^{\mu'}(b^{(1)})\otimes L^{\mu'}(b^{(2)})\big)
      \circ \Delta
        -\Delta\circ L^{\mu'}(b).
  \]
  Using eqn (\ref{EqCompPhiMorphC3CoalgByLeftMultiplication}), the fact
  that $L^{\mu'}$ is a convolution inverse of $L^{\mu}$, and the cocommutativity
  of $\Delta$, we get
  \[
       F^{\mu'}*L^{\mu}=0,~~\mathrm{hence} ~~
         0=F^{\mu'}*L^{\mu}*L^{\mu'}=F^{\mu'}*(\mathrm{id}_B\epsilon)=F^{\mu'}
  \]
  and $\mu'$ preserves comultiplications. A similar reasoning where
  $B\otimes B$ is replaced by $K$ shows that $\mu'$ preserves counits. Finally,
  it is obvious that $L^{\mu'}(\mathbf{1})$ is the inverse
  of the $K$-linear map $L^\mu(\mathbf{1})$, and since the latter fixes 
  $\mathbf{1}$ so does the former. The reasoning for right-regular bialgebras
  is completely analogous.\\
  {\it 2.} For left-unital bialgebras we get $L^\mu(\mathbf{1})=\mathrm{id}_B$,
  and the generalized Takeuchi-Sweedler argument, see Appendix \ref{AppCoalgebras},
  shows that $L^\mu$ has a convolution inverse. Right-unital bialgebras are
  treated in  an analogous manner.
\end{prooof}

Note that any $C^3$-coalgebra $(C,\Delta,\epsilon,\mathbf{1})$ becomes
a left-unital (resp.~right-unital) associative $C^3I$-bialgebra by
equipping with the \emph{left-trivial} (resp.~\emph{right-trivial}) 
multiplication
\begin{equation}\label{EqDefLeftTrivialMultiplication}
   \mu_0(a\otimes b):= \epsilon(a)b~~~~~\big(\mathrm{resp.}~~\mu_0(a\otimes b)
       := \epsilon(b)a\big).  
\end{equation}
We shall call an element $c\in B$ a \emph{generalized idempotent} iff
$\sum_{(c)}c^{(1)}c^{(2)}=c$. Moreover $c\in B$ will be called
a \emph{generalized left (resp.~right) unit element} iff for all $b\in B$
we have $cb=\epsilon(c)b$ (resp. $bc=\epsilon(c)b$).

\subsection{Rack bialgebras, augmented rack bialgebras and Leibniz algebras}
 \label{SubSecRackBialgebrasETC}

\begin{defi}
  A \textbf{rack bialgebra} $(B,\Delta,\epsilon,\mathbf{1},\mu)$ 
  is a nonassociative
  $C^3I$-bialgebra (where we write
  for all $a,b\in B$
  $\mu(a\otimes b)=:a\triangleright b$) such that the following identities
  hold for all $a,b,c\in B$
  \begin{eqnarray}
   \mathbf{1}\triangleright a &=& a,  \label{EqRackBialgOnePreservesAll}\\
   a\triangleright \mathbf{1} & = & \epsilon(a)\mathbf{1}, \label{EqRackBialgAllPreserveOne} \\
   a \triangleright (b\triangleright c) &=&
     \sum_{(a)} (a^{(1)}\triangleright b)\triangleright
               (a^{(2)}\triangleright c). \label{EqDefAlgSelfDistributivity}
 \end{eqnarray}
 The last condition (\ref{EqDefAlgSelfDistributivity}) is called the 
 \textbf{self-distributivity} condition.
\end{defi}
Note that we do not demand that the $C^3$-coalgebra $B$ should be cocommutative
nor connected. 

\begin{exem}
Any $C^3$ coalgebra $(C,\Delta,\epsilon,\mathbf{1})$ carries
a \emph{trivial rack bialgebra structure} defined by the left-trivial 
multiplicaton
\begin{equation}
    a \triangleright_0 b := \epsilon(a)b
\end{equation}
which in addition is easily seen to be associative and left-unital, but in general
not unital.
\end{exem}

Another method of constructing rack bialgebras is \emph{gauging}:
Let $(B,\Delta,\epsilon,\mathbf{1},\mu)$ a rack bialgebra --where we write 
$\mu(a\otimes b)=a\triangleright b$ for all $a,b\in B$--, and let
$f:B\to B$ a morphism of $C^3$-coalgebras such that for all $a,b\in B$
\begin{equation}
   f(a\triangleright b)= a\triangleright \big(f(b)\big),
\end{equation}
i.e. $f$ is $\mu$-equivariant. 
It is a routine check that $(B,\Delta,\epsilon,\mathbf{1},\mu_f)$ is a 
rack bialgebra where for all $a,b\in B$ the multiplication is defined by
\begin{equation}
   \mu_f(a\otimes b):=a \triangleright_f b := \big(f(a)\big)\triangleright b.
\end{equation}
We shall call $(B,\Delta,\epsilon,\mathbf{1},\mu_f)$ the 
$f$-\emph{gauge
of $(B,\Delta,\epsilon,\mathbf{1},\mu)$}.

\begin{exem}
Let $(H,\Delta_H,\epsilon_H,\mu_H,\mathbf{1}_H,S)$ be a 
cocommutative Hopf algebra
over $K$. Then it is easy to see (cf.~also the particular case $B=H$ and 
$\Phi=\mathrm{id}_H$ of 
Proposition \ref{PAugmentedRackBialgebraIsRackBialgebra}) 
that the new multiplication 
$\mu:H\otimes H\to H$, written $\mu(h\otimes h')=h\triangleright h'$,
defined by the usual \emph{adjoint representation}
\begin{equation}\label{EqSecondDefAdjointRep}
    h\triangleright h' := \mathrm{ad}_h(h')
       :=\sum_{(h)}h^{(1)}h'\big(S(h^{(2)})\big),
\end{equation}
equips the $C^4$-coalgebra $(H,\Delta_H,\epsilon_H,\mathbf{1}_H)$
with a rack bialgebra structure.
\end{exem}

In general, the adjoint representation does not seem to preserve the 
coalgebra structure if no cocommutativity is
assumed.

\begin{exem}  \label{rack_algebra}
Recall that a pointed set $(X,e)$ is a {\it pointed rack} in case there is a
binary operation 
$\rhd:X\times X\to X$ such that for all $x\in X$, the map $y\mapsto x\rhd y$ 
is bijective and such that for all $x,y,z\in X$,
the self-distributivity and unit relations
$$
   x\rhd(y\rhd z)\,=\,(x\rhd y)\rhd(x\rhd z),~~~e\triangleright x = x,~~~
   \mathrm{and}~~~x\triangleright e = e
$$
are satisfied. Then there is a natural rack bialgebra structure on the vector space
$K[X]$ which has the elements of $X$ as a basis. $K[X]$ carries the usual coalgebra 
structure such that all $x\in X$ are set-like: $\triangle(x)=x\otimes x$ for all 
$x\in X$. The product $\mu$ is then induced by the rack product. 
By functoriality, $\mu$ is compatible with $\triangle$ and $e$. 

Observe that this construction differs slightly from the construction in 
\cite{CCES08}, Section 3.1. 
\end{exem}    
 
More generally there is the following structure:
\begin{defi}
 An \textbf{augmented rack bialgebra} over $K$ is a quadruple
 $(B,\Phi,H,\ell)$ consisting
 of a $C^3$-coalgebra $(B,\Delta,\epsilon,\mathbf{1})$, of a 
 cocommutative (!)
 Hopf algebra
 $(H,\Delta_H,\epsilon_H,\mathbf{1}_H,\mu_H,S)$, of a morphism of
 $C^3$-coalgebras $\Phi:B\to H$, and of a left action 
 $\ell:H \otimes B\to B$ of $H$ on $B$ which is a morphism of 
 $C^3$-coalgebras (i.e. $B$ is a $H$-module-coalgebra) 
 such that for all $h\in H$ and $a\in B$
 \begin{eqnarray}
     h.\mathbf{1} & = &\epsilon_H(h)\mathbf{1} 
            \label{EqDefAugBialgHPreservesE} \\
    \Phi(h.a) & = & \mathrm{ad}_h\big(\Phi(a)\big)
         \label{EqDefAugBialgPhiIntertwinesActionWithAdjoint}.
 \end{eqnarray}
 where $\mathrm{ad}$ denotes the usual adjoint representation for  
 Hopf algebras, see e.g. eqn (\ref{EqSecondDefAdjointRep}).
 
 \noindent We shall define a \textbf{morphism} 
 $(B,\Phi,H,\ell)\to (B',\Phi',H',\ell')$
\textbf{of augmented rack bialgebras} to be a pair $(\phi,\psi)$ 
of $K$-linear maps
where $\phi:(B,\Delta,\epsilon,\mathbf{1})\to (B',\Delta',\epsilon',\mathbf{1}')$
is a morphism of $C^3$-coalgebras, and $\psi:H\to H'$ is a morphism of Hopf 
algebras such that the obvious diagrams commute:
\begin{equation}\label{EqDefMorphAugmentedRackBialg}
   \Phi'\circ \phi = \psi\circ \Phi,~~~\mathrm{and}~~
         \ell'\circ (\psi\otimes \phi) = \phi\circ \ell 
\end{equation}
\end{defi}
An immediate consequence of this definition is the following
\begin{prop}\label{PAugmentedRackBialgebraIsRackBialgebra}
 Let $(B,\Phi,H,\ell)$ be an augmented rack bialgebra. Then the 
 $C^3$-coalgebra $(B,\epsilon,\mathbf{1})$ will become a left-regular 
 rack bialgebra
by means of the multiplication
\begin{equation}\label{EqDefRackBialgMultiplicationFromAugmented}
    a\triangleright b := \Phi(a).b
\end{equation}
for all $a,b\in B$. In particular, each Hopf algebra $H$ becomes
an augmented rack bialgebra via $(H,\mathrm{id}_H,H,\mathrm{ad})$. In general,
for each augmented rack bialgebra the map
$\Phi:B\to H$ is a morphism of rack bialgebras.
\end{prop}
\begin{prooof} 
For a proof of this proposition, see \cite{ABSW}.
\end{prooof}

\begin{exem}  \label{augmented_rack_algebra}
Exactly in the same way as a pointed rack gives rise to a rack bialgebra $K[X]$,
an augmented pointed rack $p:X\to G$ gives rise to an augmented rack bialgebra
$p:K[X]\to K[G]$. 
\end{exem} 

\begin{rem}
Motivated by the fact that the augmented racks $p:X\to G$ are exactly the 
Yetter-Drinfeld modules over the (set-theoretical) Hopf algebra $G$, we may 
ask about the relation of augmented rack bialgebras to 
Yetter-Drinfeld modules, and more generally of rack bialgebras to 
the Yang-Baxter equation. For these subjects, see \cite{ABSW}. 
\end{rem}  

The link to Leibniz algebras is contained in the following

\begin{prop}   \label{proposition_primitives}
 Let $(B,\Delta,\epsilon,\mathbf{1},\mu)$ be a rack bialgebra
 over $K$.
 \begin{enumerate}
  \item Then its $K$-submodule of all primitive elements, 
  $\mathsf{Prim}(B)=:\mathfrak{h}$,
 (see eqn (\ref{EqDefPrimitives}) of Appendix \ref{AppCoalgebras}) 
 is a subalgebra with respect to $\mu$ (written $a\triangleright b$)
 satisfying the 
 \textbf{(left) Leibniz identity}
 \begin{equation}\label{EqDefLeibnizIdentity}
  x\triangleright (y\triangleright z)=
     (x\triangleright y)\triangleright z
         + y\triangleright (x\triangleright z)
 \end{equation}
 for all $x,y,z\in \mathfrak{h}=\mathsf{Prim}(B)$. Hence the pair
 $(\mathfrak{h},[~,~])$ with $[x,y]:=x\triangleright y$ for all 
 $x,y\in \mathfrak{h}$ is a \textbf{Leibniz algebra} over $K$.
 Moreover, every morphism of rack bialgebras maps primitive elements
 to primitive elements and thus induces a morphism of Leibniz algebras.
 \item More generally, $\mathfrak{h}$ and each subcoalgebra of order 
 $k\in \mathbb{N}$, $B_{(k)}$,  
 (see eqn (\ref{EqDefSubcoalgebraOfOrderK}) of Appendix \ref{AppCoalgebras}) 
 is stable by left 
   $\triangleright$-multiplications with every $a\in B$.
   In particular, each $B_{(k)}$ is a rack subbialgebra of 
   $(B,\Delta,\epsilon,\mathbf{1},\mu)$.
 \end{enumerate}
 \end{prop}
\begin{prooof}
 {\it 2.} Let $x\in \mathfrak{h}$ and $a\in B$. Since $\mu$
 is a morphism of $C^3$-coalgebras and $x$ is primitive, we get
 \begin{eqnarray*}
        \Delta(a\triangleright x)
          &=& \sum_{(a)}
               (a^{(1)}\triangleright x)
                 \otimes (a^{(2)}\triangleright\mathbf{1})
                 +
              \sum_{(a)}
              (a^{(1)}\triangleright \mathbf{1})
                 \otimes (a^{(2)}\triangleright x) \\
       & \stackrel{(\ref{EqRackBialgAllPreserveOne})}{=} &
       \sum_{(a)}
               \big((a^{(1)}\epsilon(a^{(2)}))\triangleright x\big)
                 \otimes \mathbf{1}
                 +
              \sum_{(a)}
               \mathbf{1}
                 \otimes 
                 \big((\epsilon(a^{(1)})a^{(2)})\triangleright x\big) 
                     \\
        & = &
          (a\triangleright x)\otimes \mathbf{1}
            ~ + ~
             \mathbf{1}\otimes(a\triangleright x),
 \end{eqnarray*}
 whence $a\triangleright x$ is primitive. For the statement on the $B_{(k)}$,
 we proceed by induction: For $k=0$, this is clear. Suppose the statement is true
 until $k\in \mathbb{N}$, and let $x\in B_{(k+1)}$. Then
 \begin{eqnarray*}
  \lefteqn{\Delta(a\triangleright x) - (a\triangleright x)\otimes \mathbf{1}
       - \mathbf{1}\otimes(a\triangleright x)} \\
          & = & \sum_{(a)(x)}\Big(
          (a^{(1)}\triangleright x^{(1)})\otimes
          (a^{(2)}\triangleright x^{(2)})
          -(a^{(1)}\triangleright x)\otimes
          (a^{(2)}\triangleright \mathbf{1})\\
          & &
          ~~~~~~~~~-(a^{(1)}\triangleright \mathbf{1})\otimes
          (a^{(2)}\triangleright x)\Big) \\
          & = &
           \big(\Delta(a)\big)\triangleright
            \big(\Delta(x) -x\otimes \mathbf{1}
              -\mathbf{1}\otimes x\big) \\
          & = & \sum_{(a)(x)} (a^{(1)} \triangleright x^{(1)'})\otimes
          (a^{(2)}\triangleright x^{(2)'})
 \end{eqnarray*}
 where we have used the extended multiplication (still denoted $\rhd$) 
$\rhd:(B\otimes B)\otimes (B\otimes B)\to (B\otimes B)$
and set
 \[
     \Delta(x) -x\otimes \mathbf{1}
              -\mathbf{1}\otimes x
                 =: \sum_{(x)}x^{(1)'}\otimes x^{(2)'}
                   ~~\in~B_{(k)}\otimes B_{(k)} 
 \]
 by the definition of $B_{(k+1)}$, see Appendix \ref{AppCoalgebras}.
 By the induction hypothesis, all the terms $a^{(1)}\triangleright x^{(1)'}$
 and $a^{(2)}\triangleright x^{(2)'}$ are in $B_{(k)}$, whence
 $\Delta(a\triangleright x) - (a\triangleright x)\otimes \mathbf{1}
       - \mathbf{1}\otimes(a\triangleright x)$ is in $B_{(k)}\otimes B_{(k)}$,
 implying that $a\triangleright x$ is in $B_{(k+1)}$.\\
 {\it 1.} It follows from {\it 2.} that $\mathfrak{h}$
 is a subalgebra with respect to $\mu$. Let $x,y,z\in\mathfrak{h}$. Then
 since $x$ is primitive, it follows from 
 $\Delta(x)=x\otimes \mathbf{1}+ \mathbf{1}\otimes x$ and the 
 self-distributivity identity (\ref{EqDefAlgSelfDistributivity}) that
 \[
     x\triangleright (y\triangleright z)=
       (x\triangleright y) \triangleright (\mathbf{1}\triangleright z)
       +(\mathbf{1}\triangleright y) \triangleright (x\triangleright z)
       \stackrel{(\ref{EqRackBialgOnePreservesAll})}{=}
       (x\triangleright y)\triangleright z
          + y\triangleright (x\triangleright z).
 \]
 proving the left Leibniz identity. The morphism statement is clear,
 since each morphism of rack bialgebras is a morphism of $C^3$-coalgebras
 and preserves primitives.
\end{prooof}

Leibniz algebras have been invented by A. M. Blokh \cite{Blo} in 1965, and then rediscovered
by J.-L. Loday in 1992 in the search of an explanation for the absence of periodicity in 
algebraic K-Theory \cite[p.323, eqn (10.6.1.1)']{Lod92}.

As an immediate consequence, we get that the functor $\mathsf{Prim}$
induces a functor from the category of all rack bialgebras over $K$
to the category of all Leibniz algebras over $K$.

\begin{rem}  \label{slikes}
Define {\it set-like elements} to be elements $a$ in a rack bialgebra $B$
such that $\Delta(a)=a\otimes a$. Thanks to the fact that $\rhd$ is a morphism
of coalgebras, the set of set-like elements 
${\rm Slike}(B)$ is closed under $\rhd$. In fact, ${\rm Slike}(B)$ is a rack,
and one obtains in this way a functor 
${\rm Slike}:{\rm Rack Bialg}\to{\rm Racks}$. 
\end{rem}

\begin{prop}
The functor of set-likes ${\rm Slike}:{\rm RackBialg}\to{\rm Racks}$ has the functor 
$K[-]:{\rm Racks}\to{\rm RackBialg}$ (see Example \ref{rack_algebra}) as its left-adjoint.
\end{prop}

\begin{prooof}
This follows from the adjointness of the same functors, seen as functors between the categories
of pointed sets and of $C^4$-coalgebras, observing that the $C^4$-coalgebra morphism induced by 
a morphism of racks respects the rack product.
\end{prooof}     

Observe that the restriction of ${\rm Slike}:{\rm Rack Bialg}\to{\rm Racks}$ to 
the subcategory of connected, cocommutative Hopf algebras ${\rm Hopf}$ 
(where the Hopf algebra is given the rack product defined in 
eqn (\ref{EqSecondDefAdjointRep})) gives the usual functor of group-like elements.     

\subsection{(Augmented) rack bialgebras for any Leibniz algebra}

Let $(\mathfrak{h}, [~,~])$ be a \emph{Leibniz algebra over $K$}, i.e.
$\mathfrak{h}$ is a $K$-module equipped with a $K$-linear map 
$[~,~]:\mathfrak{h}\otimes \mathfrak{h}\to\mathfrak{h}$ satisfying
the (left) Leibniz identity (\ref{LeibnizIdentity}).

Recall first that each Lie algebra over $K$ is a Leibniz algebra
giving rise to a functor from the category of all Lie algebras to the 
category of all Leibniz algebras.

Furthermore, recall that each Leibniz algebra has two canonical 
$K$-submod\-ules
\begin{eqnarray}
   Q(\mathfrak{h})&:=&
   \big\{x\in\mathfrak{h}~|~\exists~N\in\mathbb{N}\setminus\{0\},~
    \exists~\lambda_1,\ldots,\lambda_N\in K,~\exists~x_1,\ldots,x_N~
    \nonumber\\
    & &
   ~~~~~~~~~~~~~~~~~~ 
       \mathrm{such~that~}x=\sum_{r=1}^N\lambda_r[x_r,x_r]\big\}, \\
    \mathfrak{z}(\mathfrak{h})
    & := & \big\{x\in\mathfrak{h}~|~\forall~y\in\mathfrak{h}
                     :~[x,y]=0\big\}.
\end{eqnarray}
It is well-known and not hard to deduce from the Leibniz identity 
that both
$Q(\mathfrak{h})$ and $\mathfrak{z}(\mathfrak{h})$ are two-sided
abelian ideals of $(\mathfrak{h}, [~,~])$, that
$Q(\mathfrak{h})\subset \mathfrak{z}(\mathfrak{h})$, and that
the quotient Leibniz algebras
\begin{equation}
    \overline{\mathfrak{h}}:=\mathfrak{h}/Q(\mathfrak{h})
    ~~~\mathrm{and}~~~
       \mathfrak{g}(\mathfrak{h}):=
       \mathfrak{h}/\mathfrak{z}(\mathfrak{h})
\end{equation}
are Lie algebras. Since the ideal $Q(\mathfrak{h})$ is clearly
mapped
into the ideal $Q(\mathfrak{h}')$ by any morphism of Leibniz algebras
$\mathfrak{h}\to\mathfrak{h}'$ (which is a priori not the case
for $\mathfrak{z}(\mathfrak{h})$ !), there is an obvious functor
$\mathfrak{h}\to \overline{\mathfrak{h}}$ from the category of
all Leibniz algebras to the category of all Lie algebras.

In order to perform the following constructions of rack bialgebras for
any given Leibniz
algebra $(\mathfrak{h},[~,~])$, choose first a two-sided ideal
$\mathfrak{z}\subset \mathfrak{h}$ such that
\begin{equation}
    Q(\mathfrak{h})\subset \mathfrak{z} \subset 
    \mathfrak{z}(\mathfrak{h}),
\end{equation}
let $\mathfrak{g}$ denote the quotient Lie algebra 
$\mathfrak{h}/\mathfrak{z}$, and let $p:\mathfrak{h}\to \mathfrak{g}$
be the natural projection. The data of $\mathfrak{z}\subset\mathfrak{h}$,
i.e. of a Leibniz algebra $\mathfrak{h}$ together with an ideal $\mathfrak{z}$
such that $Q(\mathfrak{h})\subset\mathfrak{z}\subset\mathfrak{z}(\mathfrak{h})$,
could be called an {\it augmented Leibniz algebra}. Thus we are actually associating 
an augmented rack bialgebra to every augmented Leibniz algebra. 
In fact, we will see that this augmented rack bialgebra does not depend on the choice of the ideal 
$\mathfrak{z}$ and therefore refrain from introducing augmented Leibniz algebras in a more 
formal way.

The Lie algebra $\mathfrak{g}$ naturally
acts as derivations on $\mathfrak{h}$ by means of 
(for all $x,y\in\mathfrak{h}$)
\begin{equation} \label{EqDefLeibnizHIsAGModule}
    p(x).y := [x,y]=:\mathrm{ad}_x(y)
\end{equation}
because $\mathfrak{z} \subset \mathfrak{z}(\mathfrak{h})$. Note that
\begin{equation}\label{EqCompHModLeftCenterCongAd}
 \mathfrak{h}/\mathfrak{z}(\mathfrak{h})\cong 
   \big\{\mathrm{ad}_x\in\mathrm{Hom}_K(\mathfrak{h},\mathfrak{h})~|~
         x\in\mathfrak{h}\big\}.
\end{equation}
as Lie algebras.

Consider 
now the $C^5$-coalgebra 
$(B=\mathsf{S}(\mathfrak{h}),\Delta,\epsilon,\mathbf{1})$ which is
actually a commutative cocommutative Hopf algebra over $K$ with respect
to the symmetric multiplication $\bullet$. The linear map 
$p:\mathfrak{h}\to \mathfrak{g}$ induces a unique morphism of
Hopf algebras
\begin{equation}
   \tilde{\Phi}=\mathsf{S}(p):
     \mathsf{S}(\mathfrak{h})\to \mathsf{S}(\mathfrak{g})
\end{equation}
satisfying
\begin{equation}\label{EqDefMapTildePhi}
    \tilde{\Phi}(x_1\bullet\cdots\bullet x_k)=p(x_1)\bullet\cdots\bullet p(x_k)
\end{equation}
for any nonnegative integer $k$ and $x_1,\ldots,x_k\in \mathfrak{h}$.
In other words, the association $\mathsf{S}:V\to \mathsf{S}(V)$
is a functor from the category of all $K$-modules to the category
of all commutative unital $C^5$-coalgebras.
Consider now the universal enveloping algebra $\mathsf{U}(\mathfrak{g})$
of the Lie algebra $\mathfrak{g}$. Since $\mathbb{Q}\subset K$ by 
assumption, the Poincar\'{e}-Birkhoff-Witt Theorem (in short: PBW) holds
(see e.g. \cite[Appendix]{Qui69}). More precisely, the symmetrisation
map $\omega:\mathsf{S}(\mathfrak{g})\to \mathsf{U}(\mathfrak{g})$,
defined by
\begin{equation}\label{EqDefSymmetrisationMap}
   \omega(\mathbf{1}_{\mathsf{S}(\mathfrak{g})})
      = \mathbf{1}_{\mathsf{U}(\mathfrak{g})},~~~\mathrm{and}~~~
    \omega(\xi_1\bullet\cdots\bullet \xi_k)
    = \frac{1}{k!}\sum_{\sigma\in S_k}
        \xi_{\sigma(1)}\cdots\xi_{\sigma(k)},
\end{equation}
see e.g. \cite[p.80, eqn (3)]{Dix74}, is an isomorphism of $C^5$-coalgebras 
(in general not of associative algebras). We now need an action of the 
Hopf algebra $H=\mathsf{U}(\mathfrak{g})$ on $B$, and an intertwining
map $\Phi:B\to\mathsf{U}(\mathfrak{g})$. In order to get this, we first look
at $\mathfrak{g}$-modules: The $K$-module $\mathfrak{h}$ is a 
$\mathfrak{g}$-module by means of eqn (\ref{EqDefLeibnizHIsAGModule}), the 
Lie algebra $\mathfrak{g}$ is a $\mathfrak{g}$-module via its adjoint
representation, and the linear map $p:\mathfrak{h}\to\mathfrak{g}$ is a 
morphism of $\mathfrak{g}$-modules since $p$ is a morphism of Leibniz algebras.
Now $\mathsf{S}(\mathfrak{h})$ and $\mathsf{S}(\mathfrak{g})$ are
$\mathfrak{g}$-modules in the usual way, i.e. for all 
$k\in\mathbb{N}\setminus \{0\}$, $\xi,\xi_1,\ldots,\xi_k\in\mathfrak{g}$, and
$x_1\ldots,x_k\in\mathfrak{h}$
\begin{eqnarray}
 \xi.(x_1\bullet\cdots\bullet x_k)
 & := & \sum_{r=1}^k
  x_1\bullet\cdots\bullet (\xi.x_r)\bullet \cdots \bullet x_k, 
     \label{EqDefGActionOnSH}\\
  \xi.(\xi_1\bullet\cdots\bullet \xi_k)
 & := & \sum_{r=1}^k
  \xi_1\bullet\cdots\bullet [\xi.\xi_r]\bullet \cdots \bullet \xi_k,
  \label{EqDefGActionOnSG}
\end{eqnarray}
and of course $\xi.\mathbf{1}_{\mathsf{S}(\mathfrak{h})}=0$ and
$\xi.\mathbf{1}_{\mathsf{S}(\mathfrak{g})}=0$. Recall that 
$\mathsf{U}(\mathfrak{g})$ is a $\mathfrak{g}$-module via the adjoint 
representation $\mathrm{ad}_\xi(u)=\xi.u= \xi u - u\xi$ 
(for all $\xi\in\mathfrak{g}$ and
all $u\in \mathsf{U}(\mathfrak{g})$).\\
It is easy to see that the map $\tilde{\Phi}$ (\ref{EqDefMapTildePhi}) is a 
morphism
of $\mathfrak{g}$-modules, and it is well-known that the symmetrization map 
$\omega$ (\ref{EqDefSymmetrisationMap}) is also a morphism of
$\mathfrak{g}$-modules, see e.g. \cite[p.82, Prop. 2.4.10]{Dix74}.
Define the $K$-linear map $\Phi:\mathsf{S}(\mathfrak{h})\to 
\mathsf{U}(\mathfrak{g})$ by the composition
\begin{equation}
    \Phi:= \omega\circ \tilde{\Phi}.
\end{equation}
Then $\Phi$ is a map of $C^5$-coalgebras and a map of $\mathfrak{g}$-modules.
Thanks to the universal property of the universal enveloping algebra, it 
follows
that $\mathsf{S}(\mathfrak{h})$ and $\mathsf{U}(\mathfrak{g})$ are 
left $\mathsf{U}(\mathfrak{g})$-modules, via
(for all $\xi_1,\ldots,\xi_k\in\mathfrak{g}$, and for all $a\in 
\mathsf{S}(\mathfrak{h})$) 
\begin{equation}\label{EqDefUGActionOnSH}
  (\xi_1\cdots\xi_k).a = \xi_1.(\xi_2.(\cdots \xi_k.a)\cdots)
\end{equation}
and the usual adjoint representation (\ref{EqSecondDefAdjointRep})
(for all $u\in \mathsf{U}(\mathfrak{g})$) 
\begin{equation}
 \mathrm{ad}_{\xi_1\cdots\xi_k}(u)
 =\big(\mathrm{ad}_{\xi_1}\circ\cdots\circ\mathrm{ad}_{\xi_k}\big)(u),
\end{equation}
and that $\Phi$ intertwines the $\mathsf{U}(\mathfrak{g})$-action on
$C=\mathsf{S}(\mathfrak{h})$ with the adjoint action of 
$\mathsf{U}(\mathfrak{g})$ on itself.\\
Finally it is a routine check using the above identities 
(\ref{EqDefGActionOnSH}) and (\ref{EqSecondDefAdjointRep}) that 
$\mathsf{S}(\mathfrak{h})$ becomes a module coalgebra.\\
We can resume the preceding considerations in the following
\begin{theorem}\label{TLeibnizToUARInfinityOfLeibniz}
  Let $(\mathfrak{h},[~,~])$ be a Leibniz algebra over $K$, let
  $\mathfrak{z}$ be a two-sided ideal of $\mathfrak{h}$ such that
  $Q(\mathfrak{h})\subset \mathfrak{z} \subset \mathfrak{z}(\mathfrak{h})$,
  let $\mathfrak{g}$ denote the quotient Lie algebra 
  $\mathfrak{h}/\mathfrak{z}$ by $\mathfrak{g}$, and let 
  $p:\mathfrak{h}\to \mathfrak{g}$ be the canonical projection.
  \begin{enumerate}
  \item Then there is a canonical $\mathsf{U}(\mathfrak{g})$-action $\ell$ on
  the $C^5$-coalgebra $B:=\mathsf{S}(\mathfrak{h})$ (making it into a module
  coalgebra leaving invariant $\mathbf{1}$) and a canonical lift of
  $p$ to a map of $C^5$-coalgebras, $\Phi:\mathsf{S}(\mathfrak{h})\to
  \mathsf{U}(\mathfrak{g})$ such that eqn 
  (\ref{EqDefAugBialgPhiIntertwinesActionWithAdjoint}) holds.\\
  Hence the quadruple $(\mathsf{S}(\mathfrak{h}),\Phi,
  \mathsf{U}(\mathfrak{g}), \ell)$ is an augmented rack bialgebra 
  whose associated Leibniz algebra is equal to $(\mathfrak{h},[~,~])$
  (independently of the choice of $\mathfrak{z}$).\\
  The resulting rack multiplication $\mu$ of $\mathsf{S}(\mathfrak{h})$ 
  (written
  $\mu(a\otimes b)=a\triangleright b$) is also independent on the choice of
  $\mathfrak{z}$ and is explicitly given
  as follows for all positive integers $k,l$ and
  $x_1,\ldots,x_k$, $y_1,\ldots,y_l\in\mathfrak{h}$:
  \begin{equation}   \label{adjoint_action_s}
      \big(x_1\bullet\cdots\bullet x_k)\triangleright
      \big(y_1\bullet\cdots\bullet y_l)
      =\frac{1}{k!}\sum_{\sigma\in S_k}
        \big(\mathrm{ad}^s_{x_{\sigma(1)}}\circ\cdots\circ 
          \mathrm{ad}^s_{x_{\sigma(k)}}\big)
          \big(y_1\bullet\cdots\bullet y_l)
  \end{equation}
  where $\mathrm{ad}^s_x$ denotes the action of the Lie algebra
  $\mathfrak{h}/\mathfrak{z}(\mathfrak{h})$ (see eqn 
  (\ref{EqCompHModLeftCenterCongAd}))
  on $\mathsf{S}(\mathfrak{h})$ according to eqn (\ref{EqDefGActionOnSH}).

  \item In case $\mathfrak{z}=Q(\mathfrak{h})$, the construction mentioned
  in {\it 1.} is
  a functor $\mathfrak{h}\to\mathsf{UAR}^\infty(\mathfrak{h})$ from the 
  category of all Leibniz algebras to the 
  category of all augmented rack bialgebras associating to ${\mathfrak h}$
  the rack bialgebra 
$${UAR}^\infty(\mathfrak{h}):=(\mathsf{S}(\mathfrak{h}),\Phi,
  \mathsf{U}(\mathfrak{g}),\ell)$$ 
  and to each morphism
  $f$ of Leibniz algebras the pair $\big(\mathsf{S}(f),
  \mathsf{U}(\overline{f})\big)$ where $\overline{f}$ is the induced Lie 
  algebra morphism.
  \item For each nonnegative integer $k$, the above construction restricts
  to each subcoalgebra of order $k$, $\mathsf{S}(\mathfrak{h})_{(k)}
   =\oplus_{r=0}^k \mathsf{S}^r(\mathfrak{h})$, to define an augmented
   rack bialgebra  $(\mathsf{S}(\mathfrak{h})_{(k)},\Phi_{(k)},
  \mathsf{U}(\mathfrak{g}), 
  \ell|_{\mathsf{U}(\mathfrak{g})\otimes \mathsf{S}(\mathfrak{h})_{(k)}})$
  which in case  $\mathfrak{z}=Q(\mathfrak{h})$ defines a functor
  $\mathfrak{h}\to
  \mathsf{UAR}_{(k)}(\mathfrak{h}):=
  \big(\mathsf{UAR}^\infty(\mathfrak{h})\big)_{(k)}$ from the 
  category of all Leibniz algebras to the 
  category of all augmented rack bialgebras.
  \end{enumerate}
\end{theorem}
\begin{prooof}
For a proof of this theorem, see \cite{ABSW}. 
\end{prooof}

\begin{rem}
This theorem should be compared to Proposition 3.5 in \cite{CCES08}. In \cite{CCES08},
the authors work with the vector space $N:=K\oplus{\mathfrak h}$, while we 
work with the whole symmetric algebra on the Leibniz algebra. In some sense, we extend 
their Proposition 3.5 ``to all orders''.
However, as we shall see below, $N$ is already enough to obtain a left-adjoint to the 
functor of primitives. 
\end{rem}

The above rack bialgebra associated to a Leibniz algebra ${\mathfrak h}$ can be seen
as one version of an {\it enveloping algebra} of ${\mathfrak h}$.

\begin{defi}\label{DefUniversalRackBialgebraInfiniteOrder}
Let ${\mathfrak h}$ be a Leibniz algebra.
We will call the augmented rack bialgebra $(\mathsf{S}(\mathfrak{h}),\Phi,
  \mathsf{U}(\mathfrak{g}), \ell)$ the 
  \textbf{enveloping algebra of infinite order of ${\mathfrak h}$}.
As such, it will be denoted by $\mathsf{UAR}^\infty(\mathfrak {h})$.
\end{defi}   

This terminology is justified, 
for example, by the fact that ${\mathfrak h}$ is identified to the primitives in 
$\mathsf{S}(\mathfrak{h})$ (cf Proposition \ref{proposition_primitives}). This is also 
justified by the following theorem the goal of which
is to show that the enveloping algebra
$\mathsf{UAR}^ \infty({\mathfrak h})$
fits into the following diagram of functors:

$$
 \xymatrix{ {\rm Lie} \ar[r]^{\mathsf{U}} \ar[d]^{i} & {\rm Hopf} \ar[d]^{j} \\
{\rm Leib} \ar[r]^{\hspace{-.5cm}\mathsf{UAR}^ \infty} & {\rm RackBialg}   } 
$$

Here, $i$ is the embedding functor of Lie algebras into Leibniz algebras, and 
$j$ is the embedding
functor of the category of connected, cocommutative Hopf algebras into the category of rack 
bialgebras, using the 
adjoint action (see eqn (\ref{EqSecondDefAdjointRep})) as a rack product. 

\begin{theorem}
Let ${\mathfrak g}$ be a Lie algebra. The PBW isomorphism 
$\mathsf{U}({\mathfrak g})\cong 
\mathsf{S}({\mathfrak g})$ induces an isomorphism of functors
$$
 j\circ\mathsf{U}\,\cong\,\mathsf{UAR}^\infty\circ i.
$$
\end{theorem}

\begin{prooof}
The enveloping algebra $\mathsf{UAR}^\infty({\mathfrak h})$ is by definition the 
functorial version of the rack bialgebra $\mathsf{S}(\mathfrak{h})$, i.e. 
associated to the ideal $Q({\mathfrak h})$. But in case ${\mathfrak h}$ is a 
Lie algebra, $Q({\mathfrak h})=\{0\}$. Then the map $p$ is simply the identity,
and $\mathsf{UAR}^\infty({\mathfrak h})=j\big(U({\mathfrak h})\big)$. 
\end{prooof}

As a relatively easy corollary we obtain from the preceding construction
the computation of \emph{universal rack bialgebras}. More precisely, we look for 
a left adjoint functor for the 
functor $\mathsf{Prim}$, seen as a functor from the category of all
rack bialgebras to the category of all Leibniz algebras. For a given Leibniz 
algebra
$\big(\mathfrak{h},[~,~]\big)$ define the subcoalgebra of order $1$ of the first
component of
$\mathsf{UAR}_{(1)}(\mathfrak{h})$ (see the third statement of
Theorem \ref{TLeibnizToUARInfinityOfLeibniz}), i.e.
\begin{equation}
    \mathsf{UR}(\mathfrak{h}):=\mathsf{UAR}^\infty(\mathfrak{h})_{(1)}
    := K \oplus \mathfrak{h}
\end{equation}
with $\mathbf{1}:=1=1_K$ which is rack subbialgebra according to Proposition
\ref{proposition_primitives}. Its structure reads for all 
$\lambda,\lambda'\in K$ and for all 
$x,y\in\mathfrak{h}$
\begin{eqnarray}
   \Delta(\lambda\mathbf{1} + x) 
            & = &\lambda \mathbf{1}\otimes \mathbf{1} + x\otimes \mathbf{1}
                               + \mathbf{1}\otimes x, \\
    \epsilon(\lambda\mathbf{1}+x) &=& \lambda, \\
    \mu\big((\lambda\mathbf{1} + x)\otimes (\lambda' \mathbf{1}+ x')\big)
      & = &  \lambda\lambda'\mathbf{1} + \lambda x' + [x,x'].
\end{eqnarray}
For the particular case of a Lie algebra $\big(\mathfrak{h},[~,~]\big)$,
the above construction can be found in \cite[Section 3.2]{CCES08}.
Moreover, for any other
Leibniz algebra $\big(\mathfrak{h}', [~,~]'\big)$ and any morphism of Leibniz 
algebras $f:\mathfrak{h}\to \mathfrak{h}'$ define the $K$-linear map
$\mathsf{UR}(f):\mathsf{UR}(\mathfrak{h})\to\mathsf{UR}(\mathfrak{h}')$
as the first component of $\mathsf{UAR}_{(1)}(f)$ (cf.~the third statement of
Theorem \ref{TLeibnizToUARInfinityOfLeibniz})
by
\begin{equation}\label{EqDefUROfF}
  \mathsf{UR}(f)\big(\lambda¸\mathbf{1}+x\big) =  \lambda¸\mathbf{1} + f(x),
\end{equation}
which is clearly is a morphism of rack bialgebras. Hence $\mathsf{UR}$
 is a functor
from the category of all Leibniz algebras to the category of all rack 
bialgebras. Now let $(C,\Delta_C,\epsilon_C,\mathbf{1}_C,\mu_C)$ be a rack 
bialgebra, and let $f:\mathfrak{h}\to \mathsf{Prim}(C)$ be a morphism of 
Leibniz algebras. Define the $K$-linear map 
$\hat{f}:\mathsf{UR}(\mathfrak{h})\to C$
by
\begin{equation}
   \hat{f}(\lambda\mathbf{1}+x)= \lambda \mathbf{1}_C + f(x),
\end{equation}
and it is again a routine check that it defines a morphism of rack bialgebras.
Moreover, due to the almost trivial coalgebra structure of 
$\mathsf{UR}(\mathfrak{h})$, it is clear that any morphism of rack bialgebras
$\mathsf{UR}(\mathfrak{h})\to C$ is of the above form and is uniquely 
determined by 
$\mathsf{Prim}(\hat{f})=f$. Hence we have shown the following 
\begin{theorem}\label{TUniversalLeibnizToRackBialgebra}
 There is a left adjoint functor, $\mathsf{UR}$, for the functor 
 $\mathsf{Prim}$ (associating
 to each Rack bialgebra its Leibniz algebra of all primitive elements).
 For a given Leibniz algebra $\big(\mathfrak{h},[~,~]\big)$, the object
 $\mathsf{UR}(\mathfrak{h})$ --which we shall call the 
 \textbf{Universal Rack Bialgebra of the Leibniz algebra} 
 $\big(\mathfrak{h},[~,~]\big)$--
 has the usual universal properties.
\end{theorem}

Next we can refine the above universal construction by taking into account
the augmented rack bialgebra structure of 
$\mathsf{UAR}^\infty(\mathfrak{h})$ to define another universal object. 
Consider the more detailed category
of all augmented rack bialgebras. Again, the functor $\mathsf{Prim}$ applied
to the coalgebra $B$ (and not to the Hopf algebra $H$) gives a functor from the 
first category to the category of all Leibniz algebras, and we seek again a
left adjoint of this functor, called $\mathsf{UAR}$.
Hence, a natural candidate for a
\emph{universal augmented rack bialgebra associated
to a given Leibniz algebra $\mathfrak{h}$} is
\begin{equation}
    \mathsf{UAR}(\mathfrak{h}):= \mathsf{UAR}_{(1)}(\mathfrak{h})
                        = \big(K\oplus \mathfrak{h}, \Phi_{(1)},
                        \mathsf{U}(\overline{h}),
                        \mathrm{ad}^s|_{\mathsf{U}(\overline{h})
                                         \otimes (K\oplus \mathfrak{h})}\big) .                   
\end{equation}
The third statement of Theorem \ref{TLeibnizToUARInfinityOfLeibniz} tells us
that this is a well-defined augmented rack bialgebra, and that
$\mathsf{UAR}$ is a functor from the category of all Leibniz algebras
to the category of all augmented rack bialgebras.
Now let $(B',\Phi',H',\ell')$ be an augmented rack bialgebra, and let
$f:\mathfrak{h}\to \mathsf{Prim}(B')$ be a morphism of Leibniz algebras.
Clearly, as has been shown in Theorem \ref{TUniversalLeibnizToRackBialgebra},
the map $\hat{f}: \mathsf{UR}(\mathfrak{h})\to B'$ given by eqn 
(\ref{EqDefUROfF}) is a morphism of rack bialgebras. Observe that the morphism
of $C^3$-coalgebras
$\Phi'$ sends the Leibniz subalgebra $\mathsf{Prim}(B')$ of $B'$ into 
$K$-submodule of all primitive elements of the Hopf algebra $H'$,
$\mathsf{Prim}(H')$, which is known to be a Lie subalgebra of
$H'$ equipped with the commutator Lie bracket $[~,~]_{H'}$. Moreover this 
restriction is a morphism of Leibniz algebras. Indeed, for any 
$x',y'\in\mathsf{Prim}(B')$ we have
\begin{eqnarray*}
  \Phi'\big([x',y']'\big) & = & \Phi'(x'\triangleright y')
     = \Phi'\big((\Phi'(x')).y'\big)
     =\mathrm{ad}_{\Phi'(x')}\big(\Phi'(y')\big)\\
    & = & \Phi'(x')\Phi'(y')-\Phi'(y')\Phi'(x') 
        = \big[\Phi'(x'),\Phi'(y')\big]_{H'}.
\end{eqnarray*}
It follows that the two-sided ideal $Q\big(\mathsf{Prim}(B')\big)$
of the Leibniz algebra $\mathsf{Prim}(B')$ is in the kernel of the restriction 
of $\Phi'$ to $\mathsf{Prim}(B')$, whence the map $\Phi'$ induces a well-defined
$K$-linear morphism of Lie algebras $\overline{\Phi'|}:
\overline{\mathsf{Prim}(B')}\to
\mathsf{Prim}(H')$. It follows that the composition 
$\overline{\Phi'|}\circ \overline{f}:\overline{\mathfrak{h}}\to
\mathsf{Prim}(H')\subset H'$
is a morphism of Lie algebras, and by the universal property of universal
envelopping algebras there is a unique morphism of associative unital
algebras $\psi:=\mathsf{U}\big(\overline{\Phi'|}\circ \overline{f}\big):
\mathsf{U}(\overline{\mathfrak{h}})\to H'$. But we have
for all $\xi_1,\ldots,\xi_k\in\overline{\mathfrak{h}}$
\begin{eqnarray*}
 \lefteqn{\Delta_{H'}\big(\psi(\xi_1\cdots\xi_k)\big)
     =  \Delta_{H'}\big(\psi(\xi_1)\cdots\psi(\xi_k)\big)
 =\Delta_{H'}\big(\psi(\xi_1)\big)\cdots\Delta\big(\psi(\xi_k)\big)} \\
 & = & \big(\psi(\xi_1)\otimes \mathbf{1}_{H'} 
             + \mathbf{1}_{H'}\otimes\psi(\xi_1)\big) \cdots
        \big(\psi(\xi_k)\otimes \mathbf{1}_{H'} 
             + \mathbf{1}_{H'}\otimes\psi(\xi_k)\big) \\
 & = & \big(\psi\otimes \psi\big)
    (\xi_1\otimes \mathbf{1}_{\mathsf{U}(\overline{\mathfrak{h}})}
      + \mathbf{1}_{\mathsf{U}(\overline{\mathfrak{h}})} \otimes \xi_1)\cdots
       \big(\psi\otimes \psi\big)
    (\xi_k\otimes \mathbf{1}_{\mathsf{U}(\overline{\mathfrak{h}})}
      + \mathbf{1}_{\mathsf{U}(\overline{\mathfrak{h}})} \otimes \xi_k) \\
 & = & \big(\psi\otimes \psi\big)
   \big(\Delta_{\mathsf{U}(\overline{\mathfrak{h}})}(\xi_1)\big)\cdots
    \big(\psi\otimes \psi\big)
   \big(\Delta_{\mathsf{U}(\overline{\mathfrak{h}})}(\xi_k)\big) \\
  & = & \big(\psi\otimes \psi\big)
       \Big(\Delta_{\mathsf{U}(\overline{\mathfrak{h}})}(\xi_1\cdots \xi_k)\Big) 
\end{eqnarray*}
since $\psi$ maps primitives to primitives whence $\psi$ is a morphism 
of coalgebras. It is easy to check that $\psi$ preserves counits, whence
$\psi$ is a morphism of $C^5$-Hopf-algebras. For all $\lambda\in K$ and
$x\in\mathfrak{h}$ we get
\begin{eqnarray*}
   \big(\psi\circ \Phi_{(1)}\big)(\lambda\mathbf{1}+x)
     & = & \psi\big(\lambda\mathbf{1}_{\mathsf{U}(\overline{\mathfrak{h}})}+p(x)\big)
      =\lambda\mathbf{1}_{H'} 
           +(\overline{\Phi'|}\circ \overline{f})\big(p(x)\big) \\
     & = & \lambda\mathbf{1}_{H'} + \Phi'\big(f(x)\big)
         = \Phi'\big(\hat{f}(\lambda\mathbf{1}+x)\big),
\end{eqnarray*}
showing the first equation $\psi\circ \Phi_{(1)}=\Phi'\circ \hat{f}$ of the 
morphism equation (\ref{EqDefMorphAugmentedRackBialg}). Moreover
for all $\lambda\in K$,
$x\in\mathfrak{h}$, and $u\in \mathsf{U}(\overline{\mathfrak{h}})$ we get
\[
    \hat{f}\big(u.(\lambda\mathbf{1}+x)\big)
    = \hat{f}\big(\lambda\epsilon_{\mathsf{U}(\overline{\mathfrak{h}})}(u)\mathbf{1}
          +u.x\big)
    = \lambda\epsilon_{\mathsf{U}(\overline{\mathfrak{h}})}(u)\mathbf{1}_{C'}
        +f(u.x)
\]
Let $x_1,\ldots,x_k\in \mathfrak{h}$ such that $u=p(x_1)\cdots p(x_k)$.
Then
\begin{eqnarray*}
   f\big(u.x\big)& = & f\big([x_1,[x_2,\ldots [x_k,x]\cdots ]\big)
    = f(x_1)\triangleright \big(f(x_2)\triangleright 
           \cdots \triangleright \big(f(x_k)\triangleright f(x)\big)\cdots\big)\\
   & = & \Big(\Phi'\big(f(x_1)\big)\cdots \Phi'\big(f(x_k)\big)\Big)
            .\big(f(x)\big)\\
   & = &  \Big(\psi\big(p(x_1)\big)\cdots \psi\big(p(x_k)\big)\Big)
   .\big(f(x)\big) = \big(\psi(u)\big).\big(f(x)\big),
\end{eqnarray*}
and therefore
\[
    \hat{f}\big(u.(\lambda\mathbf{1}+x)\big)
      =\big(\psi(u)\big).\big(\hat{f}(\lambda\mathbf{1}+x)\big)
\]
showing the second equation $\psi\circ \Phi_{(1)}=\Phi'\circ \hat{f}$ of the 
morphism equation (\ref{EqDefMorphAugmentedRackBialg}).
It follows that the pair $(\hat{f},\psi)$ is a morphism of augmented
rack bialgebras. We therefore have the following
\begin{theorem}\label{TUniversalLeibnizToAugmentedRackBialgebra}
 There is a left adjoint functor, $\mathsf{UAR}$, for the functor 
 $\mathsf{Prim}$ (associating
 to each augmented rack bialgebra its Leibniz algebra of all primitive elements).
 For a given Leibniz algebra $\big(\mathfrak{h},[~,~]\big)$, the object
 $\mathsf{UAR}(\mathfrak{h})$ --which we shall call the 
 \textbf{Universal Augmented Rack Bialgebra of the Leibniz algebra} 
 $\big(\mathfrak{h},[~,~]\big)$--
 has the usual universal properties.
\end{theorem}

The relationship between the different notions (taking into account also 
Remark (\ref{slikes})) is resumed in the following diagram:

$$
\xymatrix{ 
   {\rm Lie} \ar@<1ex>[r]^{\mathsf{U}} \ar@{^{(}->}[d]^{i} & 
{\rm Hopf} \ar@{^{(}->}[d]^{j} 
\ar[l]^{\rm Prim} \ar@<1ex>[r]^{\rm Slike} & {\rm Grp} \ar[l]^{K[-]} 
                                                       \ar@{^{(}->}[d] \\
{\rm Leib} \ar@<1ex>[r]^{\hspace{-.5cm}\mathsf{UAR}^{\infty}} & {\rm RackBialg}   
\ar[l]^{\hspace{-.5cm}\rm Prim} \ar@<1ex>[r]^{\hspace{.5cm}\rm Slike} & 
{\rm Racks} \ar[l]^{\hspace{.5cm}K[-]}
} 
$$

where $\mathsf{UAR}^{\infty}$ is not left-adjoint to ${\rm Prim}$, while $\mathsf{UR}$ is,
but does not render the square commutative. There is a similar diagram for augmented notions.  

\subsection{Relation with bar-unital di(co)algebras}

In the beginning of the nineties the `enveloping structure' associated to
Leibniz algebras has been the structure of \emph{dialgebras}, 
see e.g.~\cite{Lod2001}. We shall show in this section 
that rack bialgebras and certain
cocommutative Hopf dialgebras are strongly related.

\subsubsection{Left-unital bialgebras and right Hopf algebras}

\noindent 
Let $(B,\Delta,\epsilon,\mathbf{1},\mu)$ be a
nonassociative left-unital $C^3$-bialgebra. It 
 will be called 
\emph{left-unital $C^3$-bialgebra} iff $\mu$ is associative.
In general, $(B,\Delta,\epsilon,\mathbf{1},\mu)$ need not be unital, i.e.
we do not have in general $a\mathbf{1}=a$. However, it is easy to see that
the sub-module $B\mathbf{1}$ of $B$ is a $C^3$-subcoalgebra of 
$(B,\Delta,\epsilon,\mathbf{1})$, and a subalgebra of $(B,\mu)$ such that
$(B\mathbf{1},\Delta',\epsilon',\mathbf{1},\mu')$ is a unital (i.e. left-unital
and right-unital) bialgebra. Here $\Delta'$, $\epsilon'$, and $\mu'$ denote the 
obvious restrictions and corestrictions.\\
In a completely analogous way right-unital $C^3$-bialgebras are defined.

A left-unital (resp.~right unital) \emph{cocommutative} $C^3$-bialgebra 
$(B,\Delta,\epsilon,\mathbf{1},\mu)$ will be called
a \emph{cocommutative right Hopf algebra} (resp.~a
\emph{cocommutative left Hopf algebra}),  $(B,\Delta,\epsilon,\mathbf{1},\mu,S)$,
iff there is a \emph{right antipode $S$}  (resp.~\emph{left antipode $S$}),
 i.e. there is a $K$-linear map $S:B\to B$ which is a morphism
of $C^3$-coalgebras $(B,\Delta,\epsilon,\mathbf{1})$ to itself such that
\begin{equation}
      \mathrm{id}* S =\mathbf{1}\epsilon~~~~~
      \big(\mathrm{resp.~~} S*\mathrm{id}=\mathbf{1}\epsilon\big)
\end{equation}
where $*$ denotes the convolution product (see Appendix \ref{AppCoalgebras} for 
definitions). It will become clear \emph{a posteriori} that right or left 
antipodes are always unique, see Lemma \ref{LRightAntipodeProperties}.\\
A first class of examples is of course the well-known class of
all \emph{cocommutative Hopf algebras} 
$(H,\Delta,\epsilon,\mathbf{1},\mu,S)$ for which $\mathbf{1}$ is a unit element,
and $S$ is a right and left antipode.\\
Secondly it is easy to check that every $C^4$-coalgebra $(C,\Delta,\epsilon,
\mathbf{1})$ equipped with
the left-trivial multiplication (resp.~right trivial multiplication) $\mu_0$ 
(see eqn (\ref{EqDefLeftTrivialMultiplication})) and 
\emph{trivial right antipode} (resp.~\emph{trivial left antipode}) $S_0$ 
defined
by $S_0(x)=\epsilon(x)\mathbf{1}$ for all $x\in C$ (in both cases) 
is a cocommutative right Hopf algebra (resp.~cocommutative left Hopf algebra)
called the cocommutative \emph{left-trivial right Hopf algebra} 
(resp. \emph{right-trivial left Hopf algebra}) defined by
the $C^4$-coalgebra $(C,\Delta,\epsilon,\mathbf{1})$.

We have the following elementary properties showing in particular that
each right (resp.~left) antipode is unique:
\begin{lemma}\label{LRightAntipodeProperties}
 Let $(\mathcal{H},\Delta,\epsilon,\mathbf{1},\mu,S)$ be a cocommutative right 
 Hopf algebra.
 \begin{enumerate}
  \item  \label{EqsCompElementaryPropertiesOfAntipodes}
     $S*(S\circ S)=\mathbf{1}\epsilon$,~~
        $S\circ S =\mathrm{id}*\mathbf{1}\epsilon$,
        $S*\mathbf{1}\epsilon=S$, and $S\circ S\circ S= S$,
    which for each $a\in \mathcal{H}$  implies 
    $\sum_{(a)}a^{(1)}S(a^{(2)})=\mathbf{1}\epsilon(a)
        = \sum_{(a)}S(a^{(1)})a^{(2)}\mathbf{1}$. It follows that right antipodes
        are unique.
  \item \label{EqCompAntipodeMultiplicativeAntihomomorphimsm}
       For all $a,b\in \mathcal{H}$: $S(ab)=S(b)S(a)$.
  \item \label{EqCompGenIdempotentsWithAntipode}
     For any element $c\in \mathcal{H}$, $c$ is a 
      generalized idempotent if and only if  $c=\sum_{(c)} S(c^{(1)})c^{(2)}$ iff
      there is $x\in \mathcal{H}$ with $c=\sum_{(x)} S(x^{(1)})x^{(2)}$, 
      and all these
      three statements imply that $c$ is a generalized left unit element.
 \end{enumerate}
\end{lemma}
\begin{prooof}
 {\it 1.} Since $S$ is a coalgebra morphism, it preserves convolutions
when composing from the right. This gives the first equation from statement
{\it \ref{EqsCompElementaryPropertiesOfAntipodes}}. Hence
the elements $\mathrm{id}$, $S$, and $S\circ S$ satisfy the hypotheses of
the elements $a,b,c$ of
Lemma \ref{LRightInversesHavingRightInverses} in the left-unital
convolution semigroup
$\big(\mathrm{Hom}_K(\mathcal{H},\mathcal{H}),*,\mathbf{1}\epsilon\big)$, 
whence the second and
third equations of statement
{\it \ref{EqsCompElementaryPropertiesOfAntipodes}}. 
are immediate, and the fourth follows from composing the 
second from the right with $S$ and using the third. Clearly $S$ is unique
according to Lemma \ref{LRightInversesHavingRightInverses}.\\
{\it 2.} Again the elements $\mu$, $S\circ \mu$ and $(\mathrm{id}* 1\epsilon)\circ \mu$
satisfy the hypotheses on the elements
$a,b,c$ of
Lemma \ref{LRightInversesHavingRightInverses}
in the left-unital convolution semigroup
$\big(\mathrm{Hom}_K(\mathcal{H}\otimes \mathcal{H},\mathcal{H}),*,
\mathbf{1}(\epsilon\otimes\epsilon)\big)$
(using the fact that $\mu$ is a morphism of coalgebras) whence $S\circ \mu$
is the unique right inverse of $\mu$. A computation shows that also
$\mu\circ \tau\circ (S\otimes S)$ is a right inverse of $\mu$, whence we get
statement {\it \ref{EqCompAntipodeMultiplicativeAntihomomorphimsm}}.
by uniqueness of right inverses (Lemma \ref{LRightInversesHavingRightInverses}).
\\
{\it 3.} The second statement obviously implies the third, and it is easy to 
see by straight-forward computations 
that the third statement implies the first and the second. Conversely, if
$c\in\mathcal{H}$ is a generalized idempotent, i.e. $c=(\mu\circ\Delta)(c)$, 
we get --since 
$\mu\circ\Delta$ is a morphism of 
$C^3$-coalgebras-- that
\begin{eqnarray*}
    \sum_{(c)}S(c^{(1)})c^{(2)}
   & = &\sum_{(c)}S(c^{(1)})\big((\mu\circ\Delta)(c)\big)^{(2)}
   =\sum_{(c)}S(c^{(1)})(\mu\circ\Delta)\big(c^{(2)}\big) \\
   & = & \sum_{(c)}S(c^{(1)})c^{(2)}c^{(3)}
     = \sum_{(c)}S(c^{(1)})c^{(2)}\mathbf{1}c^{(3)}
     \stackrel{\mathrm{Lemma}~\ref{LRightAntipodeProperties}, 1.}{=} c,
\end{eqnarray*}
and all the three statements are equivalent.
In order to see that every such element $c$ is a generalized left unit element
pick $y\in \mathcal{H}$ and
\[
    cy = \sum_{(x)} S(x^{(1)})x^{(2)}y
       = \sum_{(x)} S(x^{(1)})x^{(2)}\mathbf{1}y 
       \stackrel{\mathrm{Lemma}~\ref{LRightAntipodeProperties}, 1.}{=}
       \epsilon(x)y =\epsilon(c)y
\]
since obviously $\epsilon(c)=\epsilon(x)$, so $c$ is a generalized 
left unit element.
\end{prooof}

There is the following right Hopf algebra analogue of the 
\emph{Suschkewitsch decomposition theorem}
for right groups (see Appendix \ref{AppSemigroups}):
\begin{theorem}\label{TRightHopfAlgebrasSuschkewitschDecomposition}
 Let $(\mathcal{H},\Delta,\epsilon,\mathbf{1},\mu,S)$ a 
 cocommutative right Hopf algebra.
 Then the following holds:
 \begin{enumerate}
  \item The $K$-submodule $(\mathcal{H}\mathbf{1},\Delta|_{\mathcal{H}\mathbf{1}},
            \epsilon|_{\mathcal{H}\mathbf{1}},
            \mu_{\mathcal{H}\mathbf{1}\otimes \mathcal{H}\mathbf{1}},
              S|_{\mathcal{H}\mathbf{1}})$
       is a unital Hopf subalgebra of 
       $(\mathcal{H},\Delta,\epsilon,\mathbf{1},\mu,S)$.
  \item The $K$-submodule 
       $E_\mathcal{H}:=\{x\in \mathcal{H}~|~x~\mathrm{is~a~generalized~idempotent}\}$
      is a right Hopf subalgebra of $\mathcal{H}$ equal to the left-trivial right Hopf 
      algebra defined by the $C^4$-coalgebra ($E_\mathcal{H},
      \Delta|_{E_\mathcal{H}},
      \epsilon|_{E_\mathcal{H}},\mathbf{1}\epsilon|_{E_\mathcal{H}})$.
  \item The map 
     \[
         \Psi:\mathcal{H}\to \mathcal{H}\mathbf{1}\otimes E_\mathcal{H}: x\mapsto 
           \sum_{(x)}x^{(1)}\mathbf{1}\otimes S(x^{(2)})x^{(3)}
     \]
     is an isomorphism of right Hopf algebras whose inverse $\Psi^{-1}$
     is the restriction
     of the multiplication map.
 \end{enumerate}
\end{theorem}
\begin{prooof}
 {\it 1.} It is easy to see that $\mathcal{H}\mathbf{1}$ equipped with all the 
   restrictions is a unital bialgebra. Note that for all $a\in \mathcal{H}$
   \begin{eqnarray*}
   \big(S|_{\mathcal{H}\mathbf{1}}* \mathrm{id}_{\mathcal{H}\mathbf{1}})\big)
       (a\mathbf{1})
        &=& \Big(\big(S*\mathrm{id}_{\mathcal{H}}\big)(a)\Big)\mathbf{1}
       =\big(S*\mathrm{id}_{\mathcal{H}}* (\mathbf{1}\epsilon)\big)(a)\\
     &  \stackrel{\mathrm{Lemma}~\ref{LRightAntipodeProperties}}{=} &
       (\mathbf{1}\epsilon)(a)= (\mathbf{1}\epsilon)(a\mathbf{1})
       =(\mathbf{1}\epsilon|_{\mathcal{H}\mathbf{1}})(a\mathbf{1})     
   \end{eqnarray*}
   whence $S|_{\mathcal{H}\mathbf{1}}$ is also a left antipode. It follows that 
   $\mathcal{H}\mathbf{1}$ is a Hopf algebra.\\
   {\it 2.} Since the property of being an generalized idempotent
   is a $K$-linear condition, it follows that $E_\mathcal{H}$ is a $K$-submodule.
   Moreover since each $c\in E_\mathcal{H}$ is of the general form 
   $c=\sum_{(x)}S(x^{(1)})x^{(2)}$, $x\in \mathcal{H}$, and since the map 
   $\iota:\mathcal{H}\to \mathcal{H}$
   defined by $\iota(x)=\sum_{(x)}S(x^{(1)})x^{(2)}$ is an idempotent
   morphism of $C^3$-coalgebras,
   we get
   \[
       \Delta(c)= \Delta\big(\iota(c)\big)
        =(\iota\otimes \iota)\big(\Delta(x)\big)
   \]
   showing that $E_\mathcal{H}$ is a $C^3$-subcoalgebra of $\mathcal{H}$. 
   Furthermore, since every
   element of $E_\mathcal{H}$ is a generalized left unit element 
   (Lemma \ref{LRightAntipodeProperties}, {\it 3.}), the restriction of the 
   multiplication of $\mu$ of $\mathcal{H}$ to 
   $E_\mathcal{H}\otimes E_\mathcal{H}$ is left trivial. Finally,
   \begin{eqnarray*}
      S(c)=S\big(\iota(c)\big) &= &\sum_{(c)}S\big(S(c^{(1)})c^{(2)}\big)
      \stackrel{\mathrm{Lemma}~\ref{LRightAntipodeProperties}, ~2.}{=}
      \sum_{(c)}S(c^{(1)})S\big(S(c^{(2)})\big)\\
      & \stackrel{\mathrm{Lemma}~\ref{LRightAntipodeProperties},~ 1.}{=} &
      \sum_{(c)}S(c^{(1)})c^{(2)}\mathbf{1}
      \stackrel{\mathrm{Lemma}~\ref{LRightAntipodeProperties},~ 1.}{=}
        \epsilon(c)\mathbf{1} =S_0(c),
   \end{eqnarray*}
    showing that the the restriction of $S$ to $E_\mathcal{H}$ is the trivial 
    right 
    antipode.\\
    {\it 3.} It is clear from the two preceding 
    statements that $\Psi$ is a well-defined linear map into the tensor
    product of two cocommutative right Hopf algebras. We have for all 
    $x\in \mathcal{H}$
    \[
       (\mu\circ \Psi)(x)= \sum_{(x)}x^{(1)}\mathbf{1}S(x^{(2)})x^{(3)}
       \stackrel{\mathrm{Lemma}~\ref{LRightAntipodeProperties},~ 1.}{=}
       \sum_{(x)}\epsilon(x^{(1)})x^{(2)} = x
    \]
    and for all $a\in \mathcal{H}$, $c\in E_\mathcal{H}$ the term 
    $(\Psi\circ \mu)(a\mathbf{1}\otimes c)$
    is equal to
    \begin{eqnarray*}
      \lefteqn{ \sum_{(a)(c)}\big( a^{(1)}c^{(1)}\mathbf{1}\big)\otimes
               \big(S(c^{(2)})S(a^{(2)})a^{(3)}c^{(3)}\big)}\\
        & =&
         \sum_{(c)} 
         \big(a c^{(1)}\mathbf{1}\big)\otimes \big(S(c^{(2)})c^{(3)}\big)
         =  \sum_{(c)} 
         (a \mathbf{1})\otimes \big(S(c^{(1)})c^{(2)}\big)
         = (a \mathbf{1})\otimes c
    \end{eqnarray*}
    because all the terms $S(a^{(2)})a^{(3)}$ and the components $c^{(1)},\ldots$
    of iterated comultiplications of generalized idempotents can be chosen
    in $E_\mathcal{H}$ (since the latter has been shown to be a subcoalgebra), and
    are thus generalized left unit elements 
    (Lemma \ref{LRightAntipodeProperties}, {\it 3.}). Hence
    $\Psi$ is a $K$-linear isomorphism. Moreover, it is easy to
    see from its definition that $\Psi$ is a morphism of $C^3$-coalgebras.\\
    Next we compute for all $a,a'\in \mathcal{H}$ and $c,c'\in E_H$:
    \begin{eqnarray*}
      \Psi^{-1}\Big(\big((a\mathbf{1})\otimes c\big)
                     \big((a'\mathbf{1})\otimes c'\big)\Big)
                     & = &
       \Psi^{-1}\big((a\mathbf{1}a'\mathbf{1})\otimes \epsilon(c)c'\big)
          = \epsilon(c)aa'c',
    \end{eqnarray*}
    and -since $c$ is a generalized left unit element--
    \begin{eqnarray*}
     \Psi^{-1}\big((a\mathbf{1})\otimes c\big)
       \Psi^{-1}\big((a'\mathbf{1})\otimes c'\big)
       & = & acac'= \epsilon(c)aa'c',
    \end{eqnarray*}
    showing that $\Psi^{-1}$ and hence $\Psi$ is a morphism of left-unital 
    algebras. Finally we obtain
    \begin{eqnarray*}
       (S|_{\mathcal{H}\mathbf{1}}\otimes S_0)\big(\Psi(x)\big)
          &= &\sum_{(x)} \big(S(x^{(1)})\mathbf{1}\big)
             \otimes \big(\epsilon\big(S(x^{(2)})x^{(3)}\big)\mathbf{1}\big)
             =S(x)\mathbf{1}\otimes \mathbf{1}, \\
      \Psi\big(S(x)\big)     & = & 
      \sum_{(x)}\big(S(x^{(1)})\mathbf{1}\big)\otimes
              \big(S(S(x^{(2)}))S(x^{(3)})\big)=
              S(x)\mathbf{1}\otimes \mathbf{1},
    \end{eqnarray*}
   thanks to Lemma \ref{LRightAntipodeProperties}, and $\Psi$ intertwines
   right antipodes.
\end{prooof}

Note that the $K$-submodule of all \emph{generalized 
left unit elements} of a right Hopf algebra
$\mathcal{H}$ is given by 
$K\mathbf{1}\oplus\big(\Phi^{-1}(\mathcal{H}\mathbf{1}\otimes E_{\mathcal{H}}^+)
\big)$ and thus in general much bigger than the submodule $E_{\mathcal{H}}$
of all generalized idempotents.

As it is easy to see that every tensor product $H\otimes C$  of a unital 
cocommutative Hopf algebra
$H$ and a $C^4$-coalgebra $C$ (equipped with the left-trivial multiplication
and the trivial right antipode) is a right Hopf algebra, it is a fairly routine 
check --using the preceding Theorem
\ref{TRightHopfAlgebrasSuschkewitschDecomposition}-- that the category of all
cocommutative right Hopf algebras is equivalent to the product category of
all cocommutative Hopf algebras and of all $C^4$-coalgebras. \\
In the sequel, we 
shall need the dual \emph{left Hopf algebra version} where all the formulas 
have to
be put in reverse order: Here every left Hopf algebra is isomorphic to
$C\otimes H$.

\subsubsection{Dialgebras and Rack Bialgebras}

Recall (cf.~e.g.~\cite{Lod2001}) that a \emph{dialgebra} over $K$ is a 
$K$-module $D$ equipped with two
associative multiplications $\vdash,\dashv:A\otimes A\to A$ (written
$a\otimes b\mapsto a\vdash b$ and $a\otimes b\mapsto a\dashv b$)
satisfying for all $a,b,c\in A$:
\begin{eqnarray}
    (a\vdash b)\vdash c & = & (a\dashv b) \vdash c,
        \label{EqLeftVdashEqualLeftDashvDialgebra}\\
    a \dashv (b\dashv c) & = & a \dashv (b\vdash c), 
       \label{EqRightVdashEqualRightDashvDialgebra}\\
    (a\vdash b) \dashv c & = & a \vdash (b \dashv c)
      \label{EqLeftVdashCommutesWithRightDashvDialgebra}.
\end{eqnarray}
An element $\mathbf{1}$ of $A$ is called a \emph{bar-unit element}
of the dialgebra $(A,\vdash,\dashv)$ and $(A,\mathbf{1},\vdash,\dashv)$ is called
a \emph{bar-unital dialgebra} iff in addition the following holds
\begin{eqnarray}
    \mathbf{1} \vdash a & = & a, \label{EqDefBarUnitDialgebraLeft}\\
    a \dashv \mathbf{1} & = & a, \label{EqDefBarUnitDialgebraRight}
\end{eqnarray}
for all $a\in A$. Moreover, we shall call 
a bar-unital dialgebra $(A,\mathbf{1},\vdash,\dashv)$ \emph{balanced} iff in 
addition for all $a\in A$
\begin{equation} \label{EqDefBalancedDialgebra}
    a\vdash \mathbf{1} = \mathbf{1} \dashv a.
\end{equation}
Clearly each associative algebra is a dialgebra upon setting $\vdash=\dashv$
equal to the given multiplication.
The class of all (bar-unital and balanced) dialgebras forms a category
where morphisms preserve both multiplications and map the initial bar-unit to
the target bar-unit.

These algebras had been introduced to have a sort of `associative analogue'
for Leibniz algebras. More precisely, there is the following important fact,
which can easily be checked, see e.g. \cite{Lod2001}:
\begin{prop}\label{PEveryDialgebraIsLeibniz}
  Let $(A,\vdash,\dashv)$ be a dialgebra. Then the $K$-module $A$ equipped
  with the bracket $[~,~]:A\otimes A\to A$, written $[a,b]$,
  \begin{equation}\label{EqDefLeibnizBracketForDialgebras}
        [a,b]:=a\vdash b - b \dashv a
  \end{equation}
  is a Leibniz algebra, denoted by $A^-$.
\end{prop}
In fact, this construction is well-known to give rise to a functor 
$A\to A^-$ from the 
category of all 
dialgebras to the category of all Leibniz algebras in complete analogy to the 
obvious functor from the category of all associative algebras to the category
of all Lie algebras.

An important construction of (bar-unital) dialgebras is the following:

\begin{exem} \label{construction_dialgebra}
Let $(B,\mathbf{1}_B)$
be a unital associative algebra over $K$, and let $A$ be a $K$-module which 
is a 
$B$-\emph{bimodule}, i.e. there are $K$-linear maps $B\otimes A\to A$ and 
$A\otimes B\to A$
(written $(b\otimes x)\mapsto bx$ and $(x\otimes b)\mapsto xb$)
equipping $A$ with the structure of a left $B$-module and a right $B$-module such
that $(bx)b'=b(xb')$ for all $b,b'\in B$ and for all $x\in A$. Suppose 
in addition that
there is a bimodule map $\Phi:A\to B$, i.e. $\Phi(bxb')=b\Phi(x)b'$
for all $b,b'\in B$ and for all $x\in A$. Then it is not hard to check
that the two multiplications $\vdash,\dashv:A\otimes A\to A$ defined by
\begin{equation}
   x\vdash y := \Phi(x)y~~~\mathrm{and}~~~x\dashv y := x\Phi(y)
\end{equation}
equip $A$ with the structure of a dialgebra. If in addition there is an element
$\mathbf{1}\in A$ such that $\Phi(\mathbf{1})=\mathbf{1}_B$, then
$(A,\mathbf{1},\vdash,\dashv)$ will be a bar-unital dialgebra. We shall call
this structure $(A,\Phi,B)$ an \emph{augmented dialgebra}. 
\end{exem}

In fact, every
dialgebra $(A,\vdash, \dashv)$ arises in that fashion: Consider the 
$K$-submodule $I\subset A$ whose elements are linear combinations of
arbitrary product expressions
\[
      p\big(a_1,\ldots, a_{r-1}, (a_r\vdash b_r-a_r\dashv b_r), a_{r+1},\ldots, 
      a_n\big)
\]
(where all reasonable parentheses and symbols $\vdash$ and $\dashv$ are 
allowed)
for any two strictly positive integer $r\leq n$, and $a_,\ldots,a_n, b_r\in A$.
It follows that the quotient module $A/I$ is equipped with an associative
multiplication induced by both $\vdash$ and $\dashv$. Let $A_{\mathrm{ass}}^1$
be equal to $A/I$ if $A$ is bar-unital: In that case, the
bar-unit $\mathbf{1}$
of $A$ projects on
the unit element of $A/I$; and let $A_{\mathrm{ass}}^1$ be equal to 
$A/I\oplus K$ (adjoining a unit element) in case $A$ does not have a bar-unit. 
Thanks to the defining equations (\ref{EqLeftVdashEqualLeftDashvDialgebra}),
(\ref{EqRightVdashEqualRightDashvDialgebra}), 
(\ref{EqLeftVdashCommutesWithRightDashvDialgebra}), it can be shown by induction
that for any strictly positive integer $n$, any $a_1,\ldots,a_n,a\in A$, and any 
product expression made of the preceding elements upon using $\vdash$ or
$\dashv$
\begin{eqnarray*}
   p(a_1,\ldots,a_n)\vdash a  & =  & a_1\vdash \cdots \vdash a_n \vdash a
     =(a_1\dashv \cdots \dashv a_n) \vdash a, \\
   a \dashv p(a_1,\ldots,a_n) & = & a \dashv  a_1\dashv \cdots \dashv a_n
   = a \dashv (a_1\vdash \cdots \vdash a_n),
\end{eqnarray*}
proving in particular that $I$ acts trivially from the left 
(via $\vdash$) and from the right (via $\dashv$) on $A$ such 
that there is a well-defined $A_{\mathrm{ass}}^1$-bimodule structure
on $A$ such that the natural map $\Phi_A:A\to A_{\mathrm{ass}}^1$ is a bimodule
morphism. Hence $(A,\Phi_A, A_{\mathrm{ass}}^1)$ is always
an augmented dialgebra, and the assignment 
$A\to (A,\Phi_A, A_{\mathrm{ass}}^1)$ is known to be a faithful functor.\\
Note also that this construction allows to \emph{adjoin a bar-unit} to a
dialgebra $(A,\vdash,\dashv)$: Consider the $K$-module 
$\tilde{A}:=A\oplus A_{\mathrm{ass}}^1$
with the obvious $A_{\mathrm{ass}}^1$-bimodule structure 
$\alpha.(b+\beta)=\alpha.b+\alpha\beta$ and $(b+\beta).\alpha=b.\alpha+\beta\alpha$ 
for all $\alpha,\beta\in A_{\mathrm{ass}}^1$ and $b\in A$. Observe that
the obvious map $\Phi_{\tilde{A}}:\tilde{A}\to A_{\mathrm{ass}}^1$ defined
by $\Phi_{\tilde{A}}(b+\beta)=\Phi_A(b)+\beta$ is an
$A_{\mathrm{ass}}^1$-bimodule map, and that $\mathbf{1}=\mathbf{1}_{A_{\mathrm{ass}}^1}$
is a bar-unit. The bar-unital augmented dialgebra 
$(\tilde{A},\Phi_{\tilde{A}},A_{\mathrm{ass}}^1)$
is easily seen to be balanced. There are nonbalanced bar-unital dialgebras as can be seen
from the augmented bar-unital dialgebra example 
$(B\otimes B,\mathbf{1}_B\otimes \mathbf{1}_B,\mu_B,B)$ where $(B,\mathbf{1},\mu_B)$
is any unital associative algebra and the bimodule action is defined by 
$b.(b_1\otimes b_2).b':= (bb_1)\otimes (b_2b')$ for all $b,b',b_1,b_2\in B$.
\\
Again, in case the dialgebra $(A,\vdash,\dashv,\mathbf{1})$ is
bar-unital and balanced,
note that $A\vdash \mathbf{1}=\mathbf{1}\dashv A$ is an associative
unital subalgebra $A'$ of $A$ whose multiplication is induced by both $\vdash$ and
$\dashv$, i.e.~$a'\vdash b'=a'\dashv b'$ for all $a',b'\in A'$. 
Since the $K$-linear map 
$\pi_A:A\to A:a\mapsto a\vdash \mathbf{1}=\mathbf{1} \dashv a$ descends
to a surjective morphism of 
associative algebras $A_{\mathrm{ass}}^1\to A'$ by the above, it is clear
that the ideal $I$ contains the kernel of $\pi_A$. On the other hand,
if $a\in \mathrm{Ker}(\pi_A)$ then $0=\pi(a)=\mathbf{1}\vdash a$, and 
obviously $a= \mathbf{1}\vdash a
-\mathbf{1}\dashv a\in I$, thus inducing a useful isomorphism 
$A_{\mathrm{ass}}^1\cong A'$, and thus a subalgebra injection 
$i_A:A_{\mathrm{ass}}^1\to A:\Phi_A(a)\mapsto a\vdash \mathbf{1}$ which is a
right inverse to the projection $\Phi_A$, i.e.~
$\Phi_A\circ i_A=\mathrm{id}_{A^1_{\mathrm{ass}}}$.

In this work, we also have to take into account coalgebra structures and thus define
the following:
\begin{defi}
 Let $(A,\Delta,\epsilon,\mathbf{1})$ be cocommutative $C^3$-coalgebra (a $C^4$-coalgebra)
 and two $K$-linear maps $\vdash,\dashv:A\otimes A\to A$. Then
  $(A,\Delta,\epsilon,\mathbf{1},\vdash,\dashv)$ will be called a \textbf{cocommutative
  bar-unital di-coalgebra} if and only if
  \begin{enumerate}
   \item $(A,\mathbf{1},\vdash,\dashv)$ is a bar-unital balanced dialgebra.
   \item Both $\vdash$ and $\dashv$ are morphisms of $C^3$-coalgebras.
  \end{enumerate}
  If in addition there is a morphism of $C^3$-coalgebras $S:A\to A$ such that
  $(A,\Delta,\epsilon,\mathbf{1},\vdash,S)$ is a cocommutative right Hopf algebra
  and $(A,\Delta,\epsilon,\mathbf{1},\dashv,S)$ is a cocommutative left Hopf algebra,
  then $(A,\Delta,\epsilon,\mathbf{1},\vdash,\dashv,S)$ is called a 
  \textbf{cocommutative
  Hopf dialgebra}.
\end{defi}
We have used a relatively simple notion of one single compatible coalgebra structure
motivated from differential geometry, see Section 
\ref{SecCoalgebraStructuresForPointedManifoldsETC}. In contrast to that, F. Goichot 
uses two a priori different coalgebra structures, see
\cite{Goi2001}. Moreover, a slightly more general context would have been to demand
 the existence of two different antipodes, a right antipode $S$ for $\vdash$,
 and a left antipode $S'$ for $\dashv$. The theory --including the classification
 in terms of ordinary Hopf algebras-- could have been done as well, but we have
 refrained from doing so since it is not hard to see that such a more general
 Hopf dialgebra is \textbf{balanced iff $S=S'$}. This fact is crucial in the
 following refinement of Proposition \ref{PEveryDialgebraIsLeibniz}:
\begin{prop}\label{PPrimOfEveryCocomHopfDialgebraIsLeibniz}
 Let $(A,\Delta,\epsilon,\mathbf{1},\vdash,\dashv,S)$ be cocommutative 
 Hopf dialgebra. Then 
  the submodule of all primitive elements of $A$,
 $\mathsf{Prim}(A)$, is a Leibniz subalgebra of $A$ equipped with 
 the bracket (\ref{EqDefLeibnizBracketForDialgebras}).
\end{prop}
\begin{prooof}
 Let $x,y\in A$ be primitive. Then, using that $\vdash$ and $\dashv$ 
are morphisms of coalgebras, we get 
 \begin{eqnarray*}
   \Delta(x\vdash y-y\dashv x)
  & = & \mathbf{1}\otimes (x\vdash y-y\dashv x) 
          + (x\vdash y-y\dashv x)\otimes \mathbf{1} \\
    & & +  (x\vdash \mathbf{1})\otimes y + y\otimes (x\vdash \mathbf{1})
          -y \otimes (\mathbf{1}\dashv x) -(\mathbf{1}\dashv x)\otimes y \\
       & = & \mathbf{1}\otimes (x\vdash y-y\dashv x)  
                + (x\vdash y-y\dashv x)\otimes \mathbf{1} 
 \end{eqnarray*}
 \textbf{because $A$ is balanced}, and therefore $x\vdash y-y\dashv x$ is 
 primitive.
\end{prooof}

The first relationship with rack bialgebras is the following simple 
generalization of a cocommutative Hopf algebra  equipped with the adjoint
representation:
\begin{prop}\label{PHopfDialgebraIsRackBialgebra}
  Let $(A,\Delta,\epsilon,\mathbf{1},\vdash,\dashv,S)$ be cocommutative 
 Hopf dialgebra. Define the following multiplication $\mu:A\otimes A\to A$
 by
  \begin{equation}
     \mu(a\otimes b):=a\triangleright b:=
        \sum_{(a)} (a^{(1)} \vdash b) \dashv \big(S(a^{(2)})\big).
  \end{equation}
  Then we have the following:
  \begin{enumerate}
  \item The map $\rhd$ defines on the $K$-module $A$ two left module structures,
  one with respect to the algebra $(A,\mathbf{1},\vdash)$, and one with respect 
  to the algebra $(A,\mathbf{1},\dashv)$, making the Hopf-dialgebra 
  $(A,\Delta,\epsilon,\mathbf{1},\vdash,\dashv,S)$
  a module-Hopf dialgebra, i.e.
  \begin{eqnarray}
     a\triangleright (b\triangleright c) & =& (a\vdash b) \triangleright c
                                                =(a\dashv b) \triangleright c 
                                       \label{EqDefModuleIdentityTriHopfDialgebra}
                                           \\
     \Delta(a\triangleright b) & =& 
     \sum_{(a),(b)} (a^{(1)} \triangleright b^{(1)})\otimes
                      (a^{(2)} \triangleright b^{(2)}) 
                   \label{EqDefTriComultiplicationHopfDialgebra}   \\
       a\triangleright (b\vdash c) & =& \sum_{(a)} (a^{(1)} \triangleright b)
                                     \vdash
                                     (a^{(2)} \triangleright c) 
                      \label{EqDefTriVdashModuleAlgebraHopfDialgebra}   \\
     a\triangleright (b\dashv c) & =& \sum_{(a)} (a^{(1)} \triangleright b)
                                     \dashv
                                     (a^{(2)} \triangleright c)
                    \label{EqDefTriDashvModuleAlgebraHopfDialgebra}   \\                 
  \end{eqnarray}
  \item $(A,\Delta,\epsilon,\mathbf{1}, \mu)$ is a cocommutative rack bialgebra.
  \end{enumerate}
\end{prop}
\begin{prooof} 
  {\it 1.} First of all we have
  \[
      \mu= \mu_\dashv\circ (\mu_\vdash \otimes S)\circ 
      (\mathrm{id}_A \otimes \tau_{A,A})\circ (\Delta \otimes \mathrm{id}_A)
  \]
  where $\mu_\vdash$ and $\mu_\dashv$ stand for the multiplication maps
  $\vdash$ and $\dashv$, and this is clearly a composition of morphisms
  of $C^3$-coalgebras whence $\mu$ is a morphism of  $C^3$-coalgebras
  proving eqn (\ref{EqDefTriComultiplicationHopfDialgebra}).
  Next, there is clearly $\mathbf{1}\triangleright b=b$ for all $b\in A$, 
  and,
   \textbf{since the dialgebra is balanced}, we get for all $a\in A$
   \begin{eqnarray*}
      a\triangleright\mathbf{1} &= &
      \sum_{(a)} (a^{(1)} \vdash \mathbf{1}) \dashv \big(S(a^{(2)})\big)
      =\sum_{(a)}\mathbf{1}\dashv a^{(1)} \dashv \big(S(a^{(2)})\big) \\ 
      & = & \sum_{(a)}\mathbf{1}\dashv 
         \Big(a^{(1)} \vdash \big(S(a^{(2)})\big)\Big)
           =\epsilon(a)\mathbf{1}\dashv \mathbf{1}  =\epsilon(a)\mathbf{1}.
   \end{eqnarray*}
   Next, let $a,b,c\in A$. Then
   \begin{eqnarray}
     a\triangleright (b\triangleright c) 
     & =& 
     \sum_{(a),(b)}
     \Big(a^{(1)}\vdash \big((b^{(1)}\vdash c)\dashv S(b^{(2)})\big)\Big)
            \dashv S(a^{(2)}) \nonumber\\
      & = & \sum_{(a),(b)} \Big(\big((a^{(1)}\vdash b^{(1)})\vdash c\big)
                     \dashv\big(S(b^{(2)})\dashv S(a^{(2)})\big) \Big)
                     \nonumber\\
       & = & \sum_{(a),(b)} \Big(\big((a^{(1)}\vdash b^{(1)})\vdash c\big)
                     \dashv\big(S(b^{(2)})\vdash S(a^{(2)})\big) \Big)
                     \nonumber\\ 
       & = & \sum_{(a),(b)} \Big(\big((a^{(1)}
                            {~}^\dashv_\vdash~
       b^{(1)})
               \vdash c\big)
                     \dashv\big(S(a^{(2)}
  {~}^\dashv_\vdash~
                     b^{(2)})\big) \Big)
                     \nonumber\\ 
      & = &  (a
       {~}^\dashv_\vdash~
      b) \triangleright c, 
   \end{eqnarray}
   proving eqs (\ref{EqDefModuleIdentityTriHopfDialgebra}). Next, for all
   $a,a',a''\in A$, we get
   \begin{eqnarray*}
    \sum_{(a)}
    (a^{(1)}\triangleright a'){~}^\dashv_\vdash~ (a^{(2)}\triangleright a'') 
      & = &
     \sum_{(a)} \Big((a^{(1)} \vdash a') \dashv S(a^{(2)})\Big)
                            {~}^\dashv_\vdash~
                  \Big((a^{(3)} \vdash a'') \dashv S(a^{(4)})\Big)  \\
     & = & 
     \sum_{(a)}
     (a^{(1)} \vdash a'){~}^\dashv_\vdash~ 
      \Big(\big(S(a^{(2)}){~}^\dashv_\vdash~ a^{(3)}\big)
        {~}^\dashv_\vdash~ \big(a''\dashv S(a^{(4)})\big)\Big)\\
       & = & 
     \sum_{(a)}
     (a^{(1)} \vdash a'){~}^\dashv_\vdash~ \big(\big(\epsilon(a^{(2)})\mathbf{1}\big)
        {~}^\dashv_\vdash~ \big(a''\dashv S(a^{(3)})\big)\big)\\  
        & = & 
     \sum_{(a)} 
        \big(a^{(1)} \vdash (a' {~}^\dashv_\vdash~ a'')\big) \dashv S(a^{(2)})
        = a\triangleright(a' {~}^\dashv_\vdash~ a''),
 \end{eqnarray*}
 where --in the second to last equation-- we have used 
 the left antipode identity for the case $\dashv$ and the fact that
 $\sum_{(a)}S(a^{(1)})\vdash a^{(2)}$ is a generalized left unit element for the 
 case $\vdash$. It follows that $(A,{~}^\dashv_\vdash~)$ is an $A$-module-algebra
 proving eqs (\ref{EqDefTriVdashModuleAlgebraHopfDialgebra}) and
 (\ref{EqDefTriDashvModuleAlgebraHopfDialgebra}).\\
{\it 2.} It remains to prove self-distributivity: For all $a,b,c\in A$,
  we get
   \begin{eqnarray*}
      \sum_{(a)} (a^{(1)}\triangleright b) \triangleright
                  (a^{(2)}\triangleright c)  
                  & \stackrel{(\ref{EqDefModuleIdentityTriHopfDialgebra})}{=} &
                  \sum_{(a)} 
                  \big((a^{(1)}\triangleright b)\dashv a^{(2)}\big)
                        \triangleright c,
   \end{eqnarray*}
   and in the end
   \begin{eqnarray*}
    \sum_{(a)} (a^{(1)}\triangleright b)\dashv a^{(2)}
    & = &
      \sum_{(a)}
         \big( (a^{(1)}\vdash b)\dashv S(a^{(2)})\big)\dashv a^{(3)} \\
    & =  &
        \sum_{(a)}
          (a^{(1)}\vdash b)\dashv \big( S(a^{(2)})\dashv  a^{(3)}\big) \\
     & =  &
        \sum_{(a)}
          (a^{(1)}\vdash b)\dashv \big(\epsilon(a^{(2)})\mathbf{1} \big)
          = a\vdash b
   \end{eqnarray*}
   proving the self-distributivity identity.
\end{prooof}

The next theorem relates augmented cocommutative rack bialgebras with cocommutative
Hopf dialgebras:
\begin{theorem}\label{TAugmentedRackBialgebraOtimesHopfIsHopfDialgebra}
 Let $(B,\Phi_B,H,\ell)$ be a cocommutative augmented rack bialgebra. 
 Then the $K$-module
 $(B\otimes H,\Delta_{B\otimes H},\epsilon_B\otimes \epsilon_H,
    \mathbf{1}_B\otimes \mathbf{1}_H,\Phi,H)$ will be an augmented cocommutative
    Hopf dialgebra by means of the following definitions. Here we use Example 
    \ref{construction_dialgebra} and take $h,h'\in H$ and 
    $b\in B$:
    \begin{enumerate}
     \item \label{EqDefPhiCOtimesH}
       $\Phi:B\otimes H\to H:(b\otimes h)\mapsto \Phi(b\otimes h):=\Phi_B(b)h$.
     \item \label{EqDefBimoduleStructureForCOtimesH}
     $h'.(b\otimes h):=\sum_{(h')}((h')^{(1)}.b)\otimes ((h')^{(2)}h)$ and
     $(b\otimes h).h':= b\otimes (hh')$.
     \item \label{EqDefDiAntipodeForCOtimesH}
       $S(b\otimes h):= \mathbf{1}_B\otimes S_H\big(\Phi_B(b)h\big)$.
    \end{enumerate}
     Moreover, the Leibniz bracket on the $K$-module of all primitive elements of
       $B\otimes H$, $\mathfrak{a}:=\mathsf{Prim}(B)\oplus \mathsf{Prim}(H)$, 
       is computed
       as follows for all $x,y\in\mathsf{Prim}(B)$ and all 
       $\xi,\eta$ in the Lie algebra $\mathsf{Prim}(H)$
       (writing $x$ and $\xi$ for the more precise $x\otimes \mathbf{1}_H$
       and $\mathbf{1}_B\otimes \xi$)
       \begin{equation}\label{EqCompLeibnizBracketPrimitivesComposedDialgebra}
          [x+\xi,y +\eta] = \Big([x,y] + \xi.y\Big) + 
          \Big([\Phi_B(x),\eta]+ [\xi,\eta]\Big)
       \end{equation}
       where each bracket is of the form 
       (\ref{EqDefLeibnizBracketForDialgebras})\footnote{For an explicit formula, see the end of the proof of the theorem.}. Note that this Leibniz algebra
       is split over the Lie subalgebra $\mathsf{Prim}(H)$, the complementary
       two-sided ideal $\{x-\Phi_B(x)~|~x\in
       \mathsf{Prim}(B)\}$ being in the left center of $\mathfrak{a}$.
\end{theorem}
\begin{prooof}
 It is clear from the definitions that condition 
 {\it \ref{EqDefBimoduleStructureForCOtimesH}}
 defines a $H$-bimodule structure on $C\otimes H$ making it into a module $C^3$-coalgebra.
 Moreover, we compute for all $h,h',h''\in H$ and $b\in B$
 \begin{eqnarray*}
    \Phi\big(h'.(b\otimes h).h''\big) & = &
       \sum_{(h')} \Phi\big(((h')^{(1)}.b)\otimes ((h')^{(2)}hh'')
       =\sum_{(h')} \Phi_B((h')^{(1)}.b)(h')^{(2)}hh''\\
      & =&\sum_{(h')} \mathrm{ad}_{(h')^{(1)}}\big(\Phi_B(b)\big)(h')^{(2)}hh''
       =\sum_{(h')} (h')^{(1)}\Phi_B(b)S_H((h')^{(2)})(h')^{(3)}hh''\\
      & =& \sum_{(h')} (h')^{(1)}\Phi_B(b)\epsilon_H((h')^{(2)})hh'' 
         =  h'\Phi_B(b)hh'' =  h'\Phi(b\otimes h)h'',
 \end{eqnarray*}
 whence $\Phi$ is a morphism of $H$-bimodules. Next, we 
 get for all $b\in B$ and $h\in H$:
 \begin{eqnarray*}
  \big(\mathrm{id}_{B\otimes H} *_{\vdash} S\big)(b\otimes h)
     & = & \sum_{(b)(h)} (b^{(1)}\otimes h^{(1)})\vdash S(b^{(2)}\otimes h^{(2)}) \\
     & = &\sum_{(b)(h)} \Phi(b^{(1)}\otimes h^{(1)})
            .\Big(\mathbf{1}_B\otimes S_H\big(\Phi_B(b^{(2)})h^{(2)}\big)\Big) \\
     & = &  \sum_{(b)(h)} \big(\Phi_B(b^{(1)}) h^{(1)}\big)
                 .\Big(\mathbf{1}_B\otimes S_H\big(\Phi_B(b^{(2)})h^{(2)}\big)\Big) \\
     & = &  \sum_{(b)(h)} \epsilon_H\Big(\Phi_B(b^{(1)}) h^{(1)}\Big)\mathbf{1}_B
               \otimes 
             \Big(\Phi_B(b^{(2)}) h^{(2)} S_H\big(\Phi_B(b^{(3)})h^{(3)}\big)\Big) \\
     & = &   \mathbf{1}_B\otimes \big(\epsilon_B(b)\epsilon_H(h)\big)
                      \mathbf{1}_H 
                      =\big(\epsilon_B\otimes \epsilon_H\big)(b\otimes h)
                       \big(\mathbf{1}_B\otimes \mathbf{1}_H\big),
 \end{eqnarray*}
 proving the right antipode identity, and
 \begin{eqnarray*}
  \big(S *_{\dashv}\mathrm{id}_{B\otimes H}\big)(b\otimes h)
     & = & \sum_{(b)(h)} S(b^{(1)}\otimes h^{(1)})\dashv (b^{(2)}\otimes h^{(2)}) \\
     & = &\sum_{(b)(h)} 
            \Big(\mathbf{1}_B\otimes S_H\big(\Phi_B(b^{(1)})h^{(1)}\big)\Big). 
           \Phi(b^{(2)}\otimes h^{(2)}) \\
     & = &\sum_{(b)(h)} 
            \Big(\mathbf{1}_B\otimes S_H\big(\Phi_B(b^{(1)})h^{(1)}\big)\Big). 
           \big(\Phi_B(b^{(2)})h^{(2)}\big) \\  
       & = &\sum_{(b)(h)} 
             \mathbf{1}_B\otimes \Big(S_H\big(\Phi_B(b^{(1)})h^{(1)}\big)
                \Phi_B(b^{(2)})h^{(2)}\Big) \\               
     & = &   \mathbf{1}_B\otimes \big(\epsilon_B(b)\epsilon_H(h)\big)
                      \mathbf{1}_H =\big(\epsilon_B\otimes \epsilon_H\big)(b\otimes h)
                       \big(\mathbf{1}_B\otimes \mathbf{1}_H\big),
 \end{eqnarray*}
 proving the left antipode identity. Finally for all $h\in H$ we get
 \[
    h.(\mathbf{1}_B\otimes \mathbf{1}_H)
    =\sum_{(h)}(\epsilon_H(h^{(1)})\mathbf{1}_B)
      \otimes h^{(2)}
      =\mathbf{1}_B\otimes h
      =(\mathbf{1}_B\otimes \mathbf{1}_H).h
 \]
 implying that the bar-unital dialgebra is  balanced.\\
 Formula (\ref{EqCompLeibnizBracketPrimitivesComposedDialgebra}) is 
 straight-forward:
 \begin{eqnarray*}
  \lefteqn{[x\otimes \mathbf{1}_H +\mathbf{1}_B\otimes \xi, 
    y\otimes \mathbf{1}_H +\mathbf{1}_B\otimes \eta ]} \\
    & = &  (x\otimes \mathbf{1}_H)\vdash (y\otimes \mathbf{1}_H)
            -(y\otimes \mathbf{1}_H)\dashv (x\otimes \mathbf{1}_H) \\
    &  &   + (x\otimes \mathbf{1}_H)\vdash (\mathbf{1}_B\otimes \eta)
             - (\mathbf{1}_B\otimes \eta)\dashv 
                   (x\otimes \mathbf{1}_H)\\
    &  &    + (\mathbf{1}_B\otimes \xi)\vdash (y\otimes \mathbf{1}_H)
            - (y\otimes \mathbf{1}_H)\dashv (\mathbf{1}_B\otimes \xi) \\
    &   &   +(\mathbf{1}_B\otimes \xi)\vdash (\mathbf{1}_B\otimes \eta)
            - (\mathbf{1}_B\otimes \eta)\dashv (\mathbf{1}_B\otimes \xi)\\
   & = &  (\Phi_B(x).y)\otimes \mathbf{1}_H + y\otimes \Phi_B(x)
            - y\otimes \Phi_B(x) \\
    &  &   + \epsilon_B(x)\mathbf{1}_B\otimes \eta 
                   +\mathbf{1}_B \otimes \big(\Phi_B(x) \eta\big)
              - \mathbf{1}_B \otimes \big(\eta\Phi_B(x)\big) \\
    &  &   + (\xi.y)\otimes \mathbf{1}_H + y\otimes \xi 
                  -y\otimes \xi \\
    &   &   \epsilon_H(\xi)\mathbf{1}_B\otimes \eta 
           + \mathbf{1}_B\otimes (\xi\eta)
           -\mathbf{1}_B\otimes (\eta\xi) \\
    & = &  [x,y]\otimes \mathbf{1}_H 
              + \mathbf{1}_B\otimes \big[\Phi_B(x),\eta\big]
              + (\xi.y)\otimes \mathbf{1}_H  
              + \mathbf{1}_B\otimes [\xi,\eta],
 \end{eqnarray*}
because primitives are killed by counits, and the formula is proved.
 
\end{prooof}

\noindent A lengthy, but straight-forward reasoning shows that the above construction
assigning
$(B,\Phi_B,H,\ell)\to (B\otimes H, \Phi, H)$ defines a covariant functor from 
the category
of all cocommutative rack bialgebras to the category of all cocommutative Hopf 
dialgebras.\\
A particular case of the preceding theorem is obtained by picking any
$C^4$-coalgebra $(B,\Delta_B,\epsilon_B,\mathbf{1}_B)$ such that there is
any
$H$-module coalgebra structure $\ell$ on $B$ (such that 
$h.\mathbf{1}_B=\epsilon_H(h)\mathbf{1}_B$) and by choosing the trivial
map $\Phi_B(b)=
\epsilon_B(b)\mathbf{1}_H$: It follows that $(B,\Phi_B,H,\ell)$
is an augmented rack bialgebra with left-trivial multiplication.
It turns out that the Hopf dialgebra $B\otimes H$ formed out of this
is already isomorphic to the general cocommutative Hopf 
dialgebra:
\begin{theorem}
 Let $(A,\Delta,\epsilon,\mathbf{1},\vdash,\dashv,S)$ be a cocommutative 
 Hopf dialgebra. Let $E_A$ be the $C^3$-subcoalgebra of all generalized 
 idempotent elements\footnote{Recall that this means 
$c=(\dashv\circ\Delta)(c)$.} with respect to $\dashv$, and let 
$H_A=\mathbf{1}\dashv A$ 
 be the Hopf subalgebra according to the Suschkewitsch decomposition
 of the left Hopf algebra $(A,\Delta,\epsilon,\mathbf{1},\dashv,S)$,
 see Theorem \ref{TRightHopfAlgebrasSuschkewitschDecomposition}. Then we
 have:
 \begin{enumerate}
  \item By means of the Suschkewitsch isomorphism $A\to E_A\otimes H_A$
  for the left Hopf algebra $(A,\Delta,\epsilon,\mathbf{1},\dashv,S)$,
  we can transfer the cocommutative dialgebra structure of $A$ to
  $E_A\otimes H_A$: There is a well-defined left module-coalgebra action 
  $\ell$ of
  $H_A$ on $E_A$ defined by (for all $c\in E_H$, $h\in H_A$,
  $a\in A$ such that $h=\mathbf{1}\dashv a$) 
  \begin{equation}\label{EqDefLeftModuleOfHAOnIdempotents}
      \ell_h(c)=h.c=\sum_{(a)} (a^{(1)}\vdash c) \dashv S(a^{(2)}),
  \end{equation}
  and the transferred multiplications $\vdash'$ and $\dashv'$ and the 
  antipode $S'$ on 
  $E_A\otimes H_A$ read (for all $c,c'\in E_H$, $h,h'\in H_A$)
  \begin{eqnarray}
     (c\otimes h)\vdash' (c'\otimes h') & = &
         \epsilon(c)\sum_{(h)}(h^{(1)}.c')\otimes (h^{(2)}h'), 
           \label{EqCompVdashTransferred}\\
      (c\otimes h)\dashv' (c'\otimes h') & = &  
      \epsilon(c') c\otimes (hh'), 
         \label{EqCompDashvTransferred}\\
        S'(c\otimes h) & =& 
        \epsilon(c) \big(\mathbf{1}_{E_A}\otimes \big(S_{H_A}(h)\big).
         \label{EqCompAntipodeTransferred}
  \end{eqnarray}
  \item The $K$-linear map 
  $A\to \mathbf{1}\dashv A:a\mapsto \mathbf{1}\dashv a$ descends to an 
  isomorphism
  of associative algebras $A^1_{\mathrm{ass}}\cong H_A$.
  \item \textbf{S.Covez, 2006}: The Leibniz subalgebra $\mathsf{Prim}(A)$
   of $A$ (equipped with the Leibniz bracket 
   (\ref{EqDefLeibnizBracketForDialgebras}) is a \textbf{split semidirect sum} 
   out of the two-sided ideal $\mathsf{Prim}(E_A)\subset 
   \mathfrak{z}\big(\mathsf{Prim}(A)\big)$ and the Lie subalgebra
   $\mathsf{Prim}(H_A)$, i.e. for all $z,z'\in \mathsf{Prim}(E_A)$ and
   $\xi,\xi'\in\mathsf{Prim}(H_A)$, we have
   \begin{equation}\label{EqCompSplitLeibnizOnPrimitivesOfHopfDialgebra}
       [z+\xi,z'+\xi']=\xi.z + [\xi,\xi'].
   \end{equation}
   \item Let $(B,\Phi_B,H,\ell)$ be a cocommutative augmented rack bialgebra.
    Then for the Hopf dialgebra $B\otimes H$ of Theorem 
    \ref{TAugmentedRackBialgebraOtimesHopfIsHopfDialgebra}, we get that the Hopf 
    subalgebra $H_{B\otimes H}$ equals $\mathbf{1}_B \otimes H\cong H$, and
    \[
        E_{B\otimes H}
        =\left.\left\{\sum_{(b)}b^{(1)}\otimes \Big(S_H\big(\Phi_B(b^{(2)})\big)\Big)
             \in B\otimes H~\right|~b\in B \right\}
    \]
    which is isomorphic to $B$ as a $C^4$-coalgebra of $B\otimes H$.
 \end{enumerate}
\end{theorem}
\begin{prooof}
 {\it 1.} Note first that the right hand side of eqn 
 (\ref{EqDefLeftModuleOfHAOnIdempotents}) is just $a\triangleright c$
 of Proposition \ref{PHopfDialgebraIsRackBialgebra} which had been shown
 to be a left module-Hopf dialgebra action of $(A,\dashv)$ and of $(A,\vdash)$ 
 on $A$. 
 Observe that for all $a,a'\in A$
 \begin{eqnarray*}
       (\mathbf{1}\dashv a)\triangleright a' 
       & \stackrel{(\ref{EqDefModuleIdentityTriHopfDialgebra})}{=} &
          \mathbf{1}\triangleright (a \triangleright a') = 
          a \triangleright a'
 \end{eqnarray*}
 whence the $H_A$-action $\ell$ is well-defined on $A$. Moreover
 we compute 
 for all $h\in H_A$ and all $a,a',a''\in A$ such that $h=\mathbf{1}\dashv a$:
 \begin{eqnarray*}
    \sum_{(h)}(h^{(1)}.a')\dashv (h^{(2)}.a'') & = &
     \sum_{(a)} (a^{(1)} \triangleright a') 
                            \dashv
                  (a^{(2)} \triangleright a'')  \\
     & \stackrel{(\ref{EqDefTriDashvModuleAlgebraHopfDialgebra})}{=} & 
     a\triangleright (a'\dashv a'') = h.(a'\dashv a'')
 \end{eqnarray*}
 whence $(A,\dashv)$ is also a $H_A$-module-algebra.
 Now let $c\in E_H$. By definition, $c$ is a generalized idempotent 
 (w.r.t.~$\dashv$), hence
 $c= \sum_{(c)} c^{(1)}\dashv c^{(2)}$, and thus for all $h\in H$
 \[
     h.c=h.\sum_{(c)} c^{(1)}\dashv c^{(2)}
        =\sum_{(h),(c)} (h^{(1)}.c^{(1)})\dashv (h^{(2)}.c^{(2)})
        =\sum_{(h.c)} (h.c)^{(1)}\dashv (h.c)^{(2)}
 \]
whence $h.c$ is also in $E_A$, and $E_A$ is a $H_A$-submodule of $A$.\\
 Recall the Suschkewitsch decomposition of the left Hopf
 algebra $(A,\Delta,\epsilon,\mathbf{1},\vdash,S)$ where one can use Theorem 
 \ref{TRightHopfAlgebrasSuschkewitschDecomposition} and dualize
 all the formulas:
 \begin{eqnarray*}
       \Psi: A\to E_A\otimes H_A  & : & 
       a\mapsto \sum_{(a)}
       \big(a^{(1)}\dashv S(a^{(2)})\big)\otimes(\mathbf{1}\dashv a^{(3)}).
         \\
        \Psi^{-1}:E_A\otimes H_A \to A & : & 
          \big(c\otimes(\mathbf{1}\dashv a)\big) \mapsto c\dashv a.
 \end{eqnarray*}
 Formulas (\ref{EqCompDashvTransferred}) and (\ref{EqCompAntipodeTransferred})
 consequences of Theorem \ref{TRightHopfAlgebrasSuschkewitschDecomposition}.
 The only formula which remains to be shown is eqn 
 (\ref{EqCompVdashTransferred}). Note first that every generalized idempotent
 $c\in E_\mathcal{H}$ (w.r.t.~$\dashv$) is also a generalized idempotent with 
 respect
 to $\vdash$. Indeed, since all the components $c^{(1)}$ and $c^{(2)}$ in
 $\Delta(c)=\sum_{(c)}c^{(1)}\otimes c^{(2)}$ can be chosen in $E_\mathcal{H}$,
 we get
 \begin{eqnarray*}
  \sum_{(c)} c^{(1)}\vdash c^{(2)} 
    &=&\sum_{(c)} \big(c^{(1)}\dashv S(c^{(2)})\big)\vdash c^{(3)}
    =\sum_{(c)} \big(c^{(1)}\vdash S(c^{(2)})\big)\vdash c^{(3)}\\
    & = &\sum_{(c)} \epsilon(c^{(1)})\mathbf{1}\vdash c^{(2)} =c.
 \end{eqnarray*}
Next for
 all $c,c'\in E_\mathcal{H}$,
 $h,h'\in H_\mathcal{H}$, and $a,a'\in A$ such that $\mathbf{1}\dashv a=h$
 and $\mathbf{1}\dashv a'=h'$, we get --since $c$ is a generalized left unit 
 element (w.r.t.~$\vdash$) thanks to {\it 3.} in Lemma 
\ref{LRightAntipodeProperties}--
 \begin{eqnarray*}
   \lefteqn{ (c\otimes h)\vdash' (c'\otimes h')  = 
       \Psi\big(\Psi^{-1}(c\otimes h) \vdash\Psi^{-1}(c'\otimes h') \big)
    =\Psi\big((c\dashv a)\vdash (c'\dashv a')\big)} ~~~~~~~~~~~~~~~~~~~~~~
    ~~~~~~~\\
    & = & \Psi\big((c\vdash a)\vdash (c'\dashv a')\big) =
         \epsilon(c)\Psi\big(a \vdash (c'\dashv a')\big)
\end{eqnarray*}
and this is equal to
\begin{eqnarray*}
    \lefteqn{
       \epsilon(c) \sum_{(a),(c'),(a')} 
       \Big(\big(a^{(1)} \vdash ((c')^{(1)}\dashv (a')^{(1)})\big)\dashv
      \big( S\big(a^{(2)} \vdash ((c')^{(2)}\dashv (a')^{(2)})\big)\big)\Big)
      }
      ~~~~~~~~~~~~\\
 & &     ~~~~~~~~~~~\otimes  \big(\mathbf{1}\dashv 
           \big(a^{(3)} \vdash ((c')^{(3)}\dashv (a')^{(3)})\big)\big)\\
 & = &
     \epsilon(c) \sum_{(a),(c'),(a')} 
       \Big(a^{(1)} \vdash \big((c')^{(1)}\dashv ((a')^{(1)}\dashv
       S((a')^{(2)})) \dashv S((c')^{(2)})\dashv S(a^{(2)})\big)\Big)\\       
 & &      ~~~~~~~~~~~~ ~~~~~~~\otimes \big(\mathbf{1}\dashv 
           \big((a^{(3)} \dashv (c')^{(3)})\dashv (a')^{(3)}\big)\big)\\
 & = & \epsilon(c) \sum_{(a),(c')} 
      \Big(a^{(1)} \vdash \big(((c')^{(1)}\dashv S((c')^{(2)}))\dashv
                                           S(a^{(2)})\big)\Big)\\
  & &      ~~~~~~~~~~~~ ~~~~~~~\otimes \big(\mathbf{1}\dashv 
           \big(a^{(3)} \dashv a'\big)\big)\\   
   & = & \epsilon(c) \sum_{(a)}  
          \Big(a^{(1)} \vdash \big(c'\dashv S(a^{(2)})\big)\Big)
          \otimes \big(\mathbf{1}\dashv 
           \big(a^{(3)} \dashv a'\big)\big) \\
 & = & \epsilon(c) \sum_{(h)}
            \big(h^{(1)}.c'\big)\otimes\big(h^{(2)}h'\big)
 \end{eqnarray*}
proving eqn (\ref{EqCompVdashTransferred}).\\
{\it 2.} Clear for any bar-unital balanced dialgebra.\\
{\it 3.} Straight-forward computation using 
$\mathsf{Prim}(\mathcal{H})=\mathsf{Prim}(E_\mathcal{H})\oplus
\mathsf{Prim}(H_\mathcal{H})$ where the latter is well-known to be a Lie
algebra and the former is abelian.\\
{\it 4.} For each $b\in b$ and $h\in H$, we get 
\[ 
 (\mathbf{1}_B\otimes\mathbf{1}_H)\dashv (b\otimes h)
 =(\mathbf{1}_B\otimes \big(\Phi_B(b)h\big),
\]
proving the first statement. Moreover
\begin{eqnarray*}
   \sum_{(b),(h)} (b^{(1)}\otimes h^{(1)})\dashv 
   \big(S(b^{(2)}\otimes h^{(2)})\big))
  & = &\sum_{(b),(h)} b^{(1)}\otimes  h^{(1)}S_H(h^{(2)})
      S_H\big(\Phi_B(b^{(2)})\big) \\
  & = & \sum_{(b)}b^{(1)}\otimes S_H\big(\Phi_B(b^{(2)})\big),
\end{eqnarray*}
proving the form of the generalized idempotents, and since the 
$K$-linear map $B\to B\otimes H$ given by
$\big(\mathrm{id}_B\otimes (S_H\circ \Phi_B)\big)\circ \Delta$ is an 
injective morphism of $C^3$-coalgebra, the statement is proved.
\end{prooof}

\noindent The third statement had been proved by Simon Covez in his Master thesis
\cite{Cov1}
in the differential geometric context of digroups,  
compare with Section \ref{SecCoalgebraStructuresForPointedManifoldsETC}.

\begin{exem}
As an example, let us compute the Suschkewitsch decomposition for the augmented rack bialgebra
$K[X]$ where $p:X\to G$ is an augmented pointed rack, see Example \ref{augmented_rack_algebra}. 
By the above theorem, part {\it 4.}, its associated 
cocommutative augmented Hopf dialgebra decomposes as $B\otimes H$, where the Hopf algebra $H=K[G]$ 
is the standard group algebra and $B=K[X]$. The generalized idempotents are in this case
$$E_{B\otimes H}\,=\,\{b\otimes p(b)^{-1}\,|\,b\in B\}.$$
\end{exem}   

We finish this section with a formula relating \emph{universal 
algebras}: The functor associating to any dialgebra $A$
its Leibniz algebra $A^-$ via eqn (\ref{EqDefLeibnizBracketForDialgebras})
is well-known to have a left adjoint (see e.g.~\cite{Lod2001}) associating
to any Leibniz algebra $(\mathfrak{h},[~,~])$ 
its (in general non bar-unital) \emph{universal enveloping dialgebra 
$\mathsf{Ud}(\mathfrak{h})$ associated to $\mathfrak{h}$} defined by
\begin{equation}\label{EqDefUniversalDialgebra}
  \mathsf{Ud}(\mathfrak{h})=
    \mathfrak{h}\otimes \mathsf{U}(\overline{\mathfrak{h}}).
\end{equation}
But also in the category of bar-unital balanced dialgebras, there is such
a left adjoint: To any Leibniz algebra $(\mathfrak{h},[~,~])$, we associate 
its \emph{universal balanced
bar-unital enveloping dialgebra $\widetilde{\mathsf{Ud}}(\mathfrak{h})$}
\begin{equation}\label{EqCompUniversalBarUnitalBalancedDialgebra}
   \widetilde{\mathsf{Ud}}(\mathfrak{h})=
     \mathsf{UAR}(\mathfrak{h})\otimes \mathsf{U}(\overline{\mathfrak{h}}).
\end{equation}
Before proving the theorem, we note that
$\widetilde{\mathsf{Ud}}(\mathfrak{h})=
\mathsf{Ud}(\mathfrak{h})\oplus \mathsf{U}(\overline{\mathfrak{h}})$
is obtained by adjoining
a balanced bar-unit to $\mathsf{Ud}(\mathfrak{h})$.
\begin{theorem}
 For any Leibniz algebra $(\mathfrak{h},[~,~])$, the assignment
 $\mathfrak{h}\to \widetilde{\mathsf{Ud}}(\mathfrak{h})$ defines a 
 left adjoint
 functor to the functor associating to any bar-unital balanced
 dialgebra its commutator Leibniz algebra.
\end{theorem}
\begin{prooof}
  Clearly,  $\widetilde{\mathsf{Ud}}(\mathfrak{h})$ is the cocommutative
  Hopf dialgebra associated to the universal augmented rack bialgebra 
  $(\mathsf{UAR}(\mathfrak{h}),\Phi_\mathfrak{h},
  \mathsf{U}(\overline{\mathfrak{h}}),\ell)$ (cf.~Theorem 
  \ref{TAugmentedRackBialgebraOtimesHopfIsHopfDialgebra}) 
  which in turn is associated to the Leibniz algebra $(\mathfrak{h},[~,~])$
  (cf.~Theorem \ref{TUniversalLeibnizToAugmentedRackBialgebra}).
  Since both assignments are functorial, it follows that the 
  assignment $\mathfrak{h}\to\widetilde{\mathsf{Ud}}(\mathfrak{h})$
  is a functor. It remains to prove the universal property: Let
  $(\mathfrak{h},[~,~])$ be a Leibniz algebra, let 
  $(A,\mathbf{1},\vdash,\dashv)$ a bar-unital balanced dialgebra, and let
  $\varphi:\mathfrak{h}\to A^-$ be a morphism of Leibniz algebras. It follows
  that the $K$-linear map 
  $\Phi_A\circ \varphi:\mathfrak{h}\to A_\mathrm{ass}^1$ vanishes on
  the two-sided ideal $Q(\mathfrak{h})$ and descends to a morphism 
  $\overline{\varphi}$ of
  the quotient Lie algebra $\overline{\mathfrak{h}}$ to 
  $A_\mathrm{ass}^1$ with its 
  commutator Lie bracket such that
  $\overline{\varphi}\circ p= \Phi_A\circ \varphi$. Hence there is a unique 
  morphism
  $\mathsf{U}(\overline{\varphi}):\mathsf{U}(\overline{\mathfrak{h}})\to
  A_\mathrm{ass}^1$ of unital associative algebras extending $\overline{\varphi}$.
  Define the $K$-linear map 
  $\hat{\varphi}:\widetilde{\mathsf{Ud}}(\mathfrak{h})=
  \mathsf{U}(\overline{\mathfrak{h}})\oplus \mathsf{Ud}(\mathfrak{h})\to A$ by
  (for all $u,v\in\mathsf{U}(\overline{\mathfrak{h}})$ and $x\in \mathfrak{h}$):
  \[
    \hat{\varphi}(u)=i_A\big(\mathsf{U}(\overline{\varphi})(u)\big)
    ~~\mathrm{and}~~
    \hat{\varphi}(x\otimes v) = 
    \varphi(x)\dashv i_A\big(\mathsf{U}(\overline{\varphi})(v)\big)
  \]
  where we recall the natural injection of unital algebras
  $i_A: A_\mathrm{ass}^1\to A$ given by $i_A\big(\Phi_A(a)\big)=\mathbf{1}\dashv a$
  for all $a\in A$.
  We shall show that $\hat{\varphi}$ is a morphism of augmented dialgebras:
  We compute for all $u,u',u''\in \mathsf{U}(\overline{\mathfrak{h}})$, using that
  $i_A$ and $\mathsf{U}(\overline{\varphi})$ are morphisms of unital associative
  algebras and that in the image of $i_A$, we can use the multiplication symbols
  $\vdash$ and $\dashv$ arbitrarily:
  \[
     \hat{\varphi}(u'uu'')=
     i_A\big(\mathsf{U}(\overline{\varphi})(u')\big)
     \vdash i_A\big(\mathsf{U}(\overline{\varphi})(u)\big)\dashv
             i_A\big(\mathsf{U}(\overline{\varphi})(u'')\big)
            = \hat{\varphi}(u')\vdash  \hat{\varphi}(u) \dashv\hat{\varphi}(u'')
  \]
  showing the fact that $\hat{\varphi}$ preserves the bimodule structures on
  the first component of $\widetilde{\mathsf{Ud}}(\mathfrak{h})$. Next we have
  for all $x_1,x\in\mathfrak{h}$
  \begin{eqnarray*}
  \varphi\big(p(x_1).x\big) & = & \varphi([x_1,x])
       =\varphi(x_1)\vdash \varphi(x)-\varphi(x)\dashv \varphi(x_1) \\
       & = & i_A \big(\overline{\varphi}\big(p(x_1)\big)\big)\vdash \varphi(x)
          -\varphi(x)\dashv i_A \big(\overline{\varphi}\big(p(x_1)\big)\big)         
 \end{eqnarray*}
 and by induction on $k$ in $u'=p(x_1)\cdots p(x_k)\in
 \mathsf{U}(\overline{\mathfrak{h}})$ and $x_1,\ldots,x_k\in\mathfrak{h}$,
 we prove
 \[
 \varphi(u'.x) =  \sum_{(u')}i_A 
     \big(\overline{\varphi}(u'^{(1)})\big)\vdash \varphi(x)
      \dashv  \big(\overline{\varphi}\big(S(u'^{(2)})\big)\big).
 \]
Now, for all $u',u'',v\in \mathsf{U}(\overline{\mathfrak{h}})$ and $x\in 
  \mathfrak{h}$, we get:
  \begin{eqnarray*}
   \hat{\varphi}(u'.(x\otimes v).u'') & = & 
       \sum_{(u')}\hat{\varphi}\big(((u')^{(1)}.x)\otimes ((u')^{(2)}vu'')\big)\\
       &= &\sum_{(u')}\varphi((u')^{(1)}.x)\dashv 
       i_A\big(\mathsf{U}(\overline{\varphi})((u')^{(2)}vu'')\big) \\
     & = &\sum_{(u')} i_A \big(\mathsf{U}(\overline{\varphi})((u')^{(1)})\big)\vdash 
         \varphi(x)
       \dashv i_A\big(\mathsf{U}(\overline{\varphi})\big(S((u')^{(2)})(u')^{(3)}vu''
          \big)
                \big)\\
     & = &   i_A\big(\mathsf{U}(\overline{\varphi})(u')\big)\vdash 
         \varphi(x) \dashv i_A\big(\mathsf{U}(\overline{\varphi})(v)\big)
            \dashv i_A\big(\mathsf{U}(\overline{\varphi})( u'')\big) \\
    & = & i_A\big(\mathsf{U}(\overline{\varphi})(u')\big)\vdash 
            \hat{\varphi}(x\otimes v)\dashv
            i_A\big(\mathsf{U}(\overline{\varphi})( u'')\big)
  \end{eqnarray*}
 showing the fact that $\hat{\varphi}$ preserves the bimodule structures on
  the second component of $\widetilde{\mathsf{Ud}}(\mathfrak{h})$. Hence
  $\hat{\varphi}$ is a morphism of bar-unital (augmented) dialgebras.
  The uniqueness of $\hat{\varphi}$ follows from the universal property
  of $\mathsf{U}(\overline{\mathfrak{h}})$
\end{prooof}

\section{Coalgebra Structures for pointed manifolds with 
multiplication}
  \label{SecCoalgebraStructuresForPointedManifoldsETC}

In this section, the symbol $\mathbb{K}$ denotes either the field of all real
numbers, $\mathbb{R}$, or the field of all complex numbers, $\mathbb{C}$.
We define here the monoidal category of pointed manifolds, and exhibit the 
\emph{Serre functor}
sending a pointed manifold to the coalgebra of point-distributions supported 
in the distiguished point. We recall further that this is a strong monoidal 
functor.
Further down, we will study Lie (semi) groups, Lie racks, and Lie digroups
as examples of this 
construction, motivating geometrically the notions of a rack bialgebra and 
of a Hopf dialgebra. 

\subsection{Pointed manifolds with multiplication(s)}

Recall first the \emph{category of all pointed manifolds} $\mathcal{M}f*$ 
whose objects
consist of pairs $(M,e)$ where $M$ is a non-empty differentiable manifold and 
$e$ is an element
of $M$ and whose morphisms $(M,e)\to (M',e')$ are given by all smooth maps 
$\phi:M\to M'$ of the underlying manifolds such that $\phi(e)=e'$.
Recall that the cartesian product $\times$ makes  $\mathcal{M}f*$ into a 
\emph{monoidal
category} (see e.g. \cite[p.161-170; 251-257]{Mac98} for definitions) by setting 
$(M,e_1)\times (N,e_2):=\big(M\times N,(e_1,e_2)\big)$ with 
the one-point set $\big(\{\mathrm{pt}\},\mathrm{pt}\big)$ as unit object and 
the usual associators, left-unit and right-unit identifications borrowed from 
the category of sets.
This monoidal category is \emph{symmetric} by means of the usual (tensor) flip map
$\tau_T=\tau_{T~(M,N)}:M\times N\to N\times M: (x,y)\mapsto (y,x)$ where the 
pair of distinguished points is also interchanged.\\
By simply forgetting about the differentiable structure we get the category
of \emph{pointed sets}.

Recall that a \emph{pointed manifold with multiplication} is a triple
$(M,e,\mathbf{m})$ where $(M,e)$ is a pointed manifold, and 
$\mathbf{m}:(M,e)\times(M,e)\to (M,e)$
is a smooth map of pointed manifolds, i.e. is a smooth map $M\times M\to M$ such that
$\mathbf{m}(e,e)=e$. Moreover, a pointed manifold with multiplication
will be called \emph{left-regular} (resp.~\emph{right-regular}) if all the 
left (resp.~right) multiplication maps $y\mapsto \mathbf{m}(x,y)$
(resp.~$y\mapsto \mathbf{m}(y,x)$ are diffeomorphisms.
Morphisms of pointed manifolds with multiplication
$(M,e,\mathbf{m})\to (M',e',\mathbf{m}')$ are smooth maps of pointed manifolds
$\phi:(M,e)\to (M',e')$ such that
\begin{equation}
    \phi\circ \mathbf{m}= \mathbf{m}'\circ (\phi\times \phi).
\end{equation}
The obvious generalization are a finite
number of maps $M^{\times n}\to M$ with $n\geq 1$) arguments.\\
Again by forgetting about differentiable structures, we get the category of
\emph{pointed sets with multiplication}.

\subsection{Coalgebra Structure for distributions supported in one point}
   \label{SubSecCoalgStrucForDistributions}



For any pointed manifold $(M,e)$ recall the $\mathbb{K}$-vector space
\begin{eqnarray}
    \lefteqn{\mathcal{E}'_e(M):=} \nonumber \\
    & & \big\{T:\mathcal{C}^\infty(M,\mathbb{K})\to \mathbb{K}~|~T
        \mathrm{~is~a~continuous~linear~map~and~}\mathrm{supp}(T)=\{e\}\big\}
        \nonumber \\
        & & \label{DefEPrimeXOfM}
\end{eqnarray}
of all distributions supported in the singleton $\{e\}$, 
see e.g. \cite[Ch.6,7]{Rud73} for definitions). 
Now let $\phi:M\to M'$ be a smooth map such that
$\phi(e)=e'$. For any distribution $S\in \mathcal{E}'_e(M)$ and any smooth function
$f'\in \mathcal{C}^\infty(M',\mathbb{K})$ the well-known prescription
\begin{equation}
    (\phi_*S)(f'):=S(f'\circ \phi)
\end{equation}
gives a well-defined distribution $\phi_*S$ on the target manifold $M'$
supported in $e'=\phi(e)$, and the
map $\phi_*:\mathcal{E}'_e(M)\to \mathcal{E}'_{e'}(M')$ is a $\mathbb{K}$-linear
map (which is continuous). Clearly for three pointed manifolds $(M,e)$, $(M',e')$, and
$(M'',e'')$ with smooth maps $\phi:(M,e)\to (M',e')$ and $\psi:(M',e')\to (M'',e'')$
we get
\begin{equation}
    (\psi\circ \phi)_*=\psi_*\circ \phi_*~~~\mathrm{and}~~~
        \big(\mathrm{id}_M\big)_* =\mathrm{id}_{\mathcal{E}'_e(M)}.
\end{equation}
This defines a covariant functor $F:\mathcal{M}f*\to \mathbb{K}\mathbf{Vect}$ to
the category of all $\mathbb{K}$-vector spaces by associating to any pointed
manifold $(M,e)$ the $\mathbb{K}$-vector space $F(M,e):=\mathcal{E}'_e(M)$, and
to any smooth map $(M,e)\to (M,e')$ the linear map 
$F(\phi):=\phi_*:\mathcal{E}'_e(M)\to\mathcal{E}'_{e'}(M')$.
We call this functor {\it Serre functor} in tribute to the predominant role it plays in
\cite{Ser64}. It is one of the main objects of this article. 

There is, however, much more structure in this functor: First any distribution
space $\mathcal{E}'_e(M)$ contains a canonical linear form 
$\epsilon=\epsilon_e:\mathcal{E}'_e(M)\to\mathbb{K}$ defined by
\begin{equation}\label{EqDefCoUnitDistributions}
    \epsilon_e(T):=T(1),
\end{equation}
where $1$ denotes the constant function $M\to\mathbb{K}$ whose only value is equal to 
$1\in\mathbb{K}$. Moreover each space $\mathcal{E}'_e(M)$ contains a canonical
element $\mathbf{1}=\mathbf{1}_e$ defined by the well-known \emph{delta distribution} 
\begin{equation}\label{EqDefUnitDistributions}
   \mathbf{1}_e = \delta_e: f\mapsto f(e),  
\end{equation}
and we clearly have
\begin{equation}
    \epsilon_e(\mathbf{1}_e)=1.
\end{equation}
Moreover, both $\epsilon_e$ and $\mathbf{1}_e$ are natural in the following sense:
let $\phi:(M,e)\to (M',e')$ be a smooth map of pointed differentiable manifolds.
Then it is straight-forward to check that
\begin{equation}
    \epsilon_{e'}\circ \phi_*= \epsilon_e~~~\mathrm{and}~~~
       \phi_*(\mathbf{1}_e)=\mathbf{1}_{e'}.
\end{equation}

Recall the well-known \emph{tensor product} or \emph{direct product} of 
two distributions (cf. \cite[p.~403]{Tre67}):
More generally, let $M$ and $N$ be two differentiable 
manifolds, and let
$S\in\mathcal{D}'(M)$ and $T\in\mathcal{D}'(N)$ be two distributions
(where the symbol $\mathcal{D}'(M)$ denotes the continuous dual space of
the test function space $\mathcal{D}(M)$ of all smooth $\mathbb{K}$-valued 
functions with compact support). Let $f:M\times N\to\mathbb{K}$ be a smooth function
with compact support $K\subset M\times N$, and let 
$K_1:=\mathrm{pr}_M(K)\subset M$, $K_2=\mathrm{pr}_{N}(K)\subset N$, whence
$K$ is a subset of the compact set $K_1\times K_2$. Let 
$T^{(2)}:\mathcal{D}(M\times N)\to \mathrm{Fun}(M,\mathbb{K})$ be the following
map: For each $x\in M$, let $f_x\in \mathcal{D}(N)$ be the partial function
$y\mapsto f(x,y)$. Then we set 
\[
   \big(T^{(2)}(f)\big)(x)= T(f_x).
\]
The superscript $(2)$ means here that we see $f$ as a function of its second variable only,
when applying the distribution $T$.
 
Upon using the approximation theorem of any distribution by a sequence of 
regular distributions (see e.g. \cite[p.157, Thm.6.32]{Rud73}), one can 
show that $T^{(2)}(f):M\to\mathbb{K}$ is a smooth function
having compact support in $K_1$. It is clear that 
$T^{(2)}:\mathcal{D}(M\times N)\to \mathcal{D}(M)$ is linear, and it can 
be shown by the same approximation 
theorem that $T^{(2)}$ is continuous, see \cite[p.416]{Tre67}. 
It follows that the map
\[
     (S,T)\mapsto \Big(f\mapsto S\big(T^{(2)}(f)\big)\Big)
\]
is a well-defined $\mathbb{K}$-bilinear map 
$\mathcal{D}'(M)\times \mathcal{D}'(N)\to
\mathcal{D}'(M\times N)$, and there is thus a unique linear map
(where $\otimes$ denotes the usual algebraic tensor product over $\mathbb{K}$)
\begin{equation}
    F_{2~M,N}: \mathcal{D}'(M)\otimes \mathcal{D}'(N)\to
\mathcal{D}'(M\times N)
\end{equation}
such that for all $f\in \mathcal{D}(M\times N)$ we have
\[
    \big(F_{2~M,N}(S\otimes T)\big)(f) = S\big(T^{(2)}(f)\big).
\]
Note also that it can be shown that the right hand side is equal to  
$T\big(S^{(1)}(f)\big)$ where the 
notation is self-explanatory.
Furthermore, it is not hard to see that for two distributions supported in one 
point,
i.e. $S\in \mathcal{E}'_{e_1}(M)$ and $T\in \mathcal{E}'_{e_2}(N)$ the distribution
$F_{2~M,N}(S\otimes T)$ is supported in $(e_1,e_2)$, i.e. is an element of
$\mathcal{E}'_{(e_1,e_2)}(M\times N)$. We shall denote the restriction of 
the map $F_{2~M,N}$ to $\mathcal{E}'_{e_1}(M)\otimes\mathcal{E}'_{e_2}(N)$ by 
the same symbol $F_{2~M,N}$.
For three pointed manifolds $(M,e_1)$, $(N,e_2)$, and $(P,e_3)$, let
$\alpha_{M,N,P}:M\times (N\times P)\to (M\times N)\times P$ be the usual 
associator for the monoidal category of all sets 
(see \cite{Mac98}), and for three vector spaces $V,W,X$ over 
$\mathbb{K}$, let
$\beta_{V,W,X}:V\otimes (W\otimes X)\to (V\otimes W)\otimes X$ be the 
well-known associator
for the monoidal category of all vector spaces. By using the definitions, it is 
not hard
to see that the following identity holds
\begin{eqnarray}  \nonumber
   \lefteqn{ F_{2~(M\times N), P}\circ 
    \big(F_{2~M,N}\otimes \mathrm{id}_{\mathcal{E}'_{e_3}(P)}\big)
      \circ \beta_{\mathcal{E}'_{e_1}(M),\mathcal{E}'_{e_2}(N),\mathcal{E}'_{e_3}(P)}  }  \\
     &  = & 
      \big(\alpha_{M,N,P}\big)_*\circ 
        F_{2~M,(N\times P)} \circ 
        \big(\mathrm{id}_{\mathcal{E}'_{e_1}(M)}\otimes F_{2~N,P}\big)   
        \label{associativity_of_F}
\end{eqnarray}
hence eqn (3) of \cite[p.255]{Mac98}. In the same vein, the two diagrams 
in eqn (4) of \cite[p.256]{Mac98} are satisfied upon setting 
$F_0=\epsilon_{\mathrm{pt}}$
and $\lambda:\mathbb{K}\otimes V \to V$ and $\rho:V\otimes \mathbb{K}\to V$
the usual left-unit and right-unit identifications in the monoidal category 
of vector spaces.\\
Let $(M',e_1')$ and $(N',e_2')$ two other 
pointed differentiable 
manifolds, and let $\phi:(M,e_1)\to (M',e_1')$ and $\psi:(N,e_2)\to (N',e_2')$ 
two smooth maps
of pointed differentiable manifolds. It is a straight-forward check that
the map $F_{2~M,N}$ to 
is natural in the following sense
\begin{equation}   \label{naturality_of_F}
    F_{2~M',N'}\circ \big(\phi_*\otimes \psi_*\big)  
         = (\phi\times\psi)_*\circ F_{2~M,N}.
\end{equation}
Moreover, note that the map $F_0=\epsilon_{\mathrm{pt}}$ 
(see eqn (\ref{EqDefCoUnitDistributions})) defines an isomorphism
of $\mathcal{E}'_{\mathrm{pt}}(\{\mathrm{pt}\})$ to
 $\mathbb{K}$ which had already been seen to be natural.\\
 As a result, the functor $F$ is a \emph{monoidal functor} in the sense of 
 \cite[p.255-257]{Mac98}. Moreover, since the category $\mathcal{M}f*$ is even
 a \emph{symmetric monoidal category} by means of the canonical flip map 
 $\tau_{M,N}:M\times N\to N\times M: (x,y)\to (y,x)$, see e.g.
 \cite[p.252-253]{Mac98}, and the monoidal category $\mathbb{K}$-\textbf{vect}
 is also symmetric, 
 it is not hard to see that the monoidal functor is also symmetric,
 see e.g. \cite[p.257]{Mac98} for definitions.

 We shall now show that the monoidal functor $F$ is \emph{strong}, i.e. that
 $F_0=\epsilon_{\mathrm{pt}}$ and $F_{2~M,N}$ are isomorphisms. This is clear for 
  $\epsilon_{\mathrm{pt}}$. Recall that for each distribution $T$
  in $\mathcal{E}'_e(V)$ (where $V$ is a nonempty open set in $\mathbb{R}^m$
  containing the point $e$), 
  there is nonnegative integer $l$ (called the {\it order} of the distribution) such that
    \[
         T= \sum_{r=0}^l\sum_{\mathbf{k}\in \mathbb{N}^m, |\mathbf{k}|=r}
                     c_{\mathbf{k}} 
            \frac{\partial^{\mathbf{k}}}{\partial x_{\mathbf{k}}}\big(\delta_e)
  \]
  where $c_{\mathbf{k}}\in\mathbb{K}$ for each multi-index $\mathbf{k}$, see 
  e.g. \cite[p.150, Thm. 6.25]{Rud73}. In a slightly more algebraic manner
  we can express this as follows: let $E$ be a finite-dimensional real vector 
  space, let $V\subset E$ be an open set containing $e\in E$. Then we have the 
  following linear isomorphism $\Phi_\mathsf{S}:\mathsf{S}(E)\to \mathcal{E}'_e(V)$
  given by $\Phi_\mathsf{S}(\mathbf{1})=\delta_e$ and for any positive integer $k$
  and vectors $w^{(1)},\ldots,w_{(k)}\in E$ and 
  $f\in\mathcal{C}^\infty(V,\mathbb{K})$
  \begin{eqnarray}\label{EqCompPhiOne}
   \lefteqn{\Phi_\mathsf{S}(w^{(1)}\bullet\ldots\bullet w_{(k)})(f)} \nonumber \\
   & = & \left.\frac{\partial^k \Big(f\big(e+s_1 w^{(1)}+\cdots +s_kw_{(k)}\big)\Big)}
          {\partial s_1\cdots\partial s_k}\right|_{s_1=0,\ldots,s_k=0}
  \end{eqnarray}
  where $\bullet$ denotes the commutative multiplication in the symmetric algebra,
  see Appendix A.
  Using the fact that the inclusion map 
  $\iota_{U_\alpha}:U_\alpha\to M$ of any chart domain of $M$ such that $e\in U_\alpha$ defines
  an isomorphism 
  $\iota_{U_\alpha *}:\mathcal{E}'_e(U_\alpha)\to \mathcal{E}'_e(M)$, 
  and that any chart
  $\varphi_\alpha:U_\alpha \to V_\alpha\subset \mathbb{R}^m$ defines an isomorphism
  $\varphi_{\alpha *}:\mathcal{E}'_e(U_\alpha)\to 
  \mathcal{E}'_{\varphi_\alpha(e)}(V_\alpha)$, we can conclude that there is a 
  linear isomorphism
  \begin{equation}\label{EqDefinitionOfPhiAlpha}
      \Phi_\alpha=:
      \iota_{U_\alpha *}\circ\varphi_{\alpha *}^{-1}\circ\iota_{V_\alpha *}^{-1}
         \circ \Phi_\mathsf{S}:
         \mathsf{S}(\mathbb{R}^m)\to \mathcal{E}'_e(M)
  \end{equation}
with the symmetric coalgebra $\mathsf{S}(\mathbb{R}^m)$ on $\R^m$ 
(see Appendix A) computed as follows
  \begin{eqnarray}
      w^{(1)}\bullet\cdots\bullet w_{(k)}  & \mapsto &
      \Big(f\mapsto \sum_{i_1,\ldots,i_k=1}^m
        \frac{\partial^k (f|_{U_\alpha}\circ \varphi_\alpha^{-1})}
          {\partial x_{i_1}\cdots\partial x_{i_k}}\big(\varphi_\alpha(e)\big)
            w_{(1)i_1}\cdots w_{(k)i_k}\Big).\nonumber \\
            & & \label{EqComputationOfPhiAlpha}
  \end{eqnarray}
  where we write $\mathbb{R}^m\ni w_{(j)}=\sum_{i=1}^m w_{(j)i}e_i$ 
  where all the $w_{(j)i}$
  are real numbers and $e_1,\ldots,e_m$ is the canonical base of $\mathbb{R}^m$.
  Note that for the particular case of $M$ being an open set $V$ of 
  $\mathbb{R}^m$ and the chart $\varphi_\alpha$ being the identity map 
  the map $\Phi_\alpha$ (see eqs (\ref{EqDefinitionOfPhiAlpha}) 
  and (\ref{EqComputationOfPhiAlpha})) coincides with the 
  canonical map $\Phi_\mathsf{S}$, see 
  eqn (\ref{EqCompPhiOne}).
  
  For two pointed manifolds $(M,e_1)$ and $(N,e_2)$ and given charts
  $(U_\alpha,\varphi_\alpha)$ of $M$ such that $e_1\in U_\alpha$ and
  $(\tilde{U}_\beta,\tilde{\varphi}_\beta)$ of $N$ such that
  $e_2\in \tilde{U}_\beta$, we thus have linear isomorphisms
  $\Phi_\alpha:\mathsf{S}(\mathbb{R}^m)\to \mathcal{E}'_{e_1}(M)$,
  $\tilde{\Phi}_\beta:\mathsf{S}(\mathbb{R}^n)\to \mathcal{E}'_{e_2}(N)$,
  and 
  $\Phi_{\alpha,\beta}:\mathsf{S}(\mathbb{R}^{m+n})\to 
  \mathcal{E}'_{(e_1,e_2)}(M\times N)$ (upon using the product chart
  $(U_\alpha\times \tilde{U}_\beta,
  \varphi_\alpha\times \tilde{\varphi}_\beta)$).
  Using the above definitions, one can compute
  that
  \[
      F_{2~M,N}\circ \Big(\Phi_\alpha\otimes \tilde{\Phi}_\beta\Big)
         = \Phi_{\alpha,\beta}\circ \Theta_{m,n}
  \]
  where $\Theta_{m,n}:\mathsf{S}(\mathbb{R}^m)\otimes \mathsf{S}(\mathbb{R}^n)
     \to \mathsf{S}(\mathbb{R}^{m+n})$ denotes the natural isomorphism
     of commutative associative unital algebras
     induced by the obvious inclusions 
     $\mathbb{R}^m\hookrightarrow \mathbb{R}^{m+n}$
     (first $m$ coordinates) and $\mathbb{R}^n\hookrightarrow \mathbb{R}^{m+n}$
     (last $n$ coordinates). It follows that 
     the natural map  $F_{2~M,N}$ is equal to
     $\Phi_{\alpha,\beta}\circ \Theta_{m,n}
     \circ\big(\Phi_\alpha^{-1}\otimes \tilde{\Phi}_\beta^{-1}\big)$ and is thus
     a linear isomorphism, whence the functor $F$ is a strong monoidal functor.
     
     In order to define more structure, let us consider the well-known 
     \emph{diagonal map} $\mathrm{diag}_M:M\to M\times M$ defined by
     \[
        \mathrm{diag}_M(x)= (x,x)   
     \]
     for all $x\in M$. Clearly, $\mathrm{diag}_M$ is a smooth map of pointed
     manifolds $(M,e)\to \big(M\times M, (e,e)\big)$. Moreover, the diagonal map 
     is clearly natural in the sense that
     \[
           \mathrm{diag}_{M'}\circ \phi 
             = \big(\phi\times \phi\big)\circ\mathrm{diag}_{M} 
     \]
     for any smooth map $\phi:(M,e)\to (M',e')$ of pointed differentiable manifolds.
     In other words, the class of diagonal maps $\mathrm{diag}_M$ constitutes
     a natural transformation from the identity functor to the \emph{diagonal
     functor} $(M,e)\to \big(M\times M, (e,e)\big)$ and $\phi\mapsto \phi\times \phi$.
     Define the following linear map $\Delta=\Delta_e=\Delta_{(M,e)}:
     \mathcal{E}'_e(M)\to 
           \mathcal{E}'_e(M)\otimes\mathcal{E}'_e(M)$ by
     \begin{equation}\label{EqDefDeltaX}
        \Delta_e:= F_{2~M,M}^{-1}\circ\mathrm{diag}_{M~*}.
     \end{equation}
This definition has avatars with more than two tensor factors. Indeed, observe that 
the naturality relation 
(\ref{naturality_of_F}) implies for $\phi={\rm id}_M$ and $\psi={\rm diag}_M$ that 
\begin{eqnarray*}
\lefteqn{({\rm id}_{M}\times\mathrm{diag}_{M})_*} \\
 & = & F_{2~M,(M\times M)}\circ({\rm id}_{\mathcal{E}'_e(M)}\times F_{2~M,M})
 \circ({\rm id}_{\mathcal{E}'_e(M)}
\times\Delta_e)\circ F_{2~M,M}^{-1}.
\end{eqnarray*}
Similarly, we have relations of this type for any number of tensor factors.

In the following, we invite the reader to look again at Appendix A for definitions
and notations about coalgebras.\\
We have the following
     \begin{theorem}
      With the above notations:
      \begin{enumerate}
       \item The $\mathbb{K}$-vector space $\mathcal{E}'_e(M)$ equipped with the 
       linear maps $\Delta_e$ (cf. eqn (\ref{EqDefDeltaX})), $\epsilon_e$
       (cf. eqn (\ref{EqDefCoUnitDistributions}), and $\mathbf{1}_e$
       (cf. eqn (\ref{EqDefUnitDistributions}) is a $C^5$-coalgebra
       which is (non canonically) isomorphic
       to the standard symmetric coalgebra 
       $$\big(\mathsf{S}(\mathbb{R}^m),\epsilon,\Delta,\mathbf{1}\big).$$
       \item The above strong symmetric monoidal functor $F$ extends to a 
       functor --also denoted
       by $F$-- from
       $\mathcal{M}f*$ to the symmetric monoidal category of $C^5$-coalgebras 
       over 
       $\mathbb{K}$.
       \item The subspace of all primitive elements of the coaugmented coalgebra
          $\mathcal{E}'_e(M)$ is natural isomorphic to the tangent space
          $T_eM$. Moreover for each smooth map $\phi:(M,e)\to (M',e')$ of pointed
          manifolds the coalgebra morphism 
          $\phi_*:\mathcal{E}'_e(M)\to\mathcal{E}'_{e'}(M')$ induces the tangent
          map $T_e\phi:T_eM\to T_{e'}M'$.
      \end{enumerate}

     \end{theorem}
\begin{prooof}
\begin{enumerate}
\item[(a)]{\bf coassociativity of $\Delta_e$:} This follows from the coassociativity-diagram of 
$\mathrm{diag}_{M~*}$ by first taking the induced diagram between distribution spaces
which reads then
$$(\alpha_{M,M,M})_*\circ({\rm id}_{M}\times\mathrm{diag}_{M})_*\circ\mathrm{diag}_{M~*}\,=\,
(\mathrm{diag}_{M}\times{\rm id}_{M})_*\circ\mathrm{diag}_{M~*}.$$
Starting for example on the left hand side, one replaces the map $\mathrm{diag}_{M~*}$ by 
$F_{2~M,M}\circ\Delta_e$, and also the map 
$({\rm id}_{M}\times\mathrm{diag}_{M})_*$ by 
$$
    F_{2~M,(M\times M)}\circ({\rm id}_{\mathcal{E}'_e(M)}
        \times F_{2~M,M})\circ({\rm id}_{\mathcal{E}'_e(M)}
           \times\Delta_e)\circ F_{2~M,M}^{-1}.
$$
Now one observes that one may apply the relation (\ref{associativity_of_F}) on the left hand side.
One obtains
\begin{eqnarray*}
\lefteqn{ F_{2~(M\times M),M}\circ(F_{2~M,M}\times{\rm id}_{\mathcal{E}'_e(M)})\circ\beta\circ(
{\rm id}_{\mathcal{E}'_e(M)}\times\Delta_e)\circ
\Delta_e} \\
& = &
F_{2~(M\times M),M}\circ(F_{2~M,M}\times{\rm id}_{\mathcal{E}'_e(M)})\circ(\Delta_e\times
{\rm id}_{\mathcal{E}'_e(M)})\circ\Delta_e.
\end{eqnarray*}
One deduces coassociativity.
\item[(b)]{\bf cocommutativity of $\Delta_e$:} This follows from the symmetry of $F$ 
(already noted before) and the cocommutativity of $\mathrm{diag}_{M~*}$. 
\item[(c)]{\bf counitality of $\Delta_e$:} This follows from the counitality of 
$\mathrm{diag}_{M~*}$, i.e.
$$
  ({\rm proj}_e\times{\rm id})\circ \mathrm{diag}_{M}\,=
     \,{\rm incl}_e^1,\,\,\,\,{\rm and}\,\,\,\,
({\rm id}\times{\rm proj}_e)\circ \mathrm{diag}_{M}\,=\,{\rm incl}_e^2,
$$
where ${\rm proj}_e:M\to\{e\}$, ${\rm incl}_e^1:M\to\{e\}\times M$ and 
${\rm incl}_e^2:M\to M\times\{e\}$ are the canonical maps. Indeed, these equations induce
the corresponding equations between distribution spaces, and translating direct products 
into tensor products (and thus $({\rm proj}_e)_*$ into $\epsilon$ and $(\mathrm{diag}_{M})_*$
into $\Delta_e$), one obtains counitality.  
\item[(d)]\textbf{connectedness of $\Delta_e$:} The coalgebra
 $\mathcal{E}'_e(M)$ is isomorphic to the symmetric algebra 
 $\mathsf{S}(\mathbb{R}^n)$, and the latter is connected.
\end{enumerate}
This shows part (1) of the statement, as the isomorphy to the standard symmetric coalgebra 
has been shown above. 

The only thing which has to be shown for the second statement is the preservation of 
the coalgebra structure on the level of morphisms, which is clear.  

For the third part, consider the linear map 
$T_eM \to \mathsf{Prim}(\mathcal{E}'_e(M))$ (see Appendix A for the definition of 
the primitives $\mathsf{Prim}(C)$ of a coalgebra $C$) defined by
\[
    v\mapsto \big(f\mapsto df_e(v)\big):
\]
Indeed the right hand side is clearly in $\mathcal{E}'_e(M)$, and the Leibniz
rule for the derivative shows that this is in $\mathsf{Prim}(\mathcal{E}'_e(M))$.
Moreover the above map is clearly injective, and since 
$\mathsf{Prim}(\mathsf{S}(\mathbb{R}^n))=\mathbb{R}^n$ and $\dim(T_eM)=n$
it follows that the above map is an isomorphism of real vector spaces.
The naturality is a simple computation.
 \end{prooof}
 
The last statement means that the composed functor $\mathsf{Prim}\circ F$
of the Serre functor $F$ and the functor associating to any coalgebra $C$ its
space of primitive elements, $\mathsf{Prim}(C)$,
is naturally isomorphic to the tangent functor $T_*$ associating to any
pointed differentiable manifold $(M,e)$ its tangent space $T_eM$.

\begin{rem}  \label{remark}
   There is neither a canonically defined (i.e. not depending
   on the choice of a chart) projection from the coalgebra to its primitives,
   so the coalgebras $\mathcal{E}'_e(M)$ are isomorphic to the cofree 
   $\mathsf{S}(\mathbb{R}^m)$, but in general not naturally, 
   nor a canonically defined
   commutative multiplication (the classical convolution of distributions of
   compact support which needs the additive vector space structure).
\end{rem}

\begin{rem}
   Note also the disjoint union $\bigcup_{x\in M}\mathcal{E}'_x(M)_{(k)}$
   carries the structure of a smooth vector bundle over $M$: Its smooth sections
   coincide with the space of all differential operators of order $k$.
\end{rem}

\begin{rem}
  In case $U\subset \mathbb{R}^m$ and $V\subset \mathbb{R}^n$ are pointed open
  sets, the coalgebra morphism $\phi_*$ of a smooth map $\phi:U\to V$
  of pointed manifolds is isomorphic to the coalgebra morphism
  $\mathsf{S}(\mathbb{R}^m)\to\mathsf{S}(\mathbb{R}^n)$ induced by
  the jet of infinite order of $\phi$ at the distinguished point $e$ of $U$,
  $j^\infty(\phi)_e$, see e.g. \cite{KMS93} for further information.
  The functorial equation $(\phi\circ \psi)_*=\phi_*\circ \psi_*$ can be 
  computed out of the chain rule for higher derivatives.
\end{rem}

\subsection{Pointed manifolds with multiplication and their 
             associated bialgebras}


We can now apply the Serre functor defined in the preceding Section
\ref{SubSecCoalgStrucForDistributions} to pointed manifolds with multiplication:

\begin{theorem}
   Let $(M,e,\mathbf{m})$ be a pointed manifold with multiplication.
   Then the $C^5$-coalgebra $\mathcal{E}'_e(M)$ carries a multiplication, i.e.
   a linear map 
   $\mu=\mathbf{m}_*\circ F_{2~M,M}:
    \mathcal{E}'_e(M)\otimes\mathcal{E}'_e(M)\to\mathcal{E}'_e(M)$
   which is a morphism of $C^5$-coalgebras.\\
   In case $\mathbf{m}$ is left-unital (resp.~right unital), the 
   nonassociative $C^3I$-algebra $\mathcal{E}'_e(M)$ is left regular
   (resp.~right regular)
\end{theorem}

\begin{prooof}
The map $\mu$ exists and is linear by functoriality. We have trivially
$$
\mathrm{diag}_{M}\circ {\bf m}
      \,=\,({\bf m}\times {\bf m})\circ({\rm id}_M\times
     \tau_{M,M}\times{\rm id}_M)\circ(\mathrm{diag}_{M}\times\mathrm{diag}_{M}),
$$
and this shows that $\mu$ is a morphism of coalgebras by translating 
$\mathrm{diag}_{M}$ into $\Delta_e$ using as before the maps of type $F_{2}$.
The regularity statements are a consequence of the connectedness of
the $C^3$-coalgebra $\mathcal{E}'_e(M)$, see Lemma 
\ref{LRegularityOfLeftOrRightMultiplications}.
\end{prooof}

In the following, we shall enumerate some important (sub)categories
of pointed differentiable manifolds with multiplications.

\subsubsection{Lie groups and universal enveloping algebras}

Let $\big(G,\mathbf{m},e,(~)^{-1}\big)$ a Lie group. The following theorem is well-known
(see \cite{Ser64}):
\begin{theorem}
 The associated coalgebra with multiplication $\mu$ of the Lie group
 $\big(G,\mathbf{m},e,(~)^{-1}\big)$
 is an associative unital bialgebra (in fact, a Hopf
 algebra) isomorphic to the \textbf{universal enveloping algebra of the Lie algebra}
 $\mathfrak{g}=T_eG$ of $G$.
\end{theorem}
We just indicate the isomorphism: For any $\xi\in\mathfrak{g}$, let $\xi^+$
denote the \emph{left invariant vector field $\xi^+(g):=T_eL_g(\xi)$ 
generated by
its value $\xi\in \mathfrak{g}=T_eG$}. Then the map 
$\Phi_\mathsf{U}:\mathsf{U}(\mathfrak{g})\to F(G)$ is given by
(for all $k\in \mathbb{N}$, $\xi_1,\ldots,\xi_k\in\mathfrak{g}$ and 
$f'\in\mathcal{C}^\infty(G,\mathbb{K})$)
\begin{equation}\label{EqDefIsoUniversalEnvOfGToFOfG}
   \Phi_\mathsf{U}(\xi_1\cdots\xi_k)\big(f'\big)=
    \Big(\big(\mathbb{L}_{\xi_1^+}\circ\cdots 
              \circ\mathbb{L}_{\xi_k^+}\big)(f')\Big)(e)
\end{equation}
where $\mathbb{L}_X$ denotes the Lie derivative in the direction of the vector field 
$X$.

\noindent
Note that the identities for the inverse map $g\mapsto g^{-1}$ can be written
as
\begin{equation}
 \mathbf{m}\circ\big((~)^{-1}\times \mathrm{id}_G\big)\circ \mathrm{diag}_G
    = (g\mapsto e)
    = \mathbf{m}\circ\big(\mathrm{id}_G\times(~)^{-1}\big)\circ \mathrm{diag}_G,
\end{equation}
and an application of the functor $F$ gives the convolution identities for the 
antipode, defined by $S=\big((~)^{-1}\big)_*$.

\subsubsection{Lie semigroups and Lie monoids}


It is easy to see but presumably less known that the result of the preceding 
subsection remains true for a \emph{Lie monoid} 
$\big(G,\mathbf{m},e\big)$:

\begin{theorem}
  \begin{enumerate}
 \item The associated coalgebra with multiplication $\mu$
 of the Lie mon\-oid $\big(G,\mathbf{m},e\big)$ 
 is an associative unital bialgebra (in fact a Hopf
 algebra) isomorphic to the \textbf{universal enveloping algebra of a 
 Lie algebra}
 $\mathfrak{g}\cong T_eG$.
 \item The associated coalgebra with multiplication $\mu$
 of the right Lie group $(G,\mathbf{m},e, (~)^{-1})$ 
 is a right Hopf algebra.
 \end{enumerate}
\end{theorem}
\noindent In order to see the first statement note that it is clear that 
the associated coalgebra $C:=F(G)$ carries an 
associative unital
multiplication $\mu=\mathbf{m}_*\circ F_{2~G,G}$. The fact that the 
coalgebra is always connected implies by the Takeuchi-Sweedler argument
(see Appendix \ref{AppCoalgebras}) that the identity map $\mathrm{id}_C$
has a convolution inverse, and is thus a Hopf algebra.
Since the coalgebra $C$ is connected and cocommutative, it follows from
the \emph{Cartier-Milnor-Moore Theorem} (see e.g. \cite{Qui69}) 
that the Hopf algebra $F(G)$
is isomorphic to the universal enveloping algebra over the Lie subalgebra
$\mathfrak{g}$ of its primitive elements which is equal
to $T_eG$.\\
The second statement is an immediate consequence of the functorial properties
of $F$.

\subsubsection{(Lie) dimonoids and digroups}

Recall that a \emph{Lie dimonoid} (see e.g. $\cite{Lod2001}$)
is a pointed differentiable manifold 
$(D,e)$ equipped with two smooth \emph{associative multiplications} 
$D\times D\to D$, written $(x,y)\mapsto x\vdash y$ and 
$(x,y)\mapsto x\dashv y$ (and preserving points, i.e.
$e\vdash e=e=e\dashv e$), such that the dialgebra conditions eqs 
(\ref{EqLeftVdashEqualLeftDashvDialgebra}), 
(\ref{EqRightVdashEqualRightDashvDialgebra}), 
(\ref{EqLeftVdashCommutesWithRightDashvDialgebra}),
 (\ref{EqDefBarUnitDialgebraLeft}), and (\ref{EqDefBarUnitDialgebraRight})
 hold for all
$x,y,z\in D$ and $e$ (replacing $\mathbf{1}$): 
Hence $(D,\vdash,e)$ is a left unital Lie semigroup and 
$(D,\dashv,e)$ is a right unital Lie semigroup,
and as for dialgebras, we shall say \emph{bar-unital dimonoid}
to stress the fact that the bar-unit $e$ is among the data for the dimonoid.\\
Let us call a Lie dimonoid $(D,e,\vdash,\dashv)$ \emph{balanced} iff
in addition
for all $x\in D$ the analogue of eqn (\ref{EqDefBalancedDialgebra}) holds, i.e.
$x\vdash e=e\dashv x$.
Any Lie monoid $(G,e,\mathbf{m})$ is a Lie dimonoid by setting 
$\vdash\,=\,\dashv\,=\,\mathbf{m}$. \\
Another class of examples is obtained by the following important
\emph{augmented dimonoid}
construction (cf Example \ref{construction_dialgebra}): Let
$G$ be a Lie group, let $(D,e_D)$ be a pointed differentiable manifold, 
let $G$ smoothly
act on the left and on the right of $M$ (written $(g,x)\mapsto gx$ and
$(x,g)\mapsto xg$) such that $(gx)g'=g(xg')$ for all $g,g'\in g$ and $x\in D$,
and let $f:(D,e_D)\to (G,e)$ be a smooth map of pointed manifolds such that
for all $g,g'\in G$ and $x\in D$
\begin{equation}\label{EqInterwiningMapForBiActionsDimonoid}
         f(gxg')=gf(x)g'.
\end{equation}
Then the pointed manifold $\big(D,e_D,\vdash,\dashv\big)$
will be a (bar-unital) dimonoid by setting
\begin{equation}\label{EqDefDimonoidOutOfBiaction}
    x\vdash y := f(x)y~~~\mathrm{and}~~~x\dashv y := xf(y).
\end{equation}

A \emph{general Lie digroup} is defined (according to Liu, 
\cite[Definition 1.1]{Liu04}) to be a (bar-unital) dimonoid
$(D,e,\vdash,\dashv)$ such that the left unital Lie semigroup $(D,e,\vdash)$
is a right group and the right unital Lie semigroup $(D,e,\dashv)$ is
a left group (see Appendix \ref{AppSemigroups} for definitions): Here the 
right inverse of $x$ with respect to $\vdash$ does in general not coincide 
with the left inverse of $x$ with respect to $\dashv$. For an example, take any
Lie group $G$, set $(D,e_D)=\big((G\times G, (e,e)\big)$, define the two 
canonical $G$-actions
$g(g_1,g_2):=(gg_1,g_2)$ and $(g_1,g_2)g:=(g_1,g_2g)$
(for all $g,g_1,g_2\in g$), and let $f:G\times G\to G$
be the group multiplication. Then $(D,e_D,\vdash,\dashv)$ will be a general
digroup with $(g_1,g_2)^{-1}_\vdash=(g_2^{-1}g_1^{-1},e)$ and
$(g_1,g_2)^{-1}_\dashv=(e,g_2^{-1}g_1^{-1})$.\\
In \cite[Definition 4.1]{Kin07} Kinyon defines a \emph{Lie digroup} as a general
Lie digroup such that in addition for each $x$ its right inverse 
(w.r.t. to $\vdash$) 
is equal to its left inverse (w.r.t. $\dashv$). This can be shown to be
equivalent to demanding
that the general Lie digroup $(D,e,\vdash,\dashv)$ be balanced.\\ 
Again using the Suschkewitsch decomposition
Theorem (which applies in case the underlying manifold is connected), 
it is not hard to see that the category of all connected Lie digroups
(in the sense of Kinyon) is equivalent to the category of all left
$\mathcal{G}$-spaces, i.e. whose objects are pairs $(G,X)$ where $G$ is 
a connected Lie group and $X$ is a pointed connected left $G$-space (i.e. the distinguished
point of $X$ is a fixed point of the $G$-action) with obvious morphisms.
Recall that the Lie digroup is given by $X\times G$ equipped with the 
point $(e_X,e)$ and the two 
multiplications $(x_1,g_1)\vdash (x_2,g_2)= (g_1x_2,g_1g_2)$ and
$(x_1,g_1)\dashv (x_2,g_2)= (x_1,g_1g_2)$ for all $x_1,x_1\in X$ and
$g_1,g_2\in G$.

The following theorem is a direct consequence of the functorial properties
of the functor $F$:
\begin{theorem}
Let $(D,e,\vdash,\dashv)$ be a bar-unital Lie dimonoid $(D,e,\vdash,\dashv)$.
  \begin{enumerate}
 \item The underlying vector space of the associated coalgebra $F(D)$
 to the bar-unital Lie dimonoid $(D,e,\vdash,\dashv)$
 equipped with $\mathbf{1}$ and the multiplications $\mu_\vdash,\mu_\dashv$
 is an associative bar-unital dialgebra. \\
 In case $D$ is balanced, $F(D)$ is a cocommutative Hopf dialgebra.
 
 \item In case $(D,e,\vdash,\dashv)$ is a Lie digroup (in the sense of Kinyon),
   $F(D)$ is a cocommutative Hopf dialgebra.
 \end{enumerate}
\end{theorem}

\subsubsection{(Lie) racks}

Recall that a \emph{Lie rack} is a pointed manifold with multiplication
$(M,e,\mathbf{m})$ satisfying the following identities for all $x,y,z\in M$
where the standard notation is $\mathbf{m}(x,y)=x\triangleright y$
\begin{eqnarray}
     e \triangleright x & = & x, \label{EqDefRackETriXEqualsX}\\
     x \triangleright e & = & e, \label{EqDefRackXTriEEqualsE} \\
     x \triangleright ( y \triangleright z) & = &
       (x\triangleright y) \triangleright (x\triangleright z)
         \label{EqSelfDistributivity}
\end{eqnarray}
In addition, one demands that $(M,e,\mathbf{m})$ be left-regular, i.e. for all 
$x\in M$ the left multiplication maps
$L_x:y\mapsto x\triangleright y$
should be a diffeomorphism.\\
Note the following version of the self-distributivity identity
(\ref{EqSelfDistributivity}) in terms of maps:
\begin{eqnarray}
\lefteqn{\mathbf{m}\circ (\mathrm{id}_M \times \mathbf{m})}  \nonumber \\
  & =  & \mathbf{m}\circ (\mathbf{m}\times \mathbf{m})\circ 
     (\mathrm{id}_M\times\tau_{M,M}\times \mathrm{id}_M) \circ 
     (\mathrm{diag}_M\times\mathrm{id}_M\times \mathrm{id}_M)
       \label{rack_identity} 
\end{eqnarray}

\begin{exem}
Note that every pointed differentiable manifold $(M,e)$ carries
a \emph{trivial Lie rack structure} defined for all $x,y\in M$ by
\begin{equation}
   x\triangleright_0 y := y,
\end{equation}
and this assignment is functorial.
\end{exem}

\begin{exem}
Any Lie group $G$ becomes a Lie rack upon setting
for all $g,g'\in G$
\begin{equation}
     g\triangleright g' := gg'g^{-1},
\end{equation}
again defining a functor from the category of Lie groups to the 
category of all Lie racks. 
\end{exem}

\begin{exem}
Let $G$ be a Lie group and $V$ be a (smooth) $G$-module 
(supposed to be a real or complex vector space). On $X:=V\times G$, 
we define a binary operation $\rhd$ by
$$(v,g)\rhd(v',g')\,=\,(g(v'),gg'g^{-1})$$
for all $v,v'\in V$ and all $g,g'\in G$. $X$ is a Lie rack with unit 
$1:=(0,1)$ which is called a {\it linear Lie rack}, see \cite{Kin07}.   
\end{exem}

\begin{exem}
Let $(D,e,\vdash,\dashv)$ be a (balanced)
\emph{digroup}. Then formula
(13) of \cite{Kin07},
\[
      x\triangleright y := x\vdash y \dashv x^{-1}
\]
equips the
pointed manifold $(D,e,\triangleright)$ with the structure of a Lie rack.
\end{exem}

Any Lie rack $(M,e,\triangleright)$ can be \emph{gauged} by any
smooth map $f:(M,e)\to (M,e)$ of pointed manifolds satisfying
for all $x,y\in M$
\[
       f(x\triangleright y) = x\triangleright f(y).
\]
A straight-forward computation shows that the pointed manifold $(M,e)$
equipped with the \emph{gauged} multiplication
$\triangleright_f$ defined by
\[
        x\triangleright_f y := f(x) \triangleright y
\]
is a Lie rack $(M,e,\triangleright_f)$. 

Furthermore, recall that an \emph{augmented Lie rack} (see \cite{FR92})
$(M,\phi,G,\ell)$ consists of a pointed differentiable manifold 
$(M,e_M)$,
of a Lie group $G$, of a smooth map $\phi:M\to G$ (of pointed manifolds),
and of a
smooth left $G$-action $\ell:G\times M\to M$ 
(written $(g,x)\mapsto \ell(g,x)=\ell_g(x)=gx$) such that for all
$g\in G$, $x\in M$
\begin{eqnarray}
     ge_M & = & e_M, \label{EqDefAugRackGActionPreservesE} \\
     \phi(gx) & = & g\phi(x)g^{-1}. 
     \label{EqDefAugRackPhiIntertwinesGActionAndConjugation}
\end{eqnarray}
It is a routine check that the multiplication $\triangleright$ on $M$
defined for all $x,y\in M$ by
\begin{equation}
  x\triangleright y  :=  \ell_{\phi(x)}(y)
\end{equation}
satisfies all the axioms (\ref{EqDefRackETriXEqualsX}), 
(\ref{EqDefRackXTriEEqualsE}), and (\ref{EqSelfDistributivity})
of a Lie rack, thus making $(M,e_M,\triangleright)$ into a Lie rack
such that the map $\phi$ is a morphism of Lie racks, i.e.
for all $x,y\in M$
\begin{equation}
     \phi(x\triangleright y)= \phi(x)\phi(y)\phi(x)^{-1}.
\end{equation}
A morphism $(\Psi,\psi):(M,\phi,G,\ell)\to (M',\phi',G',\ell')$ of 
augmented
Lie racks is a pair of maps of pointed differentiable manifolds
$\Psi:M\to M'$ and $\psi:G\to G'$ such that $\psi$ is homomorphism
 of Lie groups and such that all reasonable diagrams commute, viz:
 for all $g\in G$
\begin{eqnarray}
   \phi'\circ \Psi & = & \psi\circ \phi \\
   \Psi\circ \ell_{g} & = & \ell'_{\psi(g)}\circ \Psi
\end{eqnarray}
Note that the trivial Lie rack structure of a pointed manifold $(M,e)$
comes from an augmented Lie rack over the trivial Lie group $G=\{e\}$.

Let $(M,e,\triangleright)$ be a Lie rack.
Applying the functor $F$ we get the following
\begin{theorem}\label{TFOfRackEqualsRackBialgebra}
 The associated coalgebra $F(M)$ with multiplication $\mu$
 of the Lie rack $\big(M,e,\mathbf{m}\big)$ 
 is a rack bialgebra, i.e. satisfying for all $a,b,c\in C$,
 using the same notation $a\triangleright b$ for $\mu(a\otimes b)$:
 \begin{eqnarray}
   \mathbf{1}\triangleright a &=& a,\label{equation1}\\
   a\triangleright \mathbf{1} & = & \epsilon(a)\mathbf{1},\label{equation2}\\
   a \triangleright (b\triangleright c) &=&
     \sum_{(a)} (a^{(1)}\triangleright b)\triangleright
               (a^{(2)}\triangleright c).\label{equation3}
 \end{eqnarray}
\end{theorem}

\begin{prooof}
\begin{enumerate}
\item[(\ref{equation1})] By definition, 
$S\otimes T\in {\mathcal E}_{e}'(M)\otimes {\mathcal E}_{e}'(M)$
are sent by $\mu=\rhd_*$ to the distribution $f\mapsto F_{2~M,M}(S\otimes T)(f\circ\rhd)$. 
We evaluate this formula for $S={\bf 1}$. This gives the distribution 
$f\mapsto {\bf 1}(T^{(2)}(f\circ\rhd))$. But $T^{(2)}$ means that the function is seen as function
of its second variable, i.e. $T^{(2)}(f\circ\rhd)(y)=T(f(y\rhd-))$. On the other hand, 
the delta distribution ${\bf 1}$ evaluates a function in $e$, thus 
$$
  {\bf 1}(T^{(2)}(f\circ\rhd))\,=\,T(f(e\rhd-))\,=\,T(f),
$$
because $e\rhd y=y$ for all $y\in M$. This shows ${\bf 1}\rhd T=T$.  
\item[(\ref{equation2})] Exchanging the roles of the two 
variables in the 
above computation, we obtain for $T\rhd {\bf 1}$ the distribution 
$T({\bf 1}^{(2)}(f\circ\rhd))$
or in other words ${\bf 1}(T^{(1)}(f\circ\rhd))$, i.e. 
the above element $y$ is now in the 
second place. We obtain  
$$
    {\bf 1}(T^{(1)}(f\circ\rhd))\,=\,T(f(-\rhd e))\,=\,T(f(e))
         \,=\,T(1)\,=\,\epsilon(T).
$$
This shows $T\rhd {\bf 1}=\epsilon(T){\bf 1}$. 
\item[(\ref{equation3})] As remarked before, the definition of $\Delta_e$, 
namely $\Delta_e= F_{2~M,M}^{-1}\circ\mathrm{diag}_{M~*}$, induces thanks to the naturality
relation (\ref{naturality_of_F}) relations like  
$$({\rm id}_{M}\times\mathrm{diag}_{M})_*=
F_{2~M,(M\times M)}\circ({\rm id}_{\mathcal{E}'_e(M)}\times F_{2~M,M})
\circ({\rm id}_{\mathcal{E}'_e(M)}\times\Delta_e)\circ F_{2~M,M}^{-1},$$
and 
\begin{eqnarray*}
\lefteqn{
\Delta_e\times{\rm id}_{\mathcal{E}'_e(M)}\times{\rm id}_{\mathcal{E}'_e(M)}}  \\ 
  & = &
\big(F_{2~M,M}^{-1}\times{\rm id}_{\mathcal{E}'_e(M)}
      \times{\rm id}_{\mathcal{E}'_e(M)}\big)
  \circ
      \big(F_{2~(M\times M),M}^{-1}\times{\rm id}_{\mathcal{E}'_e(M)}\big) \\
  & &          \circ F_{2~(M\times M)\times M,M}^{-1}\circ   
                 \big(\mathrm{diag}_{M}\times{\rm id}_M\times{\rm id}_M\big)_*\\
  & &               \circ 
               F_{2~(M\times M),M}
                     \circ \big(F_{2~M,M}\times{\rm id}_{\mathcal{E}'_e(M)}\big).
\end{eqnarray*}
Therefore, starting from the relation induced on $\mathcal{E}'_e(M)$ by relation
(\ref{rack_identity}), one replaces 
$(\mathrm{diag}_{M}\times{\rm id}_M\times{\rm id}_M)_*$ by the above
and obtains finally an equation equivalent to equation 
(\ref{EqDefAlgSelfDistributivity}). 
\end{enumerate}
\end{prooof}

\begin{rem}
This theorem should be compared to Proposition 3.1 in \cite{CCES08}. In \cite{CCES08},
the authors work with the vector space $K[M]$ generated by the rack $M$, while we 
work with point-distributions on a Lie rack $M$. Once again, in some sense, 
we extend their Proposition 3.1
``to all orders''. Observe however that their structure is slightly different (motivated in their 
Remark 7.2). 
\end{rem}  


We get a similar theorem for an augmented Lie rack: Let 
$\mathfrak{g}$ denote the Lie algebra of the Lie group $G$, then we have 
the
\begin{theorem}\label{TFOfAugRackEqualsAugRackBialgebra}
 The associated coalgebra $C$ with multiplication $\mu$
 of an augmented Lie rack $(M,\phi,G,\ell)$ is a cocommutative
 augmented rack bialgebra $(C,\phi_*,\mathsf{U}(\mathfrak{g}),\ell)$
\end{theorem}

We shall close the subsection with a geometric explanation of some of the 
structures appearing in Subsection \ref{SubSecRackBialgebrasETC}: Let
$\big(\mathfrak{h},[~,~]\big)$ be a real finite-dimensional Leibniz algebra.
There is the following Lie rack structure on the manifold $\mathfrak{h}$
defined by
\begin{equation}\label{EqDefStandardLieRackStructureOnLeibniz}
    x\blacktriangleright y := e^{\mathrm{ad}_x}(y)
\end{equation}
Moreover, pick a two-sided ideal $\mathfrak{z}\subset \mathfrak{h}$ with
$Q(\mathfrak{h})\subset \mathfrak{z}\subset \mathfrak{z}(\mathfrak{h})$ so
that the quotient algebra $\mathfrak{g}:=\mathfrak{h}/\mathfrak{z}$ is a Lie 
algebra. Let $p:\mathfrak{h}\to \mathfrak{g}$ be the canoncial projection.
Let $G$ be the connected simply connnected Lie group having Lie algebra
$\mathfrak{g}$. Since $\mathfrak{g}$ acts on $\mathfrak{h}$ as derivations,
there is a unique Lie group action $\ell$ of $G$ on $\mathfrak{h}$ by automorphisms
of Leibniz algebras. Consider the smooth map 
\begin{equation}\label{EqDefPhiForStandardAugmentedLieRack}
    \phi: \mathfrak{h}\to G: x\mapsto \exp\big(p(x)\big).
\end{equation}
Clearly $\phi(g.x)=g\phi(x)g^{-1}$ for all $x\in\mathfrak{h}$ and $g\in G$
whence $(\mathfrak{h},\phi,G,\ell)$ is an augmented Lie rack, and it is not hard
to see that the Lie rack structure coincides with 
(\ref{EqDefStandardLieRackStructureOnLeibniz}).
\begin{theorem}
 The $C^5$-rack bialgebra associated to the augmented Lie rack
 $(\mathfrak{h},\phi,G,\ell)$ by means of the Serre functor
 is isomorphic to the universal envelopping algebra of infinite order,
 $\mathsf{UAR}^\infty(\mathfrak{h})$, see Definition 
 \ref{DefUniversalRackBialgebraInfiniteOrder} and Theorem
 \ref{TLeibnizToUARInfinityOfLeibniz}.
\end{theorem}
\begin{prooof}
 First we compute $\phi_*= \exp_*\circ p_*$. Since 
 $p:\mathfrak{h}\to\mathfrak{g}$ is linear, it is easy
 to see using formula (\ref{EqCompPhiOne}) that for all $k\in \mathbb{N}$ and
 $x_1,\ldots,x_k\in\mathfrak{h}$
 \[
     p_*\big(\Phi_\mathsf{S}(x_1\bullet\cdots\bullet x_k)\big)
       = \Phi_\mathsf{S}\big(p(x_1)\bullet\cdots\bullet p(x_k)\big)
        =\Phi_\mathsf{S}\Big(\mathsf{S}(p)
        \big(x_1\bullet\cdots\bullet x_k\big)\Big),
 \]
 see (\ref{EqCompPhiOne}) for a definition of $\Phi_\mathsf{S}$.
 Next, for all $k\in\mathbb{N}$ and $\xi_1,\ldots,\xi_k\in\mathfrak{g}$,
 we shall show the formula (for all $f'\in\mathcal{C}^\infty(G,\mathbb{K})$)
 \begin{eqnarray*}
    \Big(\exp_*\big(\Phi_\mathsf{S}(\xi_1\bullet\cdots\bullet \xi_k)\big) \Big)(f')
    & = & \sum_{\sigma\in S_k} \frac{1}{k!}
         \big(\Phi_\mathsf{U}(\xi_{\sigma(1)}\cdots\xi_{\sigma(k)})\big)(f')\\
     & = &  \Big(\Phi_\mathsf{U}\big(\omega
         (\xi_1\bullet\cdots\bullet\xi_k)\big)\Big)(f')
 \end{eqnarray*}
 (see eqn (\ref{EqDefIsoUniversalEnvOfGToFOfG}) for a definition of 
 $\Phi_\mathsf{U}$). Both sides of this equation are symmetric $k$-linear maps
 in the arguments $\xi_1,\ldots,\xi_k$,
 hence by the Polarization Lemma (see e.g.~\cite{Tho2013}), it suffices to
 check equality in case $\xi_1=\cdots =\xi_k=\xi$. Since for each real number
 $t$ the map $g\mapsto F^\xi_t(g):=g\exp(t\xi)$ is the flow of the left invariant vector
 field $\xi^+$, we get
 \begin{eqnarray*}
  \Big(\exp_*\big(\Phi_\mathsf{S}(\xi\bullet\cdots\bullet \xi)\big) \Big)(f')
    & = & \left.\frac{\partial^k}{\partial t^k}
            \Big(f'\big(\exp(t\xi)\big)\Big)\right|_{t=0} 
     =  \left.\frac{\partial^k}{\partial t^k}
            \Big(f'\big(F^\xi_t(e))\big)\Big)\right|_{t=0} \\
    & = & \left.\frac{\partial^k}{\partial t^k}
            \Big(\big((F^\xi_t)^*f'\big)(e))\big)\Big)\right|_{t=0} \\
    & = & \left.\Big((F^\xi_t)^*
             \big(\big(\mathbb{L}_{\xi^+}\circ\cdots\circ\mathbb{L}_{\xi^+}\big)
                (f')\big)\Big)(e)\right|_{t=0} \\
    & = &  \big(\big(\mathbb{L}_{\xi^+}\circ\cdots\circ\mathbb{L}_{\xi^+}\big)
                (f')\big)(e)\\
    & = &  \Big(\Phi_\mathsf{U}\big(\omega
         (\xi\bullet\cdots\bullet\xi)\big)\Big)(f')
 \end{eqnarray*}
 proving the above formula. It follows that
 \begin{equation}\label{EqCompExpStarCircPStar}
   \phi_* = \exp_*\circ p_* = \Phi_{\mathsf{U}} \circ \omega \circ 
                                  \mathsf{S}(p) \circ \Phi_\mathsf{S}^{-1}.
 \end{equation}
 Next, we compute $\ell_*$. We get for positive integers $k,l$, 
 $\xi_1,\ldots,\xi_k\in\mathfrak{g}$, $x\in \mathfrak{h}$, and
 $f\in\mathcal{C}^\infty(\mathfrak{h},\mathbb{K})$:
 \begin{eqnarray*}
   \lefteqn{ \Big(\ell_*\big(\Phi_\mathsf{U}(\xi_1\cdots\xi_k)\otimes 
    \Phi_\mathsf{S}(x^{\bullet l})\big)\Big)(f)}\\
    & = & 
    \left.\frac{\partial^{k+l}}{\partial s_1\cdots \partial s_k \partial t^l}
    \Big(f\big(
      \big(\ell_{\exp(s_1 \xi_1)}\circ \cdots\circ \ell_{\exp(s_k \xi_k)}\big)
         (tx)
        \big)\Big)\right|_{s_1=\cdots=s_k=0=t}\\
    & = &  
    \left.\frac{\partial^{k}}{\partial s_1\cdots \partial s_k}
     \sum_{j_1\cdots j_l=1}^{\dim (\mathfrak{h})}
    \left( \frac{\partial^l f}{\partial x_{j_1}\cdots \partial  x_{j_l}}(0)
      \ell_s(x)_{j_1}\cdots \ell_s(x)_{j_l}\right)\right|_{s_1=\cdots=s_k=0}     
 \end{eqnarray*}
where in the last line we have used a basis of 
$\mathfrak{h}$, have written $y_1,\ldots,y_n$ ($n=\dim(\mathfrak{h})$)
for the components of each vector $y\in\mathfrak{h}$, and
used the notation $\ell_s$ for the linear map 
$\ell_{\exp(s_1 \xi_1)}\circ \cdots\circ \ell_{\exp(s_k \xi_k)}$.
By induction on $k$ it is easy to prove that
\[
   \left.\frac{\partial^{k}}{\partial s_1\cdots \partial s_k}\Big(
      \ell_s(x)\bullet\cdots \bullet \ell_s(x)\Big)
      \right|_{s_1=\cdots=s_k=0} 
        = \mathrm{ad}^s_{\xi_1\cdots\xi_k}\big(x^{\bullet l}\big),
\]
and using again the Polarisation Lemma, we finally get for all 
$u\in \mathsf{U}(\mathfrak{g})$ and $\alpha \in \mathsf{S}(\mathfrak{h})$
\begin{equation}
    \ell_*\big(\Phi_\mathsf{U}(u)\otimes \Phi_\mathsf{S}(\alpha)\big)
       = \Phi_\mathsf{S}\big(\mathrm{ad}^s_u(\alpha)\big)=
         \Phi_\mathsf{S}\big(u.\alpha\big),
\end{equation}
and the isomorphism with the augmented rack bialgebra
$\mathsf{UAR}^\infty(\mathfrak{h})=\mathsf{S}(\mathsf{h})$ is established.
\end{prooof}

\begin{rem}
Observe that the Serre functor can be rendered completely algebraic, i.e. 
for example for an algebraic Lie rack $R$ (meaning that the underlying pointed 
manifold is a smooth algebraic variety and the rack product is algebraic),
one can take as its Serre functor image $F(R)$ the space of derivations
along the evaluation map in the distinguished point. The composition of $F$ 
with the functor of primitives gives then the tangent functor (see text 
before Remark \ref{remark}). 
This gives a new and
completely algebraic functorial way to associate to a Lie rack its tangent 
Leibniz algebra.
\end{rem}  

\appendix
 \section{Some definitions around coalgebras}
   \label{AppCoalgebras}
 
Let $C$ be a module over a commutative associative
unital ring $K$ (which we shall assume to contain $\mathbb{Q}$). Recall that
a linear map $\Delta:C\to C\otimes_K C=C\otimes C$ is called a 
\emph{coassociative comultiplication} iff 
$\big(\Delta\otimes \mathrm{id}_C\big)\circ \Delta=
\big( \mathrm{id}_C\otimes\Delta\big)\circ \Delta$, and the pair $(C,\Delta)$
is called a (coassociative) \emph{coalgebra over $K$}. Let
$(C',\Delta')$ be another coalgebra. Recall that a $K$-linear map 
$\Phi:C\to C'$ is called a \emph{homomorphism of coalgebras} iff
$\Delta'\circ\phi=(\phi\otimes \phi)\circ \Delta$. The coalgebra 
$(C,\Delta)$ is called \emph{cocommutative} iff $\tau\circ \Delta=\Delta$
where $\tau:C\otimes C\to C\otimes C$ denotes the canonical flip map.
Recall furthermore that a
linear map $\epsilon:C\to K$ is called a \emph{counit} for the coalgebra 
$(C,\Delta)$ iff $\big(\epsilon\otimes \mathrm{id}_C\big)\circ \Delta=
\big( \mathrm{id}_C\otimes\epsilon\big)\circ \Delta=\mathrm{id}_C$. 
The triple $(C,\Delta,\epsilon)$ is called a \emph{counital coalgebra}.
Moreover, a counital coalgebra $(C,\Delta,\epsilon)$ equipped with an
element $\mathbf{1}$ is called \emph{coaugmented} iff 
$\Delta(\mathbf{1})=\mathbf{1}\otimes \mathbf{1}$ and 
$\epsilon(\mathbf{1})=1\in K$. Let $C^+\subset C$ denote the kernel of $\epsilon$.
Recall that a morphism 
$\phi:(C,\Delta,\epsilon,\mathbf{1})\to (C',\Delta',\epsilon',\mathbf{1}')$
of counital coaugmented coalgebras over $K$ is a $K$-linear map satifying
$(\phi\otimes \phi)\circ \Delta=\Delta'\circ \phi$, $\epsilon'\circ\phi=\epsilon$,
and $\phi(\mathbf{1})=\mathbf{1}'$.
Moreover, for any counital coaugmented coalgebra
the $K$-submodule of all \emph{primitive elements} is defined by
\begin{equation}\label{EqDefPrimitives}
 \mathsf{Prim}(C):=
     \{x\in C~|~\Delta(x)=x\otimes \mathbf{1}+\mathbf{1}\otimes x\}.
\end{equation}
Every morphism of counital coaugmented coalgebra clearly maps primitive elements
to primitive elements, thus defining a functor $\mathsf{Prim}$ from the category
of counital coaugmented coalgebras to the category of $K$-modules.
Finally, following Quillen \cite{Qui69},
we shall call a counital coaugmented coalgebra \emph{connected} iff the 
following holds: The
sequence of submodules $(C_{(r)})_{r\in\mathbb{N}}$ defined by
$C_{(0)}=K\mathbf{1}$ and recursively by
\begin{equation}\label{EqDefSubcoalgebraOfOrderK}
   C_{(k+1)} :=\{x\in C~|~ \Delta(x)-x\otimes\mathbf{1} - 
                             \mathbf{1}\otimes x~\in~ C_{(k)} \}
\end{equation}
is easily seen to be an ascending sequence of coaugmented counital 
subcoalgebras of $(C,\Delta,\epsilon,\mathbf{1})$, and if the union of all
the $C_{(k)}$ is equal to $C$, then  $(C,\Delta,\epsilon,\mathbf{1})$ is called
connected. We refer to each $C_{(k)}$ as the \emph{subcoalgebra of order $k$}.
Clearly, each  $C_{(k)}$ is connected, and 
$C_{(1)}=K\mathbf{1}\,\oplus\,\mathsf{Prim}(C)$. Moreover, each morphism of counital
coaugmented coalgebras maps each subcoalgebra of order $k$ to the subcoalgebra
of order $k$ of the target coalgebra thus defining a functor $C\to C_{(k)}$
from the category of coaugmented counital coalgebras to itself.
We shall use the following acronyms:
\begin{defi}
We call a coassociative, counital, coaugmented coalgebra a $C^3$-coalgebra.
In case the $C^3$-coalgebra is in addition cocommutative, we shall speak of a
$C^4$-coalgebra. Finally, a connected $C^4$-coalgebra will be coined 
a $C^5$-coalgebra.
\end{defi}

\noindent Recall also that the tensor product of two counital coaugmented
coalgebras $(C,\Delta,\epsilon,\mathbf{1})$ and
$(C',\Delta',\epsilon',\mathbf{1}')$ is given by
$(C\otimes C',(\mathrm{id}_C\otimes \tau\otimes \mathrm{id}_{C'})\circ 
(\Delta\otimes \Delta'),\epsilon\otimes\epsilon',\mathbf{1}\otimes\mathbf{1}')$.
Tensor products of connected coalgebras are connected.

 Recall the standard example: Let $V$ be a $K$-module and
$\mathsf{S}(V)=\oplus_{r=0}^\infty\mathsf{S}^r(V)$ be the 
\emph{symmetric algebra generated by $V$}, i.e. the free algebra 
$\mathsf{T}(V)$ (for which we denote the tensor multiplication by suppressing
the symbol)
modulo the two-sided ideal $I$ generated by $xy-yx$ for all 
$x,y\in V$.
Denoting the commutative associative multiplication in $\mathsf{S}(V)$ 
(which is 
induced by the free multiplication)
by $\bullet$, i.e.
\[
    x_1\bullet\cdots\bullet x_k :=
    x_1\cdots x_k ~\mathrm{mod}~I,
\]
we have $\Delta(x)=x\otimes \mathbf{1}+\mathbf{1}\otimes x$ for all
$x\in V$ and $\Delta(x_1\bullet\cdots\bullet x_k)=
(x_1\otimes \mathbf{1}+\mathbf{1}\otimes x_1)\bullet \cdots\bullet
(x_k\otimes \mathbf{1}+\mathbf{1}\otimes x_k)$ for all positive integers $k$ and
$x_1,\ldots,x_k\in V$. Recall that $S^0(V)$ is the free $K$-module 
$K\mathbf{1}$ and the counit is defined by $\epsilon(\lambda \mathbf{1})
=\lambda$
for all $\lambda\in K$ and by declaring that $\epsilon$ vanishes on
$\oplus_{r=1}^\infty \mathsf{S}^r(V)$. Moreover, the submodules
$\big(\mathsf{S}(V)\big)_{(n)}$ are given by 
$\oplus_{r=0}^n \mathsf{S}^r(V)$, whence $\mathsf{S}(V)$ is clearly connected,
so it is a $C^5$-coalgebra whose submodule of primitive elements equals $V$.

Moreover, for a given coalgebra $(C,\Delta)$ and a given nonassociative
algebra $(A,\mu)$ where $\mu:A\otimes A\to A$ is a given $K$-linear map, recall
the \emph{convolution multiplication}  in the $K$-module $\mathrm{Hom}_K(C,A)$
defined in the usual way for any two $K$-linear maps $\phi,\psi:c\to A$ by
\begin{equation}\label{EqDefConvolution}
    \phi *\psi := \mu\circ (\phi\otimes \psi)\circ \Delta.
\end{equation}
In case $\Delta$ is coassociative and $\mu$ associative, $*$ will be 
associative. The following fact is rather important: If $C$ is connected
and if the $K$-linear map $\varphi:C\to A$ vanishes on $\mathbf{1}_C$,
then \emph{any convolution power series of $\varphi$ converges}, i.e.~
the evaluation of some formal series $\sum_{r=0}^\infty \lambda_r \varphi^{*r}$ 
(with
$\lambda_r\in K$ and $\varphi^{*0}:=\mathbf{1}_A\epsilon_C$) 
on $c\in C$ always reduces to a finite number of terms. In particular, let
$\psi:C\to A$ be a $K$-linear map such that $\psi(\mathbf{1}_C)=\mathbf{1}_A$.
Then --as has been observed by Takeuchi and Sweedler (see \cite[Lemma 14]{Tak71}
or \cite[Lemma 9.2.3]{Swe69}-- $\psi$ has always a \emph{convolution inverse},
i.e.~there is a unique $K$-linear map $\psi':C\to A$ such that
$\psi*\psi'=\mathbf{1}_A\epsilon_C=\psi'*\psi$, where $\psi'$ is defined by the 
geometric series $\psi'= \sum_{r=0}^\infty (\mathbf{1}_A\epsilon_C-\psi)^{*r}$.

\section{Semigroups}
   \label{AppSemigroups}

We collect some properties of semigroups which are very old, but
a bit less well-known than properties of groups. The standard reference to these
topics is the book
\cite{CP61}
by A. H. Clifford and G. B. Preston.\\
Recall that a \emph{semigroup} $\Gamma$ is a set equipped with an associative
multiplication $\Gamma\times \Gamma\to \Gamma$, written $(x,y)\mapsto xy$.
An element $e$ of $\Gamma$ is called a \emph{left unit element} 
(resp. a \emph{right unit element} resp. a \emph{unit element}) iff for all
$x\in\Gamma$ we have $ex=x$ (resp. $xe=x$ resp. iff $e$ is both left and right
unit element). A pair $(\Gamma,e)$ of a semigroup $\Gamma$ and an element
$e$ is called \emph{left unital} (resp. \emph{right unital} resp. \emph{unital})
iff $e$ is a left unit element (resp. a right unit element resp. a unit element).
A unital semigroup is also called a \emph{monoid}. It is well-known that
the unit element of a monoid is the unique unit element (unlike left or right unit 
elements in general). Let $(\Gamma,e)$
be a right unital or a left unital semigroup. Recall that for a given
element $x\in\Gamma$ an element $y\in\Gamma$ is called a 
\emph{left inverse of $x$} (resp. a \emph{right inverse of $x$} resp. an
\emph{inverse of $x$}) iff $yx=e$ (resp. $xy=e$ resp. iff $y$ is both a left and
a right inverse of $x$). Clearly, a unital semigroup $(\Gamma,e)$ such that
every element has an inverse is a \emph{group}. In that case it is well-known that
for each $x$ there is exactly one inverse element, called 
$x^{-1}$.\\
Note that by a Lemma by \emph{L. E. Dickson} (1905, see \cite[p.4]{CP61} for the 
reference)
every \emph{left} unital semigroup such that each element has at least one 
\emph{left} inverse
is already a group which can be shown by just using the definitions. Dually,
every \emph{right} unital semigroup such that each element has at least one 
\emph{right} inverse is also a group. \\
More interesting is the case of a 
\emph{left} (resp. right)
unital semigroup $(\Gamma,e)$ such that every element $x$ has at least one
\emph{right} (resp. left) inverse element. In that case (which is an equivalent
formulation of a so-called \emph{right group} (resp. left group), 
see \cite[p.37]{CP61}), the conclusion of Dickson's Lemma does no longer hold.
In order to see what is going on, there is first the following useful
\begin{lemma}\label{LRightInversesHavingRightInverses}
   Let $(\Gamma,e)$ be a left-unital semigroup, let $a,b,c$ three elements of 
   $\Gamma$ such that
   \[
        ab=e~~\mathrm{and}~~bc=e.
   \]
   Then
   \[
        c=ae,~be=b,
   \]
   and the left multiplications $L_a:x\mapsto ax$ and $L_b:x\mapsto bx$ are 
   invertible. In particular, given the element $a$, its right inverse $b$ is 
   unique under the above hypotheses.
\end{lemma}
The proof is straight-forward.\\
The structure of right (resp.left) groups is completely settled in the 
\emph{Suschkewitsch Decomposition Theorem}, 1928: Given a right group $(\Gamma, e)$,
it can be shown --using the above Lemma and elementary manipulations, see
also \cite[p.38, Thm 1.27]{CP61}--
that all the left multiplications $L_x:y\mapsto xy$
(resp. right multiplications $R_x:y\mapsto yx$) are invertible, that for each
element there is exactly one right (resp. left) inverse (whence there is a map 
$\Gamma\to \Gamma$
assigning to each element $x$ its right (resp. left) inverse $x^{-1}$), 
that the image of this right (resp. left)
inverse map is equal to $\Gamma e$ (resp. $e\Gamma)$ 
(which turns out to be a \emph{subgroup}
of $(\Gamma,e)$), and that $(\Gamma,e)$ is isomorphic to the cartesian product
$\big(\Gamma e \times E,(e,e)\big)$ (resp. $\big(E\times e\Gamma,(e,e)\big)$ 
where $E$ is the 
set of all left 
(resp. right) unit elements in $(\Gamma,e)$ (coinciding with the set of all 
idempotent elements). For right groups, the aforementioned isomorphism
is given as follows:
\begin{eqnarray}
  \phi:\Gamma e\times E\to \Gamma & : & (a,f)\mapsto af, \\ 
  \phi^{-1}: \Gamma \to \Gamma e\times E & :&  x\mapsto (xe,x^{-1}x).
\end{eqnarray}
Note that both components of $\phi^{-1}$ are idempotent maps.
There is a completely analogous statement for left groups.

Recall that a \emph{Lie semigroup} is a differentiable manifold $\Gamma$ 
equipped
with a smooth associative multiplication 
$\mathbf{m}:\Gamma\times \Gamma \to \Gamma$.
All the other definitions of semigroups mentioned above
(such as left unital, right unital semigroups, monoids, groups, right groups, 
left groups etc.) carry over to the Lie, i.e. differentiable, case. \\
Moreover for 
\emph{right Lie groups}, 
it is easy to see that all the left multiplications are 
diffeomorphisms (since their inverse maps are left multiplications with the 
inverse elements and therefore 
smooth). This fact and the regular value theorem applied to the equation
$xy=e$ imply that the right inverse
map is smooth since its graph is a closed submanifold of 
$\Gamma\times\Gamma$ and 
the restriction of the projection on the first factor of the graph is a
diffeomorphism.
As the maps $x\mapsto xe=(x^{-1})^{-1}$ and $x\mapsto x^{-1}x$ are smooth and
idempotent, it follows that their images, the subgroup $\Gamma e$, and the 
semigroup of all left
unit elements, $E$, are both smooth submanifolds of $\Gamma$ and closed sets
provided $\Gamma$ is connected, see e.g. \cite[p.54, Satz 5.13]{BJ73} for 
a proof. Hence $\Gamma e$ is a connected Lie group, and the 
Suschkewitsch decomposition $\Gamma\cong \Gamma e \times E$, see 
Appendix \ref{AppSemigroups}, is a diffeomorphism. Conversely, any cartesian
product of a Lie group $G$ and a differentiable manifold $E$ equipped with the 
multiplication $(g,x)(h,y):= (gh,y)$ is easily seen to be a right Lie group.
An analogous statement holds for left Lie groups.\\
It is not hard to see that the category of all connected right Lie groups is equivalent
to the product of category of all connected Lie groups and the category of all pointed
connected manifolds.

\end{document}